\documentclass[a4paper,11pt]{article}
\usepackage[top=3cm,bottom=3cm,left=3cm,right=2.5cm]{geometry}
\usepackage{amsfonts}\usepackage{amssymb}
\usepackage{graphicx}\usepackage{amsmath}
\usepackage{subfigure}
\usepackage{abstract}
\usepackage{color}
\usepackage{cite}
\usepackage{graphicx}
\usepackage{amsmath}
\usepackage{amsthm}
\usepackage{float}

\numberwithin{equation}{section}
\providecommand{\keywords}[1]{\textbf{Key words.} #1}

\title{A posteriori error estimates of mixed discontinuous Galerkin method for the Stokes eigenvalue problem}

\author{Lingling Sun$^{1,2}$,  Hai Bi$^{1}$, Yidu Yang$^{1,}$\footnote{Corresponding author. Email address: ydyang@gznu.edu.cn}
 \\\\
{\small $^1$ School of Mathematical Sciences,
Guizhou Normal University,Guiyang,  $550001$,  China}\\{\small $^2$ School of Biology and Engineering, Guizhou Medical University,}{\small Guiyang, 550001, China}
}

\begin{document}
\date{}
\maketitle
\begin{abstract}
\indent In this paper, for the Stokes eigenvalue problem in $d$-dimensional case $(d=2,3)$, we present an a posteriori error estimate of residual type of the mixed discontinuous Galerkin finite element method
using $\mathbb{P}_{k}-\mathbb{P}_{k-1}$ element $(k\geq 1)$.
We give the a posteriori error estimators for approximate eigenpairs, prove their reliability
and efficiency for eigenfunctions, and also analyze their reliability
for eigenvalues. We implement adaptive calculation, and
the numerical results confirm our theoretical predictions and show that our method
can achieve the optimal convergence order $O(dof^{\frac{-2k}{d}})$.
\end{abstract}

\section{Introduction}
\indent~~~~~Stokes eigenvalue problem is of great importance because of their role
for the stability analysis in fluid mechanics. Hence, the development of efficient numerical
methods for the problem is of great interest.\\
\indent Adaptive finite element methods are favored in current science and engineering computing.
For a given tolerance, adaptive finite element methods require little degrees of freedom. So far, many
excellent works on the a posteriori error estimates and adaptive algorithm have been summarized  in
previous studies (see \cite{babuska1978,Verfurth1996,dorfler1996,morin2002,oden2011,brenner2012,Shi2013,sunzhou2016}, etc).
 For the Stokes eigenvalue problem, the a posteriori error estimates  has received much attention.
 For example, \cite{Lovadina2009, Han2015, Armentano2014} studied
 the a posteriori error estimates of conforming mixed method,
 Liu et al. \cite{Liu2013} presented some super-convergence
results and the related recovery type a posteriori error estimators for conforming mixed method,
  Jia et al. \cite{Xie2014}
 discussed the a posteriori error estimate of low-order non-conforming finite element,
 Gedicke et al. \cite{Gedicke2018} conducted the a
posteriori error analysis for the Arnold-Winther mixed finite element method using the stress-velocity formulation,  $\ddot{O}$nder T$\dot{u}$rk
et al. \cite{onder2016} researched a stabilized finite element method for the two-field (displacement-pressure) and three-field (stress-displacement-pressure) formulations of the Stokes eigenvalue problem.\\
\indent Discontinuous Galerkin finite element method (DGFEM) was first introduced by Reed and Hill \cite{reed1973} in 1973 and has been developed greatly (see, e.g.,  \cite{Cockburn1999,Hesthaven2008,Riviere2008,Pietro2009,Cangianl2010,Antonio2012,Brezzi2000,Ern2021}).
DGFEM for eigenvalue problems has also been discussed in many papers (see \cite{Antonietti2006, Zeng2017, Brenner2015, Wang2019, Geng2016, Yang2017,  Buffa2006, Buffa2007, Gedicke2019, LepeF2020}).
Among them Gedicke et al. \cite{ Gedicke2019} discussed the a posteriori error estimate for the divergence-conforming DGFEM using Raviart-Thomas element for velocity-pressure formulation of the Stokes eigenvalue
problem on shape-regular rectangular meshes.
Felipe Lepe et al. \cite{LepeF2020} analyzed symmetric and nonsymmetric discontinuous Galerkin methods for a pseudostress formulation of the Stokes eigenvalue problem.  It can be seen that solving the Stokes eigenvalue problem by DGFEM has attracted extensive attention of scholars . \\
\indent For the Stokes equations, it has been studied the mixed DGFEM using $\mathbb{P}_{k}-\mathbb{P}_{k-1}$ element (see \cite{Riviere2008,Badia2014,Hansbo2002,Girault2005,Riviere2006}) and $\mathbb{Q}_{k}-\mathbb{Q}_{k-1}$ element (see \cite{toselli2002,Schotzau2002,Houston2005}), which laid a foundation for us to further study the Stokes eigenvalue problem.
Based on the above work, in this paper,
we study the residual type a posteriori error estimate of the mixed DGFEM
using $\mathbb{P}_{k}-\mathbb{P}_{k-1}$ element $(k=1, 2, 3,\cdots)$
for velocity-pressure formulation of the Stokes eigenvalue problem on shape-regular simplex meshes in $\mathbb{R}^{d}(d=2,3)$.
In Section 2, we give the a prior error estimates for the mixed DGFEM of the Stokes eigenvalue problem based on the work of \cite{Badia2014}.
 In Section 3, we give the a posterior error estimators for approximate eigenpair, and use the method of \cite{Verfurth1996,oden2011} together with the enriching operator in \cite{Ohannes2003,Brenner2003} and the lifting operator in \cite{Perugia2002,Schotzau2002} to prove the reliability and efficiency of the estimator for eigenfunctions.
In Section 4, we implement adaptive calculation. The numerical results show that the approximate eigenvalues obtained by our method have the same accuracy as those in \cite{Gedicke2018, Gedicke2019} and achieve the optimal convergence order $O(dof^{\frac{-2k}{d}})$, which validates that our method is effective.\\
\indent The characteristic of the DGFEM discussed in this paper is that for the Stokes eigenvalue problem both in two and three-dimensional domains, it can use high-order elements so that it can not only capture smooth solutions but also achieve the optimal convergence order for local low smooth solutions (eigenfunctions have local singularity or local low smoothness) on adaptive locally refined graded meshes.\\
\indent Throughout this paper, $C$ denotes   a generic positive constant independent of the mesh size $h$, which may not be the same at each occurrence. We use the symbol $a\lesssim b$ to mean that $a\leq C b$, and $a\backsimeq b$ to mean that $a\lesssim b$ and $b\lesssim a$.

\section{Preliminary}
\indent Consider the following Stokes eigenvalue problem:
\begin{equation}\label{s2.1}
\begin{cases}
-\mu\Delta \mathbf{u} +\nabla p=\lambda \mathbf{u},~~~~~ in~~\Omega,\\
div \mathbf{u}=0,~~~~~~~~~~~~~~~~~in ~~\Omega,\\
\mathbf{u}=0,~~~~~~~~~~~~~~~~~~~~~on~~\partial\Omega,
\end{cases}
\end{equation}
where $\Omega\subset \mathbb{R}^{d}(d=2,3)$ is a bounded polyhedral domain,
$\mathbf{u}=(u_{1},...,u_{d})$ is the
velocity of the flow, $p$ is the pressure and $\mu>0$
is the kinematic viscosity parameter of the fluid.\\
\indent Note that for constant viscosity $\mu$, the
velocity eigenfunctions do not change in $\mu$ and thus the eigenvalues as well as the
pressure eigenfunctions scale linearly in $\mu$, i.e., the eigenpair for arbitrary constant
$\mu$ is $(\mu\lambda, \mathbf{u},\mu p)$, where  $(\lambda,\mathbf{u},p)$ denotes the eigenpair for $\mu=1$.
Hence, in this paper we only consider the case of $\mu=1$.\\
\indent  In this paper, denote by $H^{\rho}(\Omega)$ the Sobolev space on $\Omega$ of order $\rho\geq 0$ equipped with the norm $\|\cdot\|_{\rho,\Omega}$ (denoted by $\|\cdot\|_{\rho}$ for simplicity). $H_{0}^{1}(\Omega)=\{v\in H^{1}(\Omega), v|_{\partial\Omega}=0\}$.
For $\mathbf{v}=(v_{1},\cdots,v_{d})\in H^{\rho}(\Omega)^{d}$, denote $  \|\mathbf{v}\|_{\rho}=\sum\limits_{i=1}^{d}\|v_{i}\|_{\rho}$.
We also use the notation $(\cdot,\cdot)$ to denote the inner product in $L^{2}(\Omega)^{d}$ which is given by $(u,v)=\int_{\Omega}uv dx$ for $d=1$ and $(\mathbf{u},\mathbf{v})=\int_{\Omega}\mathbf{u}\cdot\mathbf{v} dx$ for $d=2,3$.
Define $\mathbf{V}=H_{0}^{1}(\Omega)^{d}$ with the norm $\|\mathbf{v}\|_{\mathbf{V}}=(\nabla\mathbf{v},\nabla\mathbf{v})^{\frac{1}{2}}$, and define $Q = L_{0}^{2}(\Omega)=\{q\in L^{2}(\Omega):(q,1)=0\}$.\\
\indent The weak formulation of (\ref{s2.1}) is given by:
 find $(\lambda,\mathbf{u},p)\in\mathbb{R}\times \mathbf{V} \times Q$, $\|\mathbf{u}\|_{0}=1$, such that
\begin{align}
\mathcal{A}(\mathbf{u},\mathbf{v})+\mathcal{B}(\mathbf{v},p)&=\lambda (\mathbf{u},\mathbf{v}),~~~\forall \mathbf{v}\in \mathbf{V}\label{s2.2}\\
\mathcal{B}(\mathbf{u},q)&=0, ~~~~\forall q\in Q,\label{s2.3}
\end{align}
where
\begin{align*}
&\mathcal{A}(\mathbf{u},\mathbf{v})=(\nabla\mathbf{u}, \nabla\mathbf{v}),\\
&\mathcal{B}(\mathbf{v},q)=-(div\mathbf{v},q).
\end{align*}
\indent The existence and uniqueness of the velocity $\mathbf{u}$ follow from the Lax-Milgram lemma in the space
$\mathbf{Z}=\{\mathbf{v\in\mathbf{V}}: b(\mathbf{v},q)=0, \forall q\in Q\}$. The stability of the pressure can be obtained by the well-known inf-sup
condition (see \cite{girault1979}):
\begin{align*}
\beta\|q\|_{L^{2}(\Omega)}\leq \sup\limits_{\mathbf{v}\in \mathbf{V}, \mathbf{v}\neq 0}\frac{\mathcal{B}(\mathbf{v},q)}{\|\mathbf{v}\|_{\mathbf{V}}},  \forall q\in Q.
\end{align*}

\indent Let $\pi_{h}=\{\kappa\}$ be a regular partition of $\Omega$ with the mesh diameter $h=\max\limits_{\kappa\in \pi_{h}} h_{\kappa}$
where $h_{\kappa}$ is the diameter of element $\kappa$.
Let $\varepsilon_{h}=\varepsilon_{h}^{i}\cup\varepsilon_{h}^{b}$ where $\varepsilon_{h}^{i}$ denotes the interior faces (edges) set and $\varepsilon_{h}^{b}$ denotes the set of faces (edges) lying on the boundary $\partial\Omega$.
 We denote by $|\kappa|$  and $|E|$ the measure of $\kappa$ and $E\in \varepsilon_{h}$, respectively.
 Let $(\cdot, \cdot)_{\kappa}$ and $(\cdot, \cdot)_{E}$ denote the inner product in $L^{2}(\kappa)$ and $L^{2}(E)$, respectively.
  We denote by $\omega(\kappa)$ the union of all elements  having at least one face (edge) in common with $\kappa$,
 and denote by $\omega(E)$ the union of the elements having in common with $E$.\\
\indent Define a broken Sobolev space
$$H^{1}(\Omega, \pi_{h})=\{v\in L^{2}(\Omega): v|_{\kappa}\in H^{1}(\kappa),~\forall \kappa\in\pi_{h}\}.$$
 For any $E\in\varepsilon_{h}^{i}$, there are two simplices $\kappa^{+}$ and $\kappa^{-}$ such that $E=\kappa^{+}\cap\kappa^{-}$. Let
$\mathbf{n}^{+}$ be the unit normal of $E$ pointing from $\kappa^{+}$ to $\kappa^{-}$ and let $\mathbf{n}^{-}=-\mathbf{n}^{+}$. \\
\indent For any $v\in H^{1}(\Omega, \pi_{h})$, we define the jump and mean of $\mathbf{v}$ on $E$ by
\begin{align*}
[\![v]\!]=v^{+}\mathbf{n}^{+}+v^{-}\mathbf{n}^{-}, ~~~\mathbf{\{}v\mathbf{\}}=\frac{1}{2}(v^{+}+v^{-}),
\end{align*}
where $v^{\pm}=v|_{\kappa^{\pm}}$.\\
For $\mathbf{v}\in H^{1}(\Omega, \pi_{h})^{d}$, we define the jump and mean of $\mathbf{v}$ on $E\in\varepsilon_{h}^{i}$ by
\begin{align*}
[\![\mathbf{v}]\!]=\mathbf{v}^{+}\cdot\mathbf{n}^{+}+\mathbf{v}^{-}\cdot\mathbf{n}^{-},~~~\{\mathbf{v}\}=\frac{1}{2}(\mathbf{v}^{+}+\mathbf{v}^{-}).
\end{align*}
We also require the full jump of vector-valued functions. For $\mathbf{v}\in H^{1}(\Omega, \pi_{h})^{d}$ , we define the full
jump by
\begin{align*}
[\![\underline{\mathbf{v}}]\!]=\mathbf{v}^{+}\otimes \mathbf{n}^{+}+\mathbf{v}^{-}\otimes \mathbf{n}^{-},
\end{align*}
where for two vectors in Cartesian coordinates $\mathbf{a}=(a_{i})$ and $\mathbf{b}=(b_{j})$, we define the matrix
$\mathbf{a}\otimes \mathbf{b}=[a_{i}b_{j}]_{1\leq i,j\leq d}$. Similarly, for tensors $\tau\in H^{1}(\Omega, \pi_{h})^{d\times d}$, the jump and mean on $E\in\varepsilon_{h}^{i}$ are
defined as follows, respectively:
\begin{align*}
[\![\mathbf{\tau}]\!]=\tau^{+}\mathbf{n}^{+}+\tau^{-}\mathbf{n}^{-}, ~~~\{\tau\}=\frac{1}{2}(\tau^{+}+\tau^{-}).
\end{align*}
For notational convenience, we also define the jump and mean on the boundary faces $E\in\varepsilon_{h}^{b}$ by modifying the above definitions appropriately. We use the definition of jump by understanding that $v^{-}=0$ (similarly, $\mathbf{v}^{-}=0$ and $\tau^{-}=0$) and the definition of mean by understanding that $v^{-}=v^{+}$(similarly, $\mathbf{v}^{-}=\mathbf{v}^{+}$ and $\tau^{-}=\tau^{+}$).\\
\indent We define the following discrete velocity and pressure spaces:
\begin{align*}
&\mathbf{V}_{h}=\{\mathbf{v}_{h}\in L^{2}(\Omega)^{d}: \mathbf{v}_{h}|_{\kappa}\in\mathbb{P}_{k}(\kappa)^{d},~\forall \kappa \in\pi_{h}\},\\
&Q_{h}=\{q_{h}\in Q:q_{h}|_{\kappa}\in \mathbb{P}_{k-1}(\kappa),~\forall \kappa \in\pi_{h}\},
\end{align*}
where $\mathbb{P}_{k}(\kappa)$ is the space of polynomials of degree less than or equal to $k\geq 1$ on $\kappa$.\\
\indent The DGFEM for the problem (\ref{s2.1}) is to find $(\lambda_{h}, \mathbf{u}_{h}, p_{h})\in \mathbb{R}^{+}\times\mathbf{V}_{h}\times Q_{h}$,
$\|\mathbf{u}_{h}\|_{0}=1$ such that
\begin{align}
\mathcal{A}_{h}(\mathbf{u}_{h},\mathbf{v}_{h})+\mathcal{B}_{h}(\mathbf{v}_{h},p_{h})&=\lambda_{h}(\mathbf{u}_{h},\mathbf{v}_{h}),  ~~\forall \mathbf{v}_{h}\in\mathbf{V}_{h},\label{s2.4}\\
\mathcal{B}_{h}(\mathbf{u}_{h},q_{h})&=0,~~~\forall q_{h}\in Q_{h},\label{s2.5}
\end{align}
where
\begin{align}
\mathcal{A}_{h}(\mathbf{u}_{h},\mathbf{v}_{h})&=\sum\limits_{\kappa\in\pi_{h}}\int_{\kappa}\nabla\mathbf{u}_{h}:\nabla\mathbf{v}_{h}dx-\sum\limits_{E\in\varepsilon_{h}}\int_{E}\{\nabla\mathbf{u}_{h}\}:[\![\underline{\mathbf{v}_{h}}]\!]ds\nonumber\\
&~~~-\sum\limits_{E\in\varepsilon_{h}}\int_{E}\{\nabla\mathbf{v}_{h}\}:[\![\underline{\mathbf{u}_{h}}]\!]ds
+\sum\limits_{E\in\varepsilon_{h}}\int_{E}\frac{\gamma}{h_{E}}[\![\underline{\mathbf{u}_{h}}]\!]:[\![\underline{\mathbf{v}_{h}}]\!]ds,\label{s2.6}\\
\mathcal{B}_{h}(\mathbf{v}_{h},q_{h})&=
-\sum\limits_{\kappa\in\pi_{h}}\int_{\kappa}q_{h}div\mathbf{v}_{h}dx+\sum\limits_{E\in\varepsilon_{h}}\int_{E}\{q_{h}\}[\![\mathbf{v}_{h}]\!]ds.\label{s2.7}
\end{align}
 Here $\gamma/h_{E}$ is the interior penalty parameter.
  From Remark 2.1 in \cite{Houston2007}, in the actual numerical implentations we can set $\gamma=C_{1}k^{2}$ with $C_{1}=10$ and $k$ is the degree of the polynomial.\\
\indent Define the DG-norm as follows:
\begin{eqnarray}
&&\|\mathbf{v}_{h}\|_{h}^{2}=\sum\limits_{\kappa\in\pi_{h}}\|\nabla \mathbf{v}_{h}\|_{0,\kappa}^{2}+\sum\limits_{E\in\varepsilon_{h}}\int\limits_{E}\frac{\gamma}{h_{E}}[\![\underline{\mathbf{v}_{h}}]\!]^{2}ds,~~on \mathbf{V}_{h}+\mathbf{V};\label{s2.8}\\
&&|||\mathbf{v_{h}}|||^{2}=\|\mathbf{v}_{h}\|_{h}^{2}+\sum\limits_{E\in \varepsilon_{h}}\int_{E}\frac{h_{E}}{\gamma}|\nabla \mathbf{v}_{h}|^{2}ds,~~on \mathbf{V}_{h}+H^{1+\frac{1}{2}}(\Omega)^{d}.\label{s2.9}
\end{eqnarray}
Note that $\|\cdot\|_{h}$ is equivalent to $\||\cdot\||$ on $\mathbf{V}_{h}$.\\
 \indent It is easy to know that (see \cite{Riviere2008}) the following continuity and coercivity properties hold:
\begin{eqnarray*}
&&|\mathcal{A}_{h}(\mathbf{u}_{h},\mathbf{v}_{h})|\lesssim \||\mathbf{u}_{h}\||\||\mathbf{v}_{h}\||,~~~ \forall  \mathbf{u}_{h}, \mathbf{v}_{h}\in \mathbf{V}_{h}+H^{1+s}(\Omega)^{d}~(s>\frac{1}{2}) \\
&& \|\mathbf{u}_{h}\|^{2}_{h}\lesssim \mathcal{A}_{h}(\mathbf{u}_{h},\mathbf{u}_{h}),
~~~ \forall \mathbf{u}_{h} \in \mathbf{V}_{h}.
\end{eqnarray*}

\indent From \cite{Riviere2006} we obtain the discrete inf-sup condition (the stability of the pressure):
\begin{align*}
\inf\limits_{p_{h}\in Q_{h}} \sup\limits_{\mathbf{v}_{h}\in \mathbf{V}_{h}} \frac{\mathcal{B}_{h}(\mathbf{v}_{h}, p_{h})}{\|\mathbf{v}_{h}\|_{h}\|p_{h}\|_{0}} \geq \beta^{*},
\end{align*}
where $\beta^{*}$ is a positive constant independent of $h$.\\

\indent We consider the source problem associated with the Stokes eigenvalue problem (\ref{s2.1}):
Given $\mathbf{f}\in (L^{2}(\Omega))^{d}$,
\begin{equation}\label{s2.10}
\begin{cases}
-\Delta \mathbf{u}^{f} +\nabla p^{f}=\mathbf{f},~~~ in~~\Omega,\\
div \mathbf{u}^{f}=0,~~~~~~~~~~~~~in ~~\Omega,\\
\mathbf{u}^{f}=0,~~~~~~~~~~~~~~~~~on~~\partial\Omega.
\end{cases}
\end{equation}
The weak formulation of (\ref{s2.10}) is given by: find $(\mathbf{u}^{f},p^{f})\in \mathbf{V}\times Q$ such that
\begin{align}
\mathcal{A}(\mathbf{u}^{f},\mathbf{v})+\mathcal{B}(\mathbf{v},p^{f})&=(\mathbf{f},\mathbf{v}),~~~~\forall\mathbf{v}\in\mathbf{V},\label{s2.11}\\
\mathcal{B}(\mathbf{u}^{f},q)&=0,~~~~~\forall q\in Q,\label{s2.12}
\end{align}
and its discontinuous Galerkin finite element form are as follows:
find  $(\mathbf{u}^{f}_{h},p^{f}_{h})\in \mathbf{V}_{h}\times Q_{h}$ such that
\begin{align}
\mathcal{A}_{h}(\mathbf{u}^{f}_{h},\mathbf{v}_{h})+\mathcal{B}_{h}(\mathbf{v}_{h},p^{f}_{h})&=(\mathbf{f},\mathbf{v}_{h}),
~~~~\forall\mathbf{v}_{h}\in\mathbf{V}_{h},\label{s2.13}\\
\mathcal{B}_{h}(\mathbf{u}^{f}_{h},q_{h})&=0,~~~~~~~~~\forall q_{h}\in Q_{h}.\label{s2.14}
\end{align}

\indent We assume that the following regularity is valid: for any $\mathbf{f}\in (L^{2}(\Omega))^{d}(d=2,3)$, there exists $(\mathbf{u}^{f}, p^{f})\in (H^{1+r}(\Omega)^{d}\times H^{r}(\Omega))\cap (W^{2,p}(\Omega)^{d}\times W^{1,p}(\Omega))~(\frac{1}{2}<r\leq 1,~p>\frac{2d}{d+1})$ satisfying (\ref{s2.10}) and
\begin{eqnarray}\label{s2.15}
\|\mathbf{u}^{f}\|_{1+r}+\|p^{f}\|_{r}\leq C\|\mathbf{f}\|_{0},
\end{eqnarray}
where $C$ is a positive constant independent of $\mathbf{f}$.\\

From Lemma 6.5 in \cite{Riviere2008}
we can obtain the consistency of the DGFEM, that is to say,
when $(\mathbf{u}^{f},p^{f})$ is the solution of the source problem (\ref{s2.10}), there hold the following equations:
\begin{align}
&\mathcal{A}_{h}(\mathbf{u}^{f},\mathbf{v}_{h})+\mathcal{B}_{h}(\mathbf{v}_{h},p^{f})=(\mathbf{f},\mathbf{v}_{h}),  ~~\forall \mathbf{v}_{h}\in\mathbf{V}_{h},\label{s2.16}\\
&\mathcal{B}_{h}(\mathbf{u}^{f},q_{h})=0,~~~\forall q_{h}\in Q_{h}.\label{s2.17}
\end{align}
From ({\ref{s2.13}})-(\ref{s2.14}) and ({\ref{s2.16}})-(\ref{s2.17}), we have
\begin{align}\label{s2.18}
&\mathcal{A}_{h}(\mathbf{u}^{f}-\mathbf{u}^{f}_{h},\mathbf{v}_{h})+\mathcal{B}_{h}(\mathbf{v}_{h}, p^{f}-p^{f}_{h})=0,~~~\forall \mathbf{v}_{h}\in \mathbf{V}_{h},\\
\label{s2.19}
&\mathcal{B}_{h}(\mathbf{u}^{f}-\mathbf{u}^{f}_{h},q_{h})=0,~~~\forall q_{h}\in Q_{h}.
\end{align}
\indent Since (\ref{s2.11})-(\ref{s2.12}) and (\ref{s2.13})-(\ref{s2.14}) are both uniquely solvable for each $\mathbf{f}\in L^{2}(\Omega)^{d}~(d=2,3)$ (see, e.g., Lemma 2.4 in \cite{Riviere2006}, and Lemma 7 and  Proposition 10 in \cite{Hansbo2002}), we can define the corresponding solution operators as follows:
\begin{eqnarray*}
&&~~~~~~~T: L^{2}(\Omega)^{d}\rightarrow \mathbf{V}, ~~~~T\mathbf{f}=\mathbf{u}^{f},\\
&&~~~~~~~T_{h}: L^{2}(\Omega)^{d}\rightarrow \mathbf{V}_{h}, ~~~~T_{h}\mathbf{f}=\mathbf{u}^{f}_{h},\\
&&~~~~~~~S: L^{2}(\Omega)^{d}\rightarrow Q, ~~~~~S\mathbf{f}=p^{f},\\
&&~~~~~~~S_{h}: L^{2}(\Omega)^{d}\rightarrow Q_{h}, ~~~~~S_{h}\mathbf{f}=p^{f}_{h}.
\end{eqnarray*}
Then (\ref{s2.2})-(\ref{s2.3}) and (\ref{s2.4})-(\ref{s2.5}) can be written in the following equivalent operator forms:
\begin{eqnarray}
  &\lambda T \mathbf{u}=\mathbf{u},~~S (\lambda \mathbf{u})=p,\label{s2.20}\\
    &\lambda_{h} T_{h} \mathbf{u}_{h}=\mathbf{u}_{h},~~S_{h} (\lambda_{h} \mathbf{u}_{h})=p_{h}.\label{s2.21}
\end{eqnarray}
It is easy to know that both $T$ and $T_{h}$ are self-adjoint and completely continuous and satisfy
\begin{eqnarray}
\|T\mathbf{f}\|_{1}+\|S\mathbf{f}\|_{0}\lesssim\|\mathbf{f}\|_{0},~~~\||T_{h}\mathbf{f}\||+\|S_{h}\mathbf{f}\|_{0}\lesssim\|\mathbf{f}\|_{0}~~~\label{s2.22}.
\end{eqnarray}

From Corollary 3.3 and Theorem 4.1 in \cite{Badia2014} we the following lemma. \\
\noindent{\bf Lemma 2.1.}~~
Assume $(\mathbf{u}^{f}, p^{f})\in H^{1+s}(\Omega)^{d}\times H^{s}(\Omega)$ for $ r< s\leq k$ and $\mathbf{f}\in H^{l}(\Omega)^{d}$ for $0\leq l\leq k+1$, then
\begin{eqnarray}\label{s2.23}
\|\mathbf{u}^{f}-\mathbf{u}^{f}_{h}\|_{h}+\|p^{f}-p^{f}_{h}\|_{0}\lesssim h^{s}(\|\mathbf{u}^{f}\|_{1+s}+\|p^{f}\|_{s})+h^{1+l}\|\mathbf{f}\|_{l}.
\end{eqnarray}
\indent Denote $I_{h}: \mathbf{V}\cap C^{0}(\overline{\Omega})^{d}\to \mathbf{V}_{h}\cap\mathbf{V}$ as the interpolation operator,
and denote $\varrho_{h}:~H^{s}(\Omega)\rightarrow Q_{h}$ as the local $L^{2}$ projection operator satisfying
$\varrho_{h}p|_{\kappa}\in \mathbb{P}_{k-1}(\kappa)$ and
\begin{eqnarray*}
 \int_{\kappa}(p-\varrho_{h}p)vdx=0,~~\forall v\in \mathbb{P}_{k-1}(\kappa),~~~\forall \kappa\in\pi_{h}.
\end{eqnarray*}
\indent Before estimating the error of velocity in the sense of $L^{2}$ norm, we introduce an auxiliary problem:
\begin{eqnarray}
&&\mathcal{A}(\omega,\mathbf{v})+\mathcal{B}(\mathbf{v},q)=(\mathbf{u}^{f}-\mathbf{u}^{f}_{h},\mathbf{v}),~~~~\forall\mathbf{v}\in\mathbf{V},
\\\label{s2.24}
&&\mathcal{B}(\omega,z)=0,~~~~~\forall z\in Q.\label{s2.25}
\end{eqnarray}
From (\ref{s2.15}) we have
\begin{align}\label{s2.26}
\|\omega\|_{1+r}+\|q\|_{r}\lesssim\|\mathbf{u}^{f}-\mathbf{u}^{f}_{h}\|_{0}.
\end{align}
Referring to Theorem 6.12 in \cite{Riviere2008}, by Nitsche's technique we can deduce the following lemma.\\
\noindent {\bf Lemma 2.2.}~~Suppose that the conditions of Lemma 2.1 and (\ref{s2.15}) hold, then
\begin{align}\label{s2.27}
\|\mathbf{u}^{f}-\mathbf{u}^{f}_{h}\|_{0}\lesssim h^{r}(\||\mathbf{u}^{f}-\mathbf{u}^{f}_{h}\||+\|p^{f}-p^{f}_{h}\|_{0}).
\end{align}

\indent By the above error estimates of the DG method for the source problem, next we can deduce the error estimates of the DG method for the eigenvalue problem.\\
By (\ref{s2.27}), (\ref{s2.23}) and (\ref{s2.15}), we have
\begin{eqnarray}
\|T_{h}-T\|_{0}\rightarrow 0,~~(h\rightarrow 0).\label{s2.29}
\end{eqnarray}
Thus, using Babu$\breve{s}$ka-Osborn spectral approximation theory \cite{Babuska1991book,D.Boffi2010}, we can get (see Lemma 2.3 in \cite{yzl2010}):\\
\noindent {\bf Lemma 2.3.}~~Assume that the regularity estimate (\ref{s2.15}) is valid.
Let $(\lambda,\mathbf{u}, p)$ and $(\lambda_{h},\mathbf{u}_{h}, p_{h})$ be the $j$th eigenpair of (\ref{s2.2})-(\ref{s2.3})
and (\ref{s2.4})-(\ref{s2.5}), respectively. Then
\begin{eqnarray}\label{s2.30}
\|\mathbf{u}_{h}-\mathbf{u}\|_{0}&\leq& C\|(T-T_{h})u\|_{0},\\\label{s2.31}
\lambda_{h}-\lambda&=&\lambda^{-2}((T-T_{h})\mathbf{u}, \mathbf{u})+R,
\end{eqnarray}
where $|R|\lesssim\|(T-T_{h})\mathbf{u}\|_{0}^{2}$.\\

From (\ref{s2.8}) and (\ref{s2.9}) we know that $|||\cdot|||$ is a norm stronger than $\|\cdot\|_{h}$, i.e., $\|v\|_{h}\lesssim|||v|||$.
Additionally, we have
\begin{eqnarray}\label{s2.32}
|||\mathbf{u}-\mathbf{u}_{h}|||^{2} \lesssim\|\mathbf{u}-\mathbf{u}_{h}\|_{h}^{2}+\sum\limits_{\kappa\in\mathcal{T}_{h}}h_{\kappa}^{2r}|\mathbf{u}-I_{h}\mathbf{u}|_{1+r,\kappa}^{2}.
\end{eqnarray}
In fact, from the inverse estimate, the interpolation estimate and the trace inequality, we deduce
\begin{eqnarray*}
&&\sum \limits_{E\in\varepsilon_{h}}h_{E}\|\nabla (\mathbf{u}-\mathbf{u}_{h})\|^{2}_{0,E}\lesssim\sum \limits_{E\in\varepsilon_{h}}h_{E}\|\nabla (I_{h}\mathbf{u}-\mathbf{u}_{h})\|^{2}_{0,E}
+\sum \limits_{E\in\varepsilon_{h}}h_{E}\|\nabla (\mathbf{u}-I_{h}\mathbf{u})\|^{2}_{0,E}\nonumber\\
&&~~~\lesssim\sum \limits_{\kappa\in\mathcal{T}_{h}}\|\nabla (I_{h}\mathbf{u}-\mathbf{u}_{h})\|^{2}_{0,\kappa}+\sum\limits_{\kappa\in\mathcal{T}_{h}}(h_{\kappa}^{-2}\|\mathbf{u}-I_{h}\mathbf{u}\|^{2}_{0,\kappa}
+|\mathbf{u}-I_{h}\mathbf{u}|^{2}_{1,\kappa}+h_{\kappa}^{2r}|\mathbf{u}-I_{h}\mathbf{u}|_{1+r,\kappa}^{2})\nonumber\\
&&~~~\lesssim \|\mathbf{u}-\mathbf{u}_{h}\|_{h}^{2}+\sum\limits_{\kappa\in\mathcal{T}_{h}}h_{\kappa}^{2r}|\mathbf{u}-I_{h}\mathbf{u}|_{1+r,\kappa}^{2}.
\end{eqnarray*}
By the above inequality and (\ref{s2.9}) we obtain (\ref{s2.32}).\\

\noindent {\bf Theorem 2.1.}~~Let $(\lambda,\mathbf{u}, p)$ and $(\lambda_{h},\mathbf{u}_{h}, p_{h})$ be the $j$th eigenpair of (\ref{s2.2})-(\ref{s2.3})
and (\ref{s2.4})-(\ref{s2.5}), respectively. Assume that the regularity estimate (\ref{s2.15}) is valid, and $(\mathbf{u},p)\in H^{1+s}(\Omega)^{d}\times H^{s}(\Omega)$ $(r\leq s\leq k)$. Then
\begin{eqnarray}
&&\|\mathbf{u}_{h}-\mathbf{u}\|_{0}\lesssim h^{r}(\||\mathbf{u}-\mathbf{u}_{h}\||+\|p-p_{h}\|_{0}),\label{s2.33}\\
&&\|\lambda_{h}-\lambda\|_{0}\lesssim h^{2s},\label{s2.34}\\
&&\|\mathbf{u}-\mathbf{u}_{h}\|_{h}+\|p-p_{h}\|_{0}\lesssim h^{s}(\|\mathbf{u}\|_{1+s}+\|p\|_{s})\label{s2.35}.
\end{eqnarray}
\noindent{\bf Proof.}~~Taking $\mathbf{f}=\lambda\mathbf{u}$ in (\ref{s2.11})-(\ref{s2.12}) and (\ref{s2.13})-(\ref{s2.14}),
then we get $\mathbf{u}^{f}=\lambda T\mathbf{u} $, $\mathbf{u}^{f}_{h}=\lambda T_{h}\mathbf{u}$,
$p^{f}=\lambda S\mathbf{u}$ and $p^{f}_{h}=\lambda S_{h}\mathbf{u}$.
Therefore, from (\ref{s2.23}) we have
\begin{eqnarray}\label{s2.36}
\|\lambda T\mathbf{u}-\lambda T_{h}\mathbf{u}\|_{h}+\|\lambda S\mathbf{u}-\lambda S_{h}\mathbf{u}\|_{0}
\lesssim h^{s}(\|\mathbf{u}\|_{1+s}+\|p\|_{s}).
\end{eqnarray}
From (\ref{s2.16}) and (\ref{s2.36}) we deduce
\begin{align}\label{s2.37}
((T-T_{h})\mathbf{u}, \mathbf{u})&=\mathcal{A}_{h}((T-T_{h})\mathbf{u}, T\mathbf{u})+\mathcal{B}_{h}((T-T_{h})\mathbf{u}, S\mathbf{u})\nonumber\\
&=\mathcal{A}_{h}((T-T_{h})\mathbf{u}, T\mathbf{u}-T_{h}\mathbf{u})+\mathcal{B}_{h}((T-T_{h})\mathbf{u}, S\mathbf{u}-S_{h}\mathbf{u})\nonumber\\
&\lesssim h^{2s}(\|\mathbf{u}\|_{1+s}+\|p\|_{s})^{2}.
\end{align}
Substituting (\ref{s2.37}) and (\ref{s2.27}) into (\ref{s2.31}) yields (\ref{s2.34}).\\
A simple calculation shows that
\begin{eqnarray*}
&&|\||\mathbf{u}-\mathbf{u}_{h}\||-\||\lambda T\mathbf{u}-\lambda T_{h}\mathbf{u}\|||
\leq \|\lambda_{h}T_{h}\mathbf{u}_{h}-\lambda T_{h}\mathbf{u}\|_{h}=\|T_{h}(\lambda_{h}\mathbf{u}_{h}-\lambda\mathbf{u}\|_{h}\lesssim \|\lambda_{h}\mathbf{u}_{h}-\lambda\mathbf{u}\|_{0},\nonumber\\
&&|\|p-p_{h}\|_{0}-\|\lambda S\mathbf{u}-\lambda S_{h}\mathbf{u}\|_{h}|
\leq \|\lambda_{h}S_{h}\mathbf{u}_{h}-\lambda S_{h}\mathbf{u}\|_{h}\lesssim \|\lambda_{h}\mathbf{u}_{h}-\lambda\mathbf{u}\|_{0},
\end{eqnarray*}
thus, from (\ref{s2.30}), (\ref{s2.31}), (\ref{s2.27}) and the above two estimates,
we deduce
\begin{eqnarray}\label{s2.38}
&&\|\lambda_{h}\mathbf{u}_{h}-\lambda\mathbf{u}\|_{0}\lesssim |\lambda_{h}-\lambda|+\|\mathbf{u}_{h}-\mathbf{u}\|_{0}\lesssim
\|\lambda T\mathbf{u}-\lambda T_{h}\mathbf{u}\|_{0}\nonumber\\
&&~~~\lesssim h^{r}(\||\lambda T\mathbf{u}-\lambda T_{h}\mathbf{u}\||+\|\lambda S\mathbf{u}-\lambda S_{h}\mathbf{u}\|_{h})\simeq h^{r}(\||\mathbf{u}-\mathbf{u}_{h}\||+\|p-p_{h}\|_{0}).
\end{eqnarray}
Thus, we get (\ref{s2.33}).
Since it is valid the relationship $\simeq$ in (\ref{s2.38}) and (\ref{s2.36}), we get (\ref{s2.35}).~~~$\Box$

\section{A posteriori error estimate for the Stokes eigenvalue problem}
\subsection {The a posteriori error indicator and its reliability for the eigenfunctions }

\indent Let $(\lambda_{h}, \mathbf{u}_{h}, p_{h})\in R^{+}\times\mathbf{V}_{h}\times Q_{h}$ be an eigenpair approximation.
To begin with, for each element $\kappa\in\pi_{h}$ we introduce the residuals
\begin{align*}
&\eta^{2}_{R_{\kappa}}=h_{\kappa}^{2}\|\lambda_{h}\mathbf{u}_{h}+\Delta\mathbf{u}_{h}-\nabla p_{h}\|^{2}_{0,\kappa}+\|div\mathbf{u}_{h}\|_{0,\kappa}^{2},\\
&\eta^{2}_{E_{\kappa}}=\frac{1}{2}\sum\limits_{E\subset\partial\kappa\setminus\partial\Omega}h_{E}\|[\![(p_{h}\mathbf{I}- \nabla\mathbf{u}_{h})]\!]\|^{2}_{0,E},
\end{align*}
where $\mathbf{I}$ denotes the $d\times d$ $(d=2, 3)$ identity matrix.
Next, we introduce the following estimator $\eta_{J_{\kappa}}$ to measure the jump of the approximate solution $\mathbf{u}_{h}$:
\begin{align*}
\eta^{2}_{J_{\kappa}}&=\sum\limits_{E\subset\partial\kappa, E\in\varepsilon^{i}_{h}}\gamma h_{E}^{-1}|[\![\underline{\mathbf{u}_{h}]\!]}|_{0,E}^{2}+ \sum\limits_{E\subset\partial\kappa, E\in\varepsilon^{b}_{h}}\gamma h_{E}^{-1}|\mathbf{u}_{h}\otimes\mathbf{n}|_{0,E}^{2}.
\end{align*}
The local error indictor 
is defined as
\begin{align*}
\eta^{2}_{\kappa}=\eta^{2}_{R_{\kappa}}+\eta^{2}_{E_{\kappa}}+\eta^{2}_{J_{\kappa}}.
\end{align*}
Finally, we introduce the global a posteriori error estimator
\begin{align*}
\eta_{h}=(\sum\limits_{\kappa\in\pi_{h}}\eta^{2}_{\kappa})^{\frac{1}{2}}.
\end{align*}

For $\kappa\in\pi_{h}$, denote $\theta_{\kappa}=$int$\{\bigcup\limits_{\overline{\kappa}_{i}\cap \overline{\kappa}\not=\emptyset}\bar{\kappa}_{i}, \kappa_{i}\in\pi_{h}\}$ and $\theta_{E}$ is the set of all elements which share at least one node with face $E$.
Let $\mathbf{v}^{I}$ be the Scott-Zhang interpolation function  \cite{ScottZhang1990}, then $\mathbf{v}^{I}\in \mathbf{V}\cap\mathbf{V}_{h}$ and
\begin{eqnarray}
&&\|\mathbf{v}-\mathbf{v}^{I}\|_{0,\kappa}+h_{\kappa}\|\nabla(\mathbf{v}-\mathbf{v}^{I})\|_{0,\kappa}\lesssim h_{\kappa}|\mathbf{v}|_{1,\theta_{\kappa}},~~~~~~~\forall \kappa\in\pi_{h},\label{s3.1}\\
&&\|\mathbf{v}-\mathbf{v}^{I}\|_{0,E} \lesssim h_{E}^{\frac{1}{2}}|\mathbf{v}|_{1,\theta_{E}},~~~~~\forall E\subset\partial\kappa.\label{s3.2}
\end{eqnarray}
\noindent Denote
\begin{eqnarray*}
\underline{\sum}_{h}=\{\underline{\tau}\in L^{2}(\Omega)^{d\times d}: \underline{\tau}|_{\kappa}\in \mathbb{P}_{k}(\kappa)^{d\times d}, \kappa\in\pi_{h}\}.
\end{eqnarray*}
We introduce the lifting operator $\mathcal{L}:\mathbf{V}+\mathbf{V}_{h}\rightarrow \underline{\sum}_{h}$ by
\begin{align}\label{s3.3}
\int_{\Omega}\mathcal{L}(\mathbf{v}):  \underline{\tau}dx=\sum\limits_{E\in\varepsilon_{h}^{i}}\int_{E}[\![\underline{\mathbf{v}}]\!]:\{\underline{\tau}\}ds,~~\forall \underline{\tau}\in \underline{\sum}_{h}.
\end{align}
Moreover, from \cite{Perugia2002,Schotzau2002}, the lifting operator has the ability property
\begin{align}\label{s3.4}
\|\mathcal{L}(\mathbf{v})\|_{0}^{2}\lesssim \sum\limits_{E\in\varepsilon_{h}^{i}}\|h_{E}^{-\frac{1}{2}}[\![\underline{\mathbf{v}}]\!]\|_{0,E}^{2},~~~\forall \mathbf{v}\in \mathbf{V}+\mathbf{V}_{h}.
\end{align}
Using this operator, we introduce an auxiliary bilinear form
\begin{align}\label{s3.5}
\widetilde{\mathcal{A}}_{h}(\cdot, \cdot): V(h)\times V(h)\rightarrow R
\end{align}
defined by
\begin{align}\label{s3.6}
\widetilde{\mathcal{A}}_{h}(\mathbf{w},\mathbf{v})=\sum\limits_{\kappa\in\pi_{h}}\int_{\kappa}\nabla\mathbf{w}:\nabla\mathbf{v}dx
-\sum\limits_{\kappa\in\pi_{h}}\int_{\kappa}\mathcal{L}(\mathbf{v}):\nabla\mathbf{w}dx
-\sum\limits_{\kappa\in\pi_{h}}\int_{\kappa}\mathcal{L}(\mathbf{w}):\nabla\mathbf{v}dx   +\sum\limits_{E\in\varepsilon_{h}}\int_{E}\frac{\gamma}{h_{E}}[\![\underline{\mathbf{w}}]\!]:[\![\underline{\mathbf{v}}]\!]ds.
\end{align}
Since $\widetilde{\mathcal{A}}_{h}=\mathcal{A}_{h}$ on $V_{h}\times V_{h}$, the DGFEM presented in (\ref{s2.4})-(\ref{s2.5}) is equivalent to finding $(\lambda_{h}, \mathbf{u}_{h}, p_{h})\in R^{+}\times V_{h}\times Q_{h}$ and satisfying
\begin{align}
&\widetilde{\mathcal{A}}_{h}(\mathbf{u}_{h}, \mathbf{v}_{h})+\mathcal{B}_{h}(\mathbf{v}_{h}, p_{h})=\lambda_{h}(\mathbf{u}_{h}, \mathbf{v}_{h}),  ~~\forall \mathbf{v}_{h}\in\mathbf{V}_{h},\label{s3.7}\\
&\mathcal{B}_{h}(\mathbf{u}_{h},q_{h})=0,~~~\forall q_{h}\in Q_{h}.\nonumber
\end{align}

\noindent{\bf Lemma 3.1.}~~Let $(\mathbf{u}^{f},p^{f})$ and $(\mathbf{u}^{f}_{h},p^{f}_{h})$ be the solutions of (\ref{s2.11})-(\ref{s2.12}) and (\ref{s2.13})-(\ref{s2.14}), respectively. Then
\begin{align}\label{s3.8}
\|\mathbf{u}^{f}-\mathbf{u}^{f}_{h}\|_{h}+\|p^{f}-p^{f}_{h}\|_{0}\backsimeq \sup\limits_{0\not=\mathbf{v}\in \mathbf{V}}\frac{|(\mathbf{f},\mathbf{v})-\widetilde{\mathcal{A}}_{h}(\mathbf{u}^{f}_{h}, \mathbf{v})-\mathcal{B}_{h}(\mathbf{v}, p^{f}_{h})|}{\|\mathbf{v}\|_{h}}+\inf\limits_{\mathbf{v}\in \mathbf{V}}\|\mathbf{u}^{f}_{h}-\mathbf{v}\|_{h}.
\end{align}
\noindent{\bf Proof.}~~
For $\forall \mathbf{\bar{u}}\in \mathbf{V}$, from (\ref{s2.11}) we have
\begin{eqnarray*}
&&~~~\widetilde{\mathcal{A}}_{h}(\mathbf{u}^{f}-\mathbf{\bar{u}}, \mathbf{u}^{f}-\mathbf{\bar{u}})=\widetilde{\mathcal{A}}_{h}(\mathbf{u}^{f}, \mathbf{u}^{f}-\mathbf{\bar{u}})-\widetilde{\mathcal{A}}_{h}(\mathbf{\bar{u}}, \mathbf{u}^{f}-\mathbf{\bar{u}})\nonumber\\
&&=(\mathbf{f}, \mathbf{u}^{f}-\mathbf{\bar{u}})-\mathcal{B}_{h}( \mathbf{u}^{f}-\mathbf{\bar{u}}, p^{f})-\widetilde{\mathcal{A}}_{h}(\mathbf{u}^{f}_{h}, \mathbf{u}^{f}-\mathbf{\bar{u}})+\widetilde{\mathcal{A}}_{h}(\mathbf{u}^{f}_{h}-\mathbf{\bar{u}},  \mathbf{u}^{f}-\mathbf{\bar{u}}).
\end{eqnarray*}
For $\forall \mathbf{\bar{u}}\in \mathbf{V}$, $\overline{p}\in Q$, we have
\begin{eqnarray*}
\mathcal{B}_{h}(\mathbf{u}^{f}-\mathbf{\bar{u}}, p^{f}-\bar{p})=\mathcal{B}_{h}(\mathbf{u}^{f}-\mathbf{\bar{u}}, p^{f})-\mathcal{B}_{h}(\mathbf{u}^{f}-\mathbf{\bar{u}}, p^{f}_{h})-\mathcal{B}_{h}(\mathbf{u}^{f}-\mathbf{\bar{u}},\bar{p}- p^{f}_{h}).
\end{eqnarray*}
Suming the above two equations and taking $\mathbf{v}=\mathbf{u}^{f}-\mathbf{\bar{u}}$, we deduce
\begin{align}\label{s3.9}
&\|\mathbf{u}^{f}-\mathbf{\bar{u}}\|_{h}\|\mathbf{v}\|_{h}+\mathcal{B}_{h}(\mathbf{v}, p^{f}-\bar{p})
=(\mathbf{f}, \mathbf{v})-\mathcal{B}_{h}(\mathbf{v}, p^{f})-\widetilde{\mathcal{A}}_{h}(\mathbf{u}^{f}_{h}, \mathbf{v})\nonumber\\
&~~~~~~+\widetilde{\mathcal{A}}_{h}(\mathbf{u}^{f}_{h}-\mathbf{\bar{u}}, \mathbf{v})+
\mathcal{B}_{h}(\mathbf{v}, p^{f})-\mathcal{B}_{h}(\mathbf{v}, p^{f}_{h})-\mathcal{B}_{h}(\mathbf{v}, \bar{p}- p^{f}_{h})\nonumber\\
&~~~=(\mathbf{f}, \mathbf{v})-\widetilde{\mathcal{A}}_{h}(\mathbf{u}^{f}_{h}, \mathbf{v})
-\mathcal{B}_{h}(\mathbf{v}, p^{f}_{h})+\widetilde{\mathcal{A}}_{h}(\mathbf{u}^{f}_{h}-\mathbf{\bar{u}},\mathbf{v})
-\mathcal{B}_{h}(\mathbf{v},\bar{p}- p^{f}_{h}).
\end{align}
By the inf-sup condition we obtain
$$\sup\limits_{\mathbf{v}\in \mathbf{V}}\frac{\mathcal{B}_{h}(\mathbf{v}, p^{f}-\bar{p})}{\|\mathbf{v}\|_{h}}\gtrsim\|p^{f}-\bar{p}\|_{0}, $$
then, dividing both sides of (\ref{s3.9}) by $\|\mathbf{v}\|_{h}$ and taking supremum for $\mathbf{v}\in \mathbf{V}$, we get
\begin{eqnarray}\label{s3.10}
&&\|\mathbf{u}^{f}-\mathbf{\bar{u}}\|_{h}+\| p^{f}-\bar{p}\|_{0}\lesssim\sup\limits_{\mathbf{v}\in\mathbf{V}}\frac{(\mathbf{f}, \mathbf{v})-\widetilde{\mathcal{A}}_{h}(\mathbf{u}^{f}_{h}, \mathbf{v})
-\mathcal{B}_{h}(\mathbf{v}, p^{f}_{h})}{\|\mathbf{v}\|_{h}}\nonumber\\
&&~~~~~~+\|\mathbf{u}^{f}_{h}-\mathbf{\bar{u}}\|_{h}+\|\bar{p}- p^{f}_{h}\|_{0},~~~\forall (\mathbf{\bar{u}}, \bar{p})\in \mathbf{V}\times Q.
\end{eqnarray}
From the triangle inequality we have
\begin{eqnarray}\label{s3.11}
&&\|\mathbf{u}^{f}-\mathbf{u}^{f}_{h}\|_{h}+\| p^{f}-p^{f}_{h}\|_{0}\lesssim\sup\limits_{\mathbf{v}\in\mathbf{V}}\frac{(\mathbf{f}, \mathbf{v})-\widetilde{\mathcal{A}}_{h}(\mathbf{u}^{f}_{h}, \mathbf{v})
-\mathcal{B}_{h}(\mathbf{v}, p_{h}^{f})}{\|\mathbf{v}\|_{h}}\nonumber\\
&&~~~~~~+\|\mathbf{u}^{f}_{h}-\mathbf{\bar{u}}\|_{h}+\|\bar{p}- p^{f}_{h}\|_{0},~~~\forall (\mathbf{\bar{u}}, \bar{p})\in \mathbf{V}\times Q.
\end{eqnarray}
Since $(\mathbf{\bar{u}}, \bar{p})$ is arbitrary and $\inf\limits_{\bar{p}\in Q}\|\bar{p}- p^{f}_{h}\|_{0}=0$, the part $\lesssim$ in (\ref{s3.8}) is valid.
The other part $\gtrsim$ in (\ref{s3.8}) is obvious.
~~$\Box$\\
\indent Lemma 3.1 can be extended to the eigenvalue problems.\\
\noindent {\bf Theorem 3.1.}~~Let $(\lambda,\mathbf{u}, p)$ and $(\lambda_{h},\mathbf{u}_{h}, p_{h})$ be the $j$th eigenpair of (\ref{s2.2})-(\ref{s2.3})
and (\ref{s2.4})-(\ref{s2.5}), respectively. Then
\begin{eqnarray}\label{s3.12}
\|\mathbf{u}-\mathbf{u}_{h}\|_{h}+\|p-p_{h}\|_{0}\backsimeq \sup\limits_{0\not=\mathbf{v}\in \mathbf{V}}\frac{|\widetilde{\mathcal{A}}_{h}(\mathbf{u}-\mathbf{u}_{h},\mathbf{v})+\mathcal{B}_{h}(\mathbf{v},p-p_{h})|}{\|\mathbf{v}\|_{h}}+\inf\limits_{\mathbf{v}\in \mathbf{V}}\|\mathbf{u}_{h}-\mathbf{v}\|_{h}.
\end{eqnarray}
\noindent{\bf Proof.}~~By ({\ref{s2.22}}) we deduce
\begin{eqnarray}
&&\|\mathbf{u}-\mathbf{u}_{h}\|_{h}+\|p-p_{h}\|_{0}=
\|\lambda T\mathbf{u}-\lambda_{h} T\mathbf{u}_{h }+\lambda_{h} T\mathbf{u}_{h }-\lambda_{h}T_{h}\mathbf{u}_{h}\|_{h}\\\nonumber
&&~~~~~~+\|\lambda S\mathbf{u}-\lambda_{h} S\mathbf{u}_{h }+\lambda_{h} S\mathbf{u}_{h }-\lambda_{h}S_{h}\mathbf{u}_{h}\|_{h}\\\nonumber
&&~~~\leq \|\lambda_{h} T\mathbf{u}_{h }-\lambda_{h}T_{h}\mathbf{u}_{h}\|_{h}
+\|\lambda_{h} S\mathbf{u}_{h }-\lambda_{h}S_{h}\mathbf{u}_{h}\|_{h}
+\|\lambda \mathbf{u}-\lambda_{h} \mathbf{u}_{h }\|_{0}.\label{s3.13}
\end{eqnarray}
Taking $\mathbf{f}=\lambda_{h}\mathbf{u}_{h}$ in (\ref{s2.11})-(\ref{s2.12}) and (\ref{s2.13})-(\ref{s2.14}),
then we get $\mathbf{u}^{f}=\lambda_{h} T\mathbf{u}_{h }$, $\mathbf{u}^{f}_{h}=\lambda_{h}T_{h}\mathbf{u}_{h}$,
$p^{f}=\lambda_{h}S\mathbf{u}_{h }$ and $p^{f}_{h}=\lambda_{h}S_{h}\mathbf{u}_{h}$.
 Therefore, from (\ref{s3.8}) we have
\begin{eqnarray}\label{s3.14}
&&~~~\|\lambda_{h} T\mathbf{u}_{h }-\lambda_{h}T_{h}\mathbf{u}_{h}\|_{h}+\|\lambda_{h} S\mathbf{u}_{h }-\lambda_{h}S_{h}\mathbf{u}_{h}\|_{0}\nonumber\\
&&\lesssim \sup\limits_{0\not=\mathbf{v}\in \mathbf{V}}\frac{|(\mathbf{\lambda_{h}u_{h}},\mathbf{v})-\widetilde{\mathcal{A}}_{h}(\lambda_{h}T_{h}\mathbf{u}_{h}, \mathbf{v})-\mathcal{B}_{h}(\mathbf{v}, S_{h}(\lambda_{h}\mathbf{u}_{h}) )|}{\|\mathbf{v}\|_{h}}
+\inf\limits_{\mathbf{v}\in \mathbf{V}}\|\mathbf{u}_{h}-\mathbf{v}\|_{h}.
\end{eqnarray}
From (\ref{s2.11}) with $\mathbf{f}=\lambda_{h}\mathbf{u}_{h}$ and (\ref{s2.22}), we deduce
\begin{eqnarray}\label{s3.15}
&&~~~|(\mathbf{\lambda_{h}u_{h}},\mathbf{v})-\widetilde{\mathcal{A}}_{h}(\lambda_{h}T_{h}\mathbf{u}_{h}, \mathbf{v})-\mathcal{B}_{h}(\mathbf{v}, S_{h}(\lambda_{h}\mathbf{u}_{h}) )|\nonumber\\
&&=|\widetilde{\mathcal{A}}_{h}(\lambda_{h}T\mathbf{u}_{h}, \mathbf{v})+\mathcal{B}_{h}(\mathbf{v}, S(\lambda_{h}\mathbf{u}_{h}))-\widetilde{\mathcal{A}}_{h}(\lambda_{h}T_{h}\mathbf{u}_{h}, \mathbf{v})-\mathcal{B}_{h}(\mathbf{v}, S_{h}(\lambda_{h}\mathbf{u}_{h}))|\nonumber\\
&&=|\widetilde{\mathcal{A}}_{h}(\lambda_{h}T\mathbf{u}_{h}-\lambda_{h}T_{h}\mathbf{u}_{h}, \mathbf{v})+\mathcal{B}_{h}(\mathbf{v},S(\lambda_{h}\mathbf{u}_{h})- S_{h}(\lambda_{h}\mathbf{u}_{h}) )|\nonumber\\
&&=|\widetilde{\mathcal{A}}_{h}(\lambda_{h}T\mathbf{u}_{h}-\lambda T\mathbf{u}+\mathbf{u}-\mathbf{u}_{h}, \mathbf{v})+\mathcal{B}_{h}(\mathbf{v},S(\lambda_{h}\mathbf{u}_{h})-S(\lambda\mathbf{u})+p- p_{h})|\nonumber\\
&&\leq |\widetilde{\mathcal{A}}_{h}(\mathbf{u}-\mathbf{u}_{h}, \mathbf{v})+\mathcal{B}_{h}(\mathbf{v},p- p_{h}) )|
+C \|\lambda_{h}\mathbf{u}_{h}-\lambda\mathbf{u}\|_{0}\|\mathbf{v}\|_{h}.
\end{eqnarray}
Substituting (\ref{s3.15}) into (\ref{s3.14}), we get
\begin{eqnarray}\label{s3.16}
&&~~~\|\lambda_{h} T\mathbf{u}_{h }-\lambda_{h}T_{h}\mathbf{u}_{h}\|_{h}+\|\lambda_{h} S\mathbf{u}_{h }-\lambda_{h}S_{h}\mathbf{u}_{h}\|_{0}\nonumber\\
&&\lesssim \sup\limits_{0\not=\mathbf{v}\in \mathbf{V}}\frac{|\widetilde{\mathcal{A}}_{h}(\mathbf{u}-\mathbf{u}_{h},\mathbf{v})+\mathcal{B}_{h}(\mathbf{v},p- p_{h}) |}{\|\mathbf{v}\|_{h}}+C( \|\lambda_{h}\mathbf{u}_{h}-\lambda\mathbf{u}\|_{0}+\inf\limits_{\mathbf{v}\in \mathbf{V}}\|\mathbf{u}_{h}-\mathbf{v}\|_{h}).
\end{eqnarray}
From  Theorem 2.1 we know that $\|\lambda_{h}\mathbf{u}_{h}-\lambda\mathbf{u}\|_{0}$ is a small quantity of higher order compared with $\|\mathbf{u}-\mathbf{u}_{h}\|_{h}+\|p-p_{h}\|_{0}$. Substituting (\ref{s3.16}) into (\ref{s3.13}), the side $\lesssim$ in (\ref{s3.12}) is true.
The other side $\gtrsim$ in (\ref{s3.12}) is obvious.
~~$\Box$\\
\noindent {\bf Lemma 3.2.}~~Under the conditions of Theorem 2.1, there holds
\begin{align}\label{s3.17}
\widetilde{\mathcal{A}}_{h}(\mathbf{u}-\mathbf{u}_{h}, \mathbf{v})+\mathcal{B}_{h}(\mathbf{v}, p-p_{h})
\lesssim \sum\limits_{\kappa\in\pi_{h}}\left(\eta_{R_{\kappa}}+\eta_{E_{\kappa}}+\eta_{J_{\kappa}}\right)\|\mathbf{v}\|_{h}+\|\lambda \mathbf{u}-\lambda_{h} \mathbf{u}_{h}\|_{0}\|\mathbf{v}\|_{h},
~~~\forall v\in \mathbf{V}.
\end{align}
\noindent {\bf Proof.}~~ Using (\ref{s3.7}), (\ref{s3.6}), (\ref{s2.7}) and Green's formula we deduce that
\begin{align}\label{s3.18}
&\widetilde{\mathcal{A}}_{h}(\mathbf{u}-\mathbf{u}_{h},\mathbf{v})
+\mathcal{B}_{h}(\mathbf{v},p-p_{h})=\widetilde{\mathcal{A}}_{h}(\mathbf{u},\mathbf{v})-\widetilde{\mathcal{A}}_{h}(\mathbf{u}_{h},\mathbf{v})+\mathcal{B}_{h}(\mathbf{v},p)-\mathcal{B}_{h}(\mathbf{v},p_{h})\nonumber\\
&~~~=\lambda (\mathbf{u},\mathbf{v})-\widetilde{\mathcal{A}}_{h}(\mathbf{u}_{h},\mathbf{v})-\mathcal{B}_{h}(\mathbf{v},p_{h})\nonumber\\
&~~~=\lambda\sum\limits_{\kappa\in\pi_{h}}\int_{\kappa}\mathbf{u}\mathbf{v}dx-\sum\limits_{\kappa\in\pi_{h}}\int_{\kappa}\nabla\mathbf{u}_{h}:\nabla\mathbf{v}dx
+\sum\limits_{\kappa\in\pi_{h}}\int_{\kappa}\mathcal{L}(\mathbf{v}):\nabla\mathbf{u}_{h}dx\nonumber\\
&~~~~~+\sum\limits_{\kappa\in\pi_{h}}\int_{\kappa}\mathcal{L}(\mathbf{u}_{h}):\nabla\mathbf{v}dx
+\sum\limits_{E\in\varepsilon_{h}}\int_{E}\frac{\gamma}{h_{E}}[\![\underline{\mathbf{u}_{h}}]\!]:[\![\underline{\mathbf{v}}]\!]ds+\sum\limits_{\kappa\in\pi_{h}}
\int_{\kappa}div\mathbf{v}p_{h}dx\nonumber\\
&~~~=\lambda\sum\limits_{\kappa\in\pi_{h}}\int_{\kappa}\mathbf{u}\mathbf{v}dx+\sum\limits_{\kappa\in\pi_{h}}\int_{\kappa}\Delta\mathbf{u}_{h}\cdot\mathbf{v}dx- \sum\limits_{\kappa\in\pi_{h}}\sum\limits_{E\subset\partial\kappa}\int_{E}\frac{\partial\mathbf{u}_{h}}{\partial\mathbf{n}}\cdot\mathbf{v}ds     +\sum\limits_{\kappa\in\pi_{h}}\int_{\kappa}\mathcal{L}(\mathbf{v}):\nabla\mathbf{u}_{h}dx\nonumber\\
&~~~~~~+\sum\limits_{\kappa\in\pi_{h}}
\int_{\kappa}\mathcal{L}(\mathbf{u}_{h}):\nabla\mathbf{v}dx -\sum\limits_{\kappa\in\pi_{h}}\int_{\kappa}\nabla p_{h}\cdot\mathbf{v}dx
+\sum\limits_{\kappa\in\pi_{h}}\sum\limits_{E\subset\partial\kappa}\int_{E}p_{_h}\mathbf{v}\cdot\mathbf{n}ds.
\end{align}
By $\mathbf{v}^{I}\in \mathbf{V}\cap\mathbf{V}_{h}$, (\ref{s2.2})-(\ref{s2.3}) and (\ref{s2.4})-(\ref{s2.5}), we obtain
\begin{align*}
\widetilde{\mathcal{A}}_{h}(\mathbf{u}-\mathbf{u}_{h}, \mathbf{v})+\mathcal{B}_{h}(\mathbf{v}, p-p_{h})=\widetilde{\mathcal{A}}_{h}(\mathbf{u}-\mathbf{u}_{h}, \mathbf{v}-\mathbf{v}^{I})+\mathcal{B}_{h}(\mathbf{v}-\mathbf{v}^{I}, p-p_{h})+\sum\limits_{\kappa\in\pi_{h}}\int_{\kappa}(\lambda \mathbf{u}-\lambda_{h}\mathbf{u}_{h})\cdot\mathbf{v}^{I}dx.
\end{align*}
Using (\ref{s3.6}), Cauchy-Schwartz inequality, (\ref{s3.1}) and (\ref{s3.2}),  (\ref{s3.18}) can be written as follows:
\begin{align}\label{s3.20}
&~~~\widetilde{\mathcal{A}}_{h}(\mathbf{u}-\mathbf{u}_{h},\mathbf{v})+\mathcal{B}_{h}(\mathbf{v},p-p_{h})\nonumber\\
&=\lambda\sum\limits_{\kappa\in\pi_{h}}\int_{\kappa}\mathbf{u}(\mathbf{v}-\mathbf{v}^{I})dx+\sum\limits_{\kappa\in\pi_{h}}
\int_{\kappa}\Delta\mathbf{u}_{h}\cdot(\mathbf{v}-\mathbf{v}^{I})dx- \sum\limits_{\kappa\in\pi_{h}}\sum\limits_{E\subset\partial\kappa}\int_{E}\frac{\partial\mathbf{u}_{h}}{\partial\mathbf{n}}\cdot(\mathbf{v}-\mathbf{v}^{I})ds     \nonumber\\
&~~~+\sum\limits_{\kappa\in\pi_{h}}\int_{\kappa}\mathcal{L}(\mathbf{v}-\mathbf{v}^{I}):\nabla\mathbf{u}_{h}dx+\sum\limits_{\kappa\in\pi_{h}}\int_{\kappa}\mathcal{L}(\mathbf{u}_{h}):\nabla(\mathbf{v}-\mathbf{v}^{I})dx -\sum\limits_{\kappa\in\pi_{h}}\int_{\kappa}\nabla p_{h}\cdot(\mathbf{v}-\mathbf{v}^{I})dx\nonumber\\
&~~~+\sum\limits_{\kappa\in\pi_{h}}\sum\limits_{E\subset\partial\kappa}\int_{E}p_{_h}(\mathbf{v}-\mathbf{v}^{I})\cdot\mathbf{n}ds+\sum\limits_{\kappa\in\pi_{h}}\int_{\kappa}(\lambda \mathbf{u}-\lambda_{h}\mathbf{u}_{h})\cdot\mathbf{v}^{I}dx\nonumber\\
&=\sum\limits_{\kappa\in\pi_{h}}\int_{\kappa}(\Delta\mathbf{u}_{h}+\lambda_{h}\mathbf{u}_{h}-\nabla p_{h})\cdot(\mathbf{v}-\mathbf{v}^{I})dx- \sum\limits_{\kappa\in\pi_{h}}\sum\limits_{E\subset\partial\kappa}\int_{E}\frac{\partial\mathbf{u}_{h}}{\partial\mathbf{n}}\cdot(\mathbf{v}-\mathbf{v}^{I})ds     \nonumber\\
&~~~+\sum\limits_{\kappa\in\pi_{h}}\int_{\kappa}\mathcal{L}(\mathbf{v}-\mathbf{v}^{I}):\nabla\mathbf{u}_{h}dx+\sum\limits_{\kappa\in\pi_{h}}\int_{\kappa}\mathcal{L}(\mathbf{u}_{h}):\nabla(\mathbf{v}-\mathbf{v}^{I})dx \nonumber\\
&~~~+\sum\limits_{\kappa\in\pi_{h}}\sum\limits_{E\subset\partial\kappa}\int_{E}p_{_h}(\mathbf{v}-\mathbf{v}^{I})\cdot\mathbf{n}ds
+\sum\limits_{\kappa\in\pi_{h}}\int_{\kappa}(\lambda \mathbf{u}-\lambda_{h}\mathbf{u}_{h})\cdot\mathbf{v}dx\nonumber\\
&=\sum\limits_{\kappa\in\pi_{h}}\int_{\kappa}(\Delta\mathbf{u}_{h}+\lambda_{h}\mathbf{u}_{h}-\nabla p_{h})\cdot(\mathbf{v}-\mathbf{v}^{I})dx+\sum\limits_{\kappa\in\pi_{h}}\sum\limits_{E\subset\partial\kappa\backslash\partial\Omega }\int_{E}(p_{h}\mathbf{I}-\nabla\mathbf{u}_{h})\mathbf{n}\cdot(\mathbf{v}-\mathbf{v}^{I})ds \nonumber\\
&~~~+\sum\limits_{\kappa\in\pi_{h}}\int_{\kappa}\mathcal{L}(\mathbf{v}-\mathbf{v}^{I}):\nabla\mathbf{u}_{h}dx
+\sum\limits_{\kappa\in\pi_{h}}\int_{\kappa}\mathcal{L}(\mathbf{u}_{h}):\nabla(\mathbf{v}-\mathbf{v}^{I})dx +\sum\limits_{\kappa\in\pi_{h}}\int_{\kappa}(\lambda \mathbf{u}-\lambda_{h}\mathbf{u}_{h})\cdot\mathbf{v}dx\nonumber\\
&\equiv B_{1}+B_{2}+B_{3}+B_{4}+B_{5}.
\end{align}
Next, we will analyze each item on the right-hand side of (\ref{s3.20}). Using the Cauchy-Schwarz inequality and the approximation property (\ref{s3.1}) and (\ref{s3.2}), we have
\begin{align*}
|B_{1}|&\leq \sum\limits_{\kappa\in\pi_{h}}\|\Delta\mathbf{u}_{h}+\lambda_{h}\mathbf{u}_{h}-\nabla p_{h}\|_{0,\kappa}
\|\mathbf{v}-\mathbf{v}^{I}  \|_{0,\kappa}\\
&\lesssim\sum\limits_{\kappa\in\pi_{h}}h_{\kappa}\|\Delta\mathbf{u}_{h}+\lambda_{h}\mathbf{u}_{h}-\nabla p_{h}\|_{0,\kappa}
\|\mathbf{v}\|^{2}_{1,\theta_{\kappa}}\\
&\lesssim\left(\sum\limits_{\kappa\in\pi_{h}}h_{\kappa}^{2}\|\Delta\mathbf{u}_{h}+\lambda_{h}\mathbf{u}_{h}-\nabla p_{h}\|_{0,\kappa}^{2}\right)^{\frac{1}{2}}
\|\mathbf{v}\|_{h}.
\end{align*}
For the second term on the right-hand side of (\ref{s3.20}), from (\ref{s3.2}) we obtain
\begin{align*}
|B_{2}|&=|\frac{1}{2}\sum\limits_{\kappa\in\pi_{h}}\sum\limits_{E\subset\partial\kappa\backslash\partial\Omega} \int_{E}[\![p_{h}\mathbf{I}-\nabla\mathbf{u}_{h}]\!]\cdot(\mathbf{v}-\mathbf{v}^{I})ds|\\
&\leq \sum\limits_{\kappa\in\pi_{h}}\sum\limits_{E\subset\partial\kappa\backslash\partial\Omega} \|[\![p_{h}\mathbf{I}-\nabla\mathbf{u}_{h}]\!]\|_{0,E}
Ch_{E}^{\frac{1}{2}}\|\mathbf{v}\|_{1,\theta_{E}}\\
&\lesssim\left(\sum\limits_{\kappa\in\pi_{h}}\sum\limits_{E\subset\partial\kappa\backslash\partial\Omega}( h_{E}^{\frac{1}{2}}\|[\![p_{h}\mathbf{I}-\nabla\mathbf{u}_{h}]\!]\|_{0,E})^{2}    \right)^{\frac{1}{2}}
\|\mathbf{v}\|_{h}.
\end{align*}
For the third term,
by the properties of the interpolation function $\mathbf{v}^{I}$, we know $[\![\mathbf{v}-\mathbf{v}^{I}]\!]=0$. Therefore, from the definition of lifting operation $\mathcal{L}$ we have
\begin{align*}
B_{3}=\sum\limits_{\kappa\in\pi_{h}}\int_{\kappa}\mathcal{L}(\mathbf{v}-\mathbf{v}^{I}):\nabla\mathbf{u}_{h}ds
=\sum\limits_{E\in\varepsilon_{h}}\int_{E}\{\underline{\nabla\mathbf{u}_{h}}\}:[\![\underline{\mathbf{v}-\mathbf{v}^{I}}]\!]ds=0.
\end{align*}
For the fourth term, using the Cauchy-Schwarz inequality, (\ref{s3.4}) and (\ref{s3.1}) we get
\begin{align*}
|B_{4}|&\leq \left(\sum\limits_{\kappa\in\pi_{h}}\|\mathcal{L}(\mathbf{u}_{h})\|_{0,\kappa}^{2}\right) ^{\frac{1}{2}}
\left(\sum\limits_{\kappa\in\pi_{h}}\|\nabla(\mathbf{v}-\mathbf{v}^{I})\|_{0,\kappa}^{2}\right) ^{\frac{1}{2}}\\
&\lesssim\left(\sum\limits_{E\in\varepsilon_{h}^{i}}\|h_{E}^{-\frac{1}{2}}[\![\underline{\mathbf{u}_{h}}]\!]\|_{0,E}^{2} \right)^{\frac{1}{2}}  \left(\sum\limits_{\kappa\in\pi_{h}} \|\nabla(\mathbf{v}-\mathbf{v}^{I})\|_{0,\kappa}^{2}\right) ^{\frac{1}{2}}
\lesssim\left(\sum\limits_{E\in\varepsilon_{h}^{i}}\|h_{E}^{-\frac{1}{2}}[\![\underline{\mathbf{u}_{h}]\!]}\|_{0,E}^{2} \right)^{\frac{1}{2}}\|\mathbf{v}\|_{h}.
\end{align*}
For the last term on the right-hand side of (\ref{s3.20}), we have
\begin{align*}
B_{5}=\sum\limits_{\kappa\in\pi_{h}}\int_{\kappa}(\lambda \mathbf{u}-\lambda_{h}\mathbf{u}_{h})\cdot\mathbf{v}dx
\leq \|\lambda\mathbf{u}-\lambda_{h}\mathbf{u}_{h}\|_{0}\|\mathbf{v}\|_{0}.
\end{align*}
Substituting $B_{1}, B_{2}, B_{3}, B_{4}, B_{5}$ into (\ref{s3.20}), we obtain the desired result (\ref{s3.17}).  ~~$\Box$\\

In \cite{Ohannes2003,Brenner2003} the authors construct the enriching operator $E_{h}:\mathbf{V}_{h}\to \mathbf{V}_{h}\cap\mathbf{V}$ by averaging and prove the following lemma.\\
\noindent{\bf Lemma 3.3.}~It is valid the following estimate:
\begin{eqnarray}
 \|\mathbf{u}_{h}-E_{h}\mathbf{u}_{h}\|_{h}\lesssim\sum\limits_{E\in\varepsilon_{h}^{i}}\gamma h_{E}^{-1}|[\![\underline{\mathbf{u}_{h}}]\!]|_{0,E}^{2}+ \sum\limits_{E\in\varepsilon_{h}^{b}}\gamma h_{E}^{-1}|\mathbf{u}_{h}\otimes\mathbf{n}|_{0,E}^{2}.\label{s3.21}
\end{eqnarray}

\noindent {\bf Theorem 3.2.}~~Suppose that the conditions of Theorem 2.1 hold, then
 there holds
\begin{eqnarray}\label{s3.22}
~~~~~~\|\mathbf{u}-\mathbf{u}_{h}\|_{h}+\|p-p_{h}\|_{0}
\lesssim\eta_{h}.
\end{eqnarray}
\noindent{\bf Proof.}~~Substituting (\ref{s3.17}) and (\ref{s3.21}) into (\ref{s3.12}), we obtain (\ref{s3.22}).~~~  $\Box$

\subsection {The efficiency of the indicators for eigenfunctions}
\indent This section is devoted to prove an efficiency bound for $\eta$. To prove the results, we
use the bubble function technique which was introduced in \cite{Verfurth1996}.\\
\indent Let $\kappa$ be an element of $\pi_{h}$. Let $b_{\kappa}$ and $b_{E}$  be the standard bubble function
on element $\kappa$ and face $E$ ($d=3$) or edge $E$ ($d=2$) of $\kappa$, respectively.
Then, from \cite{Verfurth1996,oden2011, Kanschat2008} we have the following results.\\
\noindent{\bf Lemma 3.4.}~~For any vector-valued polynomial function $\mathbf{v}_{h}$ on $\kappa$, there hold
\begin{align}
\|\mathbf{v}_{h}\|_{0,\kappa}&\lesssim \|b_{\kappa}^{1/2}\mathbf{v}_{h}\|_{0,\kappa}, \label{s3.23}\\
\|b_{\kappa}\mathbf{v}_{h}\|_{0,\kappa}&\lesssim \|\mathbf{v}_{h}\|_{0,\kappa},\label{s3.24}\\
\|\nabla (b_{\kappa}\mathbf{v}_{h})\|_{0,\kappa}&\lesssim h_{\kappa}^{-1}\|\mathbf{v}_{h}\|_{0,\kappa}.\label{s3.25}
\end{align}
For any vector-valued polynomial function $\sigma$ on $E$, it is valid that
\begin{align}
\|b_{E}\sigma\|_{0,E}&\lesssim \|\sigma\|_{0,E},\label{s3.26}\\
\|\sigma\|_{0,E}&\lesssim \|b_{E}^{1/2}\sigma\|_{0,E}.\label{s3.27}
\end{align}
Furthermore, for each $b_{E}\sigma$, there exists an extension $\sigma_{b}\in H_{0}^{1}(\omega(E))$ satisfying $\sigma_{b}|_{E}=b_{E}\sigma$ and
\begin{align}
\|\sigma_{b}\|_{0,\kappa}&\lesssim h_{E}^{1/2}\|\sigma\|_{0,E},~~~\forall \kappa\in\omega(E),\label{s3.28}\\
\|\nabla \sigma_{b}\|_{0,\kappa}&\lesssim h_{E}^{-1/2}\|\sigma\|_{0,E},~~~\forall \kappa\in\omega(E).\label{s3.29}
\end{align}
\indent From the above lemma and using the standard arguments (see \cite{Zeng2017}, Lemma 3.13), we can prove the following local bounds.\\
\noindent{\bf Lemma 3.5.}~~Under the conditions of Theorem 2.1, there holds
\begin{align}
\eta_{R_{\kappa}}\lesssim
 \|\nabla (\mathbf{u}-\mathbf{u}_{h})\|_{0,\kappa}+\|p-p_{h}\|_{0,\kappa}
+h_{\kappa}\|\lambda_{h}\mathbf{u}_{h}-\lambda \mathbf{u}\|_{0,\kappa}.\label{s3.30}
\end{align}
\noindent {\bf Proof.}~~ For any $\kappa\in \pi_{h}$, define the function $R$ and $W$ locally by
\begin{align*}
R|_{\kappa}=\lambda_{h}\mathbf{u}_{h}+\triangle\mathbf{u}_{h}-\nabla p_{h}~~~and~~~W|_{\kappa}=h_{\kappa}^{2}Rb_{\kappa}.
\end{align*}
From (\ref{s3.23}) and using $\lambda\mathbf{u}+\Delta \mathbf{u}-\nabla p=0$, we have
\begin{align*}
h_{\kappa}^{2}\|R\|_{0,\kappa}^{2}&\lesssim \int_{\kappa}R\cdot (h_{\kappa}^{2}Rb_{\kappa})dx
=\int_{\kappa}(\lambda_{h}\mathbf{u}_{h}+\Delta \mathbf{u}_{h}-\nabla p_{h})\cdot Wdx\\
&=\int_{\kappa}(\lambda_{h}\mathbf{u}_{h}+\Delta \mathbf{u}_{h}-\nabla p_{h}-(\lambda\mathbf{u}+\Delta \mathbf{u}-\nabla p))\cdot Wdx\\
&=\int_{\kappa}\Delta(\mathbf{u}_{h}-\mathbf{u})\cdot Wdx-\int_{\kappa}\nabla(p_{h}-p)\cdot Wdx+\int_{\kappa}(\lambda_{h}\mathbf{u}_{h}-\lambda\mathbf{u})\cdot Wdx.
\end{align*}
Using integration by parts and $W|_{\partial\kappa}=0$, we obtain
\begin{align*}
h_{\kappa}^{2}\|R\|_{0,\kappa}^{2}&\lesssim \int_{\kappa}\nabla (\mathbf{u}-\mathbf{u}_{h})\cdot \nabla Wdx+\int_{\kappa}(p_{h}-p) divWdx+\int_{\kappa}(\lambda_{h}\mathbf{u}_{h}-\lambda\mathbf{u})\cdot Wdx.
\end{align*}
Applying Cauchy-Schwarz inequality yields
\begin{align}\label{s3.31}
h_{\kappa}^{2}\|R\|_{0,\kappa}^{2}\lesssim\left( \|\nabla (\mathbf{u}-\mathbf{u}_{h})\|_{0,\kappa}+\|p-p_{h}\|_{0,\kappa}
+h_{\kappa}\|\lambda_{h}\mathbf{u}_{h}-\lambda \mathbf{u}\|_{0,\kappa}\right)\left(\|\nabla W\|_{0,\kappa}+h_{\kappa}^{-1}\|W\|_{0,\kappa}\right).
\end{align}
From (\ref{s3.24}) and (\ref{s3.25}) we get
\begin{align*}
\|\nabla W\|_{0,\kappa}+h_{\kappa}^{-1}\|W\|_{0,\kappa}\lesssim h_{\kappa}\|R\|_{0,\kappa}.
\end{align*}
Dividing (\ref{s3.31}) by $ h_{\kappa}\|R\|_{0,\kappa}$ and noting $\|\nabla\cdot \mathbf{u}_{h}\|_{0}=\|\nabla\cdot (\mathbf{u}_{h}-\mathbf{u})\|_{0}$, we finish the proof.
~~~$\Box$\\

\noindent{\bf Lemma 3.6.~~}Under the conditions of Theorem 2.1, there holds
\begin{align*}
\eta^{2}_{E_{\kappa}}
\lesssim \|\nabla (\mathbf{u}-\mathbf{u}_{h})\|_{0,\omega(\kappa)}+\|p-p_{h}\|_{0,\omega(\kappa)}
+\left( \sum\limits_{\kappa\in\omega(\kappa)}h_{\kappa}^{2} \|\lambda\mathbf{u}-\lambda_{h}\mathbf{u}_{h}\|_{0,\kappa}^{2}  \right)^{\frac{1}{2}}.
\end{align*}
\noindent{\bf Proof.}~~ For any interior edge $E\in\varepsilon_{h}^{i}$, let the function $R$ and $\Lambda$ be such that
\begin{align*}
R|_{E}=[\![p_{h}\mathbf{I}-\nabla\mathbf{u}_{h} ]\!]|_{E}~~~and ~~~\Lambda=h_{E}Rb_{E}.
\end{align*}
Using (\ref{s3.27}) and $[\![p\mathbf{I}-\nabla\mathbf{u}]\!]|_{E}=0$, we get
\begin{align*}
h_{E}\|R\|_{0,E}^{2}\lesssim \int_{E}R\cdot (h_{E}Rb_{E})ds=\int_{E}[\![(p_{h}-p)\mathbf{I}-\nabla(\mathbf{u}_{h}-\mathbf{u})]\!]\cdot\Lambda ds.
\end{align*}
Applying Green's formula over each of the two elements of $\omega(E)$, we derive
\begin{align*}
&h_{E}\|R\|_{0,E}^{2}\lesssim\int_{E}[\![((p_{h}-p)\mathbf{I}-\nabla(\mathbf{u}_{h}-\mathbf{u}))]\!]\cdot\Lambda ds\\
&=C (\sum\limits_{\kappa\in\omega(E)}\int_{\kappa}(-\Delta (\mathbf{u}-\mathbf{u}_{h})+\nabla (p-p_{h}))\cdot\Lambda dx-\sum\limits_{\kappa\in\omega(E)}\int_{\kappa}(\nabla (\mathbf{u}-\mathbf{u}_{h})-(p-p_{h})\mathbf{I}):\nabla\Lambda dx)
\end{align*}
Using $\lambda\mathbf{u}+\Delta\mathbf{u}-\nabla p=0$, we obtain
\begin{eqnarray}
&&h_{E}\|R\|_{0,E}^{2}\lesssim \sum\limits_{\kappa\in\omega(E)}\int_{K}(\lambda_{h}\mathbf{u}_{h}+\Delta\mathbf{u}_{h}-\nabla p_{h})\cdot\Lambda dx+\sum\limits_{\kappa\in\omega(E)}\int_{K}(\lambda\mathbf{u}-\lambda_{h}\mathbf{u}_{h})\cdot\Lambda dx\nonumber\\
&&~~~+\sum\limits_{\kappa\in\omega(E)}\int_{K}(-\nabla (u-u_{h})+(p-p_{h})\mathbf{I}):\nabla\Lambda dx\equiv T_{1}+T_{2}+T_{3}.\label{s3.32}
\end{eqnarray}
Using Cauchy-Schwarz inequality, shape-regularity of the mesh, (\ref{s3.28}) and (\ref{s3.29}) yieids
\begin{eqnarray*}
&&T_{1}\lesssim \left(\sum\limits_{\kappa\in\omega(E)}\eta^{2}_{R_{\kappa}}\right)^{1/2}\left(\sum\limits_{\kappa\in\omega(E)}h^{-2}_{\kappa}\|\Lambda\|_{0,\kappa}^{2}     \right)^{1/2}\lesssim \left(\sum\limits_{\kappa\in\omega(E)}\eta^{2}_{R_{\kappa}}\right)^{1/2}h_{E}^{1/2}\|R\|_{0,E},\\
&&T_{2}\lesssim \left( \sum\limits_{\kappa\in\omega(E)}\left(h^{2}_{\kappa}\|\lambda\mathbf{u}-\lambda_{h}\mathbf{u}_{h}    \|^{2}_{0,\kappa}\right)\right)^{1/2}h_{E}^{1/2}\|R\|_{0,E},\\
&&T_{3}\lesssim \left(\sum\limits_{\kappa\in\omega(E)}(\|\nabla (\mathbf{u}-\mathbf{u}_{h}) \|^{2}_{0,\kappa}+\|p-p_{h}\| ^{2}_{0,\kappa}) \right)^{1/2}h_{E}^{1/2}\|R\|_{0,E}.
\end{eqnarray*}
Combing the above estimates of $T_{1}$, $T_{2}$ and $T_{3}$, dividing (\ref{s3.32}) by $h_{E}^{1/2}\|R\|_{0,E}$ and summing over all interior edges of $\kappa$, we get the desired result.~~~  $\Box$ \\

\noindent{\bf Lemma 3.7.}~~Under the conditions of Theorem 2.1, there holds
\begin{align*}
\eta^{2}_{J_{\kappa}}&=\sum\limits_{E\subset\partial\kappa, E\in\varepsilon^{i}_{h}}\gamma h_{E}^{-1}|[\![\underline{\mathbf{u}_{h}-\mathbf{u}}]\!]|_{0,E}^{2}+ \sum\limits_{E\subset\partial\kappa, E\in\varepsilon^{b}_{h}}\gamma h_{E}^{-1}|(\mathbf{u}_{h}-\mathbf{u})\otimes\mathbf{n}|_{0,E}^{2}.
\end{align*}
\noindent {\bf Proof.}~For any $E\in \varepsilon_{h}^{i}(\Omega)$,  $[\![\underline{\mathbf{u}}]\!]=0$, and
for any $E\in \varepsilon_{h}\cap \partial\Omega$,  $\mathbf{u}\otimes\mathbf{n}=0$. Therefore, we obtain the desired result. ~~~$\Box$\\

\noindent {\bf Theorem 3.3.}~~Suppose that the conditions of Theorem 2.1 hold. Then the a posteriori error estimator $\eta_{h}$ is efficient:
\begin{eqnarray}
&&\eta^{2}_{\kappa}\lesssim\sum\limits_{\kappa\in\omega(\kappa)}(
\|\mathbf{u}-\mathbf{u}_{h}\|^{2}_{0,\kappa}+\|p-p_{h}\|^{2}_{0,\kappa}+h_{k}^{2}
\|\lambda\mathbf{u}-\lambda_{h}\mathbf{u}_{h}\|^{2}_{0,\kappa}),\\\label{s3.33}
&&\eta^{2}_{h}\lesssim
\|\mathbf{u}-\mathbf{u}_{h}\|^{2}_{h}+\|p-p_{h}\|^{2}+\sum\limits_{\kappa\in\pi_{h}}h_{k}^{2}\|\lambda\mathbf{u}-\lambda_{h}\mathbf{u}_{h}\|^{2}_{0,\kappa}      .\label{s3.34}
\end{eqnarray}
\noindent{\bf Proof.}~~ The conclusions follow from a combination of Lemmas 3.5-3.7.   ~~~$\Box$

\subsection{The reliability of the indicators for the eigenvalues}

\noindent{\bf Lemma 3.8.}~~Let $(\lambda, \mathbf{u},p)$ and $(\lambda_{h}, \mathbf{u}_{h},p_{h})$ be the eigenpairs of  (\ref{s2.2})-(\ref{s2.3})  and (\ref{s2.4})-(\ref{s2.5}), respectively,
then
\begin{eqnarray}\label{s3.35}
\lambda_{h}-\lambda=\mathcal{A}_{h}(\mathbf{u}-\mathbf{u}_{h}, \mathbf{u}-\mathbf{u}_{h})+2\mathcal{B}_{h}(\mathbf{u}-\mathbf{u}_{h},p-p_{h})-\lambda(\mathbf{u}-\mathbf{u}_{h}, \mathbf{u}-\mathbf{u}_{h}).
\end{eqnarray}
\noindent{\bf Proof.}~~By the consistency formulas (\ref{s2.16})-(\ref{s2.17}) we get
\begin{align}
&\mathcal{A}_{h}(\mathbf{u},\mathbf{v}_{h})+\mathcal{B}_{h}(\mathbf{v}_{h},p)=\lambda(\mathbf{u},\mathbf{v}_{h}),  ~~\forall \mathbf{v}_{h}\in\mathbf{V}_{h},\label{s3.36}\\
&\mathcal{B}_{h}(\mathbf{u},q_{h})=0,~~~\forall q_{h}\in Q_{h}.\label{s3.37}
\end{align}
From (\ref{s2.2})-(\ref{s2.3}) with $(\mathbf{v},q)=(\mathbf{u},p)$, (\ref{s2.4})-(\ref{s2.5}) with $(\mathbf{v}_{h},q_{h})=(\mathbf{u_{h}},p_{h})$
and (\ref{s3.36})-(\ref{s3.37}), we deduce
\begin{eqnarray*}
&&\mathcal{A}_{h}(\mathbf{u}-\mathbf{u}_{h}, \mathbf{u}-\mathbf{u}_{h})+2\mathcal{B}_{h}(\mathbf{u}-\mathbf{u}_{h},p-p_{h})-\lambda(\mathbf{u}-\mathbf{u}_{h}, \mathbf{u}-\mathbf{u}_{h})\nonumber\\
&&~~~=\mathcal{A}_{h}(\mathbf{u}, \mathbf{u})-2\mathcal{A}_{h}(\mathbf{u},\mathbf{u}_{h})+\mathcal{A}_{h}(\mathbf{u}_{h},\mathbf{u}_{h})
+2\mathcal{B}_{h}(\mathbf{u}, p)-2\mathcal{B}_{h}(\mathbf{u}_{h},p)\nonumber\\
&&~~~~~~-2\mathcal{B}_{h}(\mathbf{u},p_{h})+2\mathcal{B}_{h}(\mathbf{u}_{h},p_{h})-\lambda(\mathbf{u}, \mathbf{u})+2\lambda(\mathbf{u},\mathbf{u}_{h})-\lambda(\mathbf{u}_{h},\mathbf{u}_{h})\nonumber\\
&&~~~=\lambda_{h}(\mathbf{u}_{h},\mathbf{u}_{h})-\lambda(\mathbf{u}_{h},\mathbf{u}_{h})=\lambda_{h}-\lambda.
\end{eqnarray*}
The proof is completed. ~~~$\Box$\\

\noindent{\bf Theorem 3.4.}~~ Under the conditions of Theorem 2.1, there holds
\begin{eqnarray}\label{s3.38}
|\lambda-\lambda_{h}|\lesssim \eta_{h}^{2}+\sum\limits_{\kappa\in\mathcal{T}_{h}}h_{\kappa}^{2r}(|u-I_{h}\mathbf{u}|_{1+r,\kappa}^{2}
+\|p-\varrho_{h}p\|_{r}^{2}).
\end{eqnarray}
\noindent{\bf Proof.}~~ Theorem 2.1 shows $\|u-u_{h}\|_{0,\Omega}$  is a term of higher order than
$\||\mathbf{u}-\mathbf{u}_{h}\||+\|p-p_{h}\|_{0}$.
Hence, from (\ref{s3.35}) and (\ref{s3.22}), we obtain
\begin{eqnarray*}
|\lambda-\lambda_{h}|\lesssim \||\mathbf{u}-\mathbf{u}_{h}\||^{2}+\|p-p_{h}\|_{0}^{2}+\sum\limits_{E\in\varepsilon_{h}}h_{E}\|p-p_{h}\|_{0,E}^{2}.
\end{eqnarray*}
Thus, from (\ref{s2.32}) and (\ref{s3.22}) we obtain (\ref{s3.38}).~~~$\Box$\\

\noindent{\bf Remark 3.1.}~~From Theorems 3.2 and 3.3, we know the indicator $\eta_{h}$
 of the eigenfunction error $\|\mathbf{u}-\mathbf{u}_{h}\|_{h}+\|p-p_{h}\|_{0}$
is reliable and efficient up to data oscillation, so
the adaptive algorithm based on the indicator can generate a good graded mesh, which makes the
eigenfunction error $\|\mathbf{u}-\mathbf{u}_{h}\|_{h}+\|p-p_{h}\|_{0}$ can achieve
the optimal convergence rate $O(dof^{-\frac{2k}{d}})$. Thus, referring to \cite{Wu2007,Yang2020} we are able to expect to get
$\sum\limits_{\kappa\in\mathcal{T}_{h}}h_{\kappa}^{2r}(|u-I_{h}\mathbf{u}|_{1+r,\kappa}^{2}
+\|p-\varrho_{h}p\|_{r}^{2})
\lesssim dof^{-\frac{2k}{d}}$
, thereby from (\ref{s3.38}) we have $|\lambda-\lambda_{h}|\lesssim dof^{-\frac{2k}{d}}$. Therefore, we think that $\eta_{h}^{2}$
 can be viewed as the error indicator of $\lambda_{h}$.
 The numerical experiments in Section 5 show
$\eta_{h}^{2}$ as the error indicator of $\lambda_{h}$ is reliable and efficient.\\
\noindent{\bf Remark 3.2.}~~Assume that $\Omega$ can be subdivided into shape-regular affine meshes $\pi_{h}$ consisting of parallelograms $\kappa~(d =
2)$ or parallelepipeds $\kappa~ (d = 3)$, and the discrete velocity and pressure spaces are given by
\begin{align*}
&\mathbf{V}_{h}=\{\mathbf{v}_{h}\in L^{2}(\Omega)^{d}: \mathbf{v}_{h}|_{\kappa}\in\mathbb{Q}_{k}(\kappa)^{d},~\forall \kappa\in\pi_{h}\},\\
&Q_{h}=\{q_{h}\in Q:q_{h}|_{\kappa}\in \mathbb{Q}_{k-1}(\kappa),~\forall \kappa\in\pi_{h}\},
\end{align*}
where $\mathbb{Q}_{k}(\kappa)$ denotes the space of tensor product polynomials on $\kappa$ of
degree $k$ in each coordinate direction.\\
For the Stokes equation (\ref{s2.10}), Houston et al.\cite{Houston2005} studied the
 a posteriori error estimation of mixed DGFEM using the above $\mathbb{Q}_{k}-\mathbb{Q}_{k-1}$ element.
For the Stokes eigenvalue problem (\ref{s2.1}), all analysis and conclusions in this paper are valid for the mixed DGFEM using the above $\mathbb{Q}_{k}-\mathbb{Q}_{k-1}$ element.

\section{Numerical experiments}
\indent Using the a posteriori error indicators in this paper and consulting the existing standard algorithms (see,
e.g., \cite{Dai2008}), we present the following algorithm.\\
{\bf  Algorithm 4.1.}~~ Choose the parameter $0 <\theta< 1$.\\
{\bf Step 1}.  Pick any initial mesh $\pi_{h_{0}}$ with mesh size $h_{0}$ .\\
{\bf Step 2}. Solve (\ref{s2.4})-(\ref{s2.5}) on $\pi_{h_{0}}$ for discrete solution $(\lambda_{h_{0}},\mathbf{u}_{h_{0}}, p_{h_{0}})$.\\
{\bf Step 3}.  Let $l=0$.\\
{\bf Step 4}.  Compute the local indicator $\eta_{\kappa}^{2}$.\\
{\bf Step 5}. Construct $\widehat{\pi}_{h_{l}}\subset\pi_{h_{l}}$ by {\bf Marking Strategy E} .\\
{\bf Step 6}. Refine $\pi_{h_{l}}$ to get a new mesh $\pi_{h_{l+1}}$ by procedure {\bf REFINE} .\\
{\bf Step 7}. Solve (\ref{s2.4})-(\ref{s2.5}) on $\pi_{h_{l+1}}$ for discrete solution $(\lambda_{h_{l+1}},\mathbf{u}_{h_{l+1}}, p_{h_{l+1}})$.\\
{\bf Step 8}. Let $l=l+1$ and go to {\bf step} 4.\\
{\bf Marking Strategy E}. \\
\indent Given parameter $0<\theta <1$:\\
{\bf Step 1}. Construct a minimal subset $\widehat{\pi}_{h_{l}}\subset\pi_{h_{l}}$ by selecting some elements in $\pi_{h_{l}}$ such that
\begin{eqnarray*}
\sum\limits_{\kappa\in \widehat{\pi}_{h_{l}}}\eta_{\kappa}^{2}\geq \theta \eta_{h_{l}}^{2}.
\end{eqnarray*}
{\bf Step 2}. Mark all elements in $\widehat{\pi}_{h_{l}}$.\\
\indent The above marking strategy was introduced by D$\ddot{o}$rfler \cite{dorfler1996} (see also Morin et al. \cite{morin2002} ).\\
\indent We use the following notations in our tables:\\
$l$: the $l$th iteration in Algorithm 4.1.\\
$\lambda_{1,h_{l}}$: the first discrete eigenvalue at the $l$th iteration of Algorithm 4.1.\\
$dof$: the degrees of freedom at the $l$th iteration.\\
$-$: the calculation cannot proceed since the computer runs out of memory.\\
\indent We carry out experiments in $d$-dimensional cases ($d$=2,~3). Our program is compiled under the package of iFEM \cite{Chen2009} and we use internal command $'eigs'$ in MATLAB to solve matrix eigenvalue problem.\\
\subsection{The results in two-dimensional domains}
\indent We carry out experiments on three two-dimensional domains:  $\Omega_{C}=(-1,1)^{2}\setminus \{{0\leq x\leq 1,y=0}\}$, $\Omega_{L}=(-1,1)^{2}\setminus [0,1]\times [-1,0]$  and $\Omega_{S}=(0,1)^{2}$.
The discrete eigenvalue problems are solved in MATLAB 2018a on a DELL PC with 1.80GHZ CPU and 32GB RAM.
 We take $\theta=0.5$ and initial mesh $\pi_{h_{0}}$ ($h_{0}=\frac{\sqrt{2}}{16}$ ) for three two-dimensional domains.
 To compute the error of approximations of the first eigenvalue, we take $\lambda_{C} = 29.9168629$, $\lambda_{L} = 32.13269465$  and $\lambda_{S} = 52.344691168$ (see \cite{Gedicke2019})
as the reference values for two-dimensional domains $\Omega_{C}$, $\Omega_{L}$ and $\Omega_{S}$, respectively.\\
\indent The adaptive refined meshes and the error curves are shown in Figures 1-8.
 We show some adaptively refined meshes for $k=1,2,3$ on the left side of Figures 1-8
from which we can see the strongly refinement towards the tip
of the slit at the origin for $\Omega_{C}$ and $\Omega_{L}$ and a clear refinement near the four corners of $\Omega_{S}$. Furthermore, from Figures 1-8 we can see that the error curves and error indicators curves  for DG methods using $\mathbb{P}_{k}-\mathbb{P}_{k-1}$ $(k=1,2,3)$ element are both approximately parallel to the line with slope $-k$, which indicates the error indicators
are reliable and efficient and the adaptive algorithm can reach the optimal convergence order. It coincides with our theoretical results.
 It also can be seen from error curves that under the same $dof$, the approximations obtained by adaptive algorithm are more accurate than those computed on uniform meshes,
 and the approximations obtained by high order elements are more accurate than those computed by low order elements on both uniform meshes and adaptive meshes.\\
\indent The approximations of the first eigenvalue for $\Omega_{C}$, $\Omega_{L}$ and $\Omega_{S}$ using $\mathbb{P}_{3}-\mathbb{P}_{2}$ element are listed in Tables 1-3. These eigenvalues have the same accuracy as those \cite{Gedicke2018, Gedicke2019}, which further proves that our method is effective.


\begin{figure}[H]
\begin{minipage}[t]{0.5\textwidth}
  \centering
  \includegraphics[width=7.3cm]{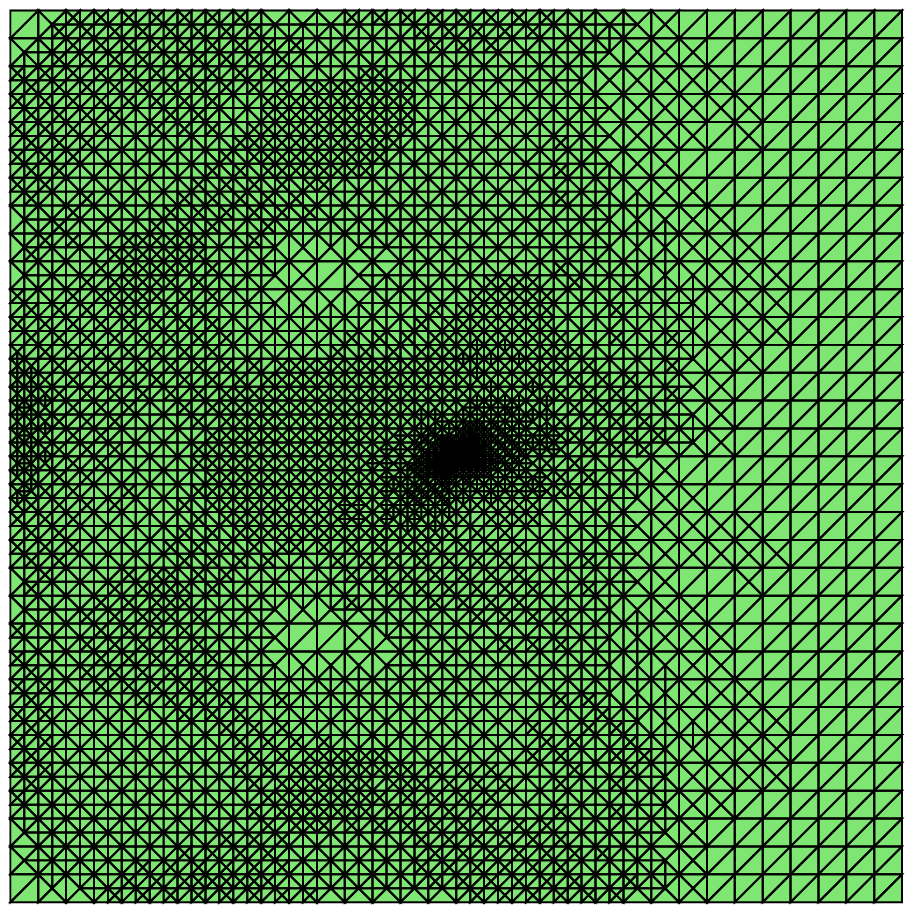}
\end{minipage}
\begin{minipage}[t]{0.5\textwidth}
  \centering
  \includegraphics[width=7.3cm]{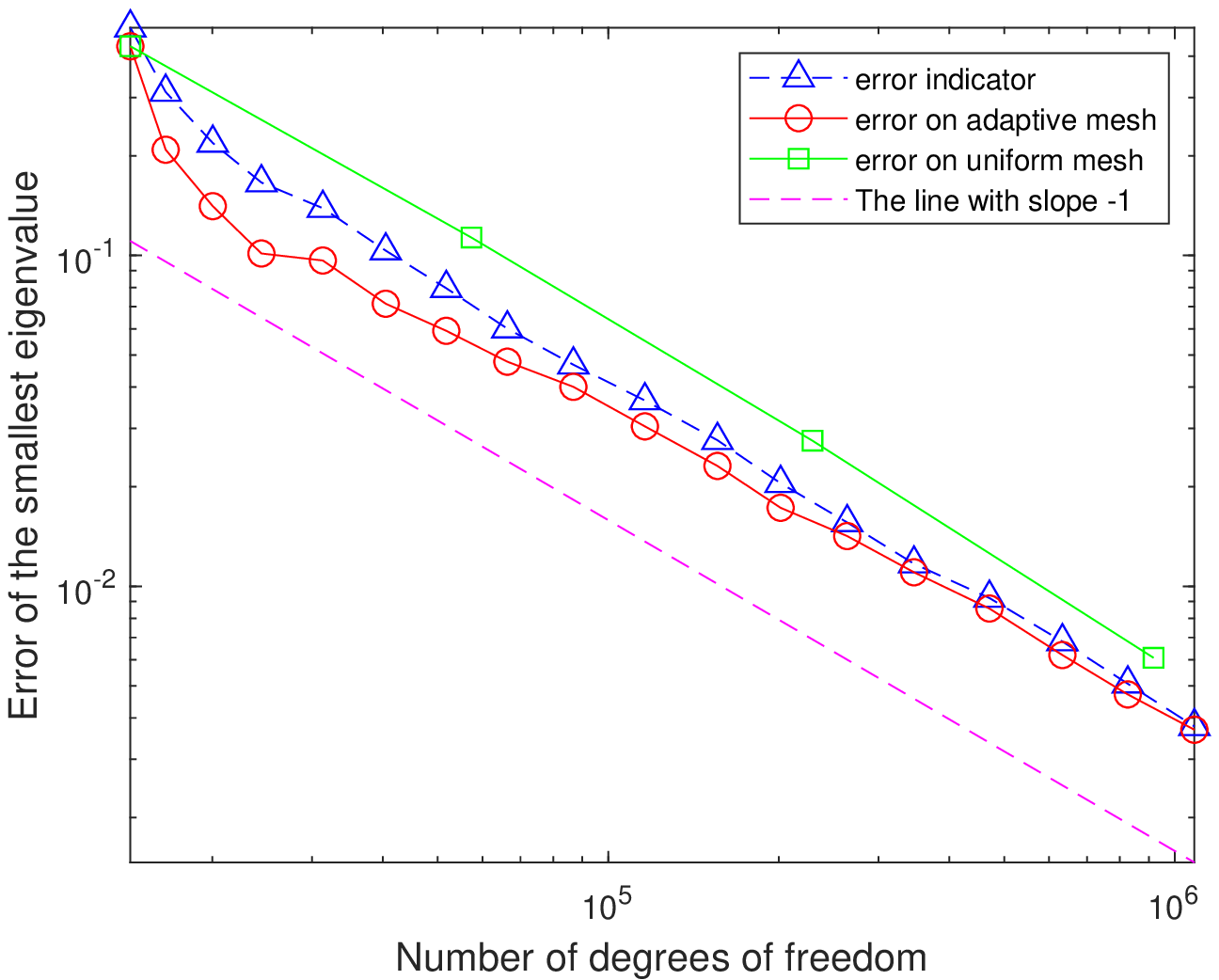}
\end{minipage}
\caption{Adaptive mesh after $l$=7 refinement times (left) and error curves (right) of the smallest eigenvalue by DGFEM using $\mathbb{P}_{1}-\mathbb{P}_{0}$ element  on $\Omega_{C}$.}
\end{figure}
\begin{figure}[H]
\begin{minipage}[t]{0.5\textwidth}
  \centering
   \includegraphics[width=7.3cm]{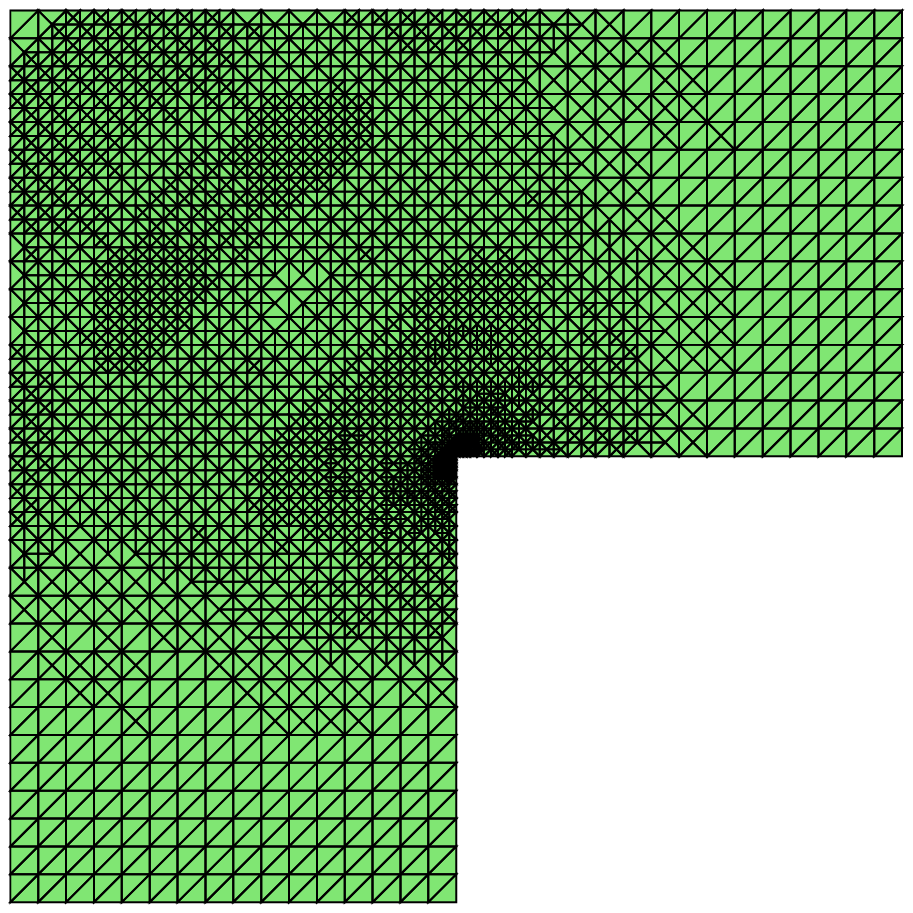}
\end{minipage}
\begin{minipage}[t]{0.5\textwidth}
  \centering
  \includegraphics[width=7.3cm]{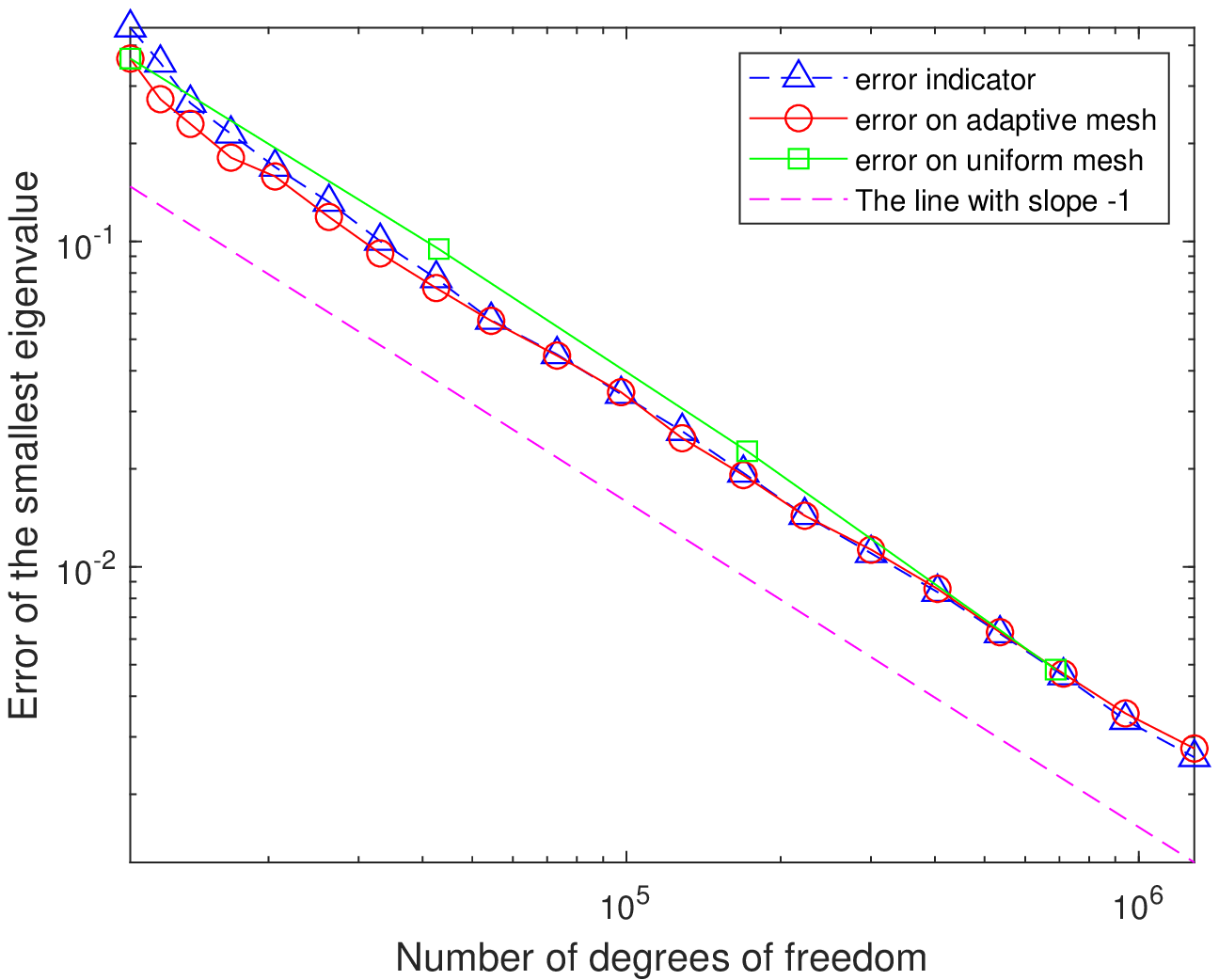}
\end{minipage}
\caption{Adaptive mesh after $l$=7 refinement times (left) and the error curves (right) of the smallest eigenvalue
by DGFEM using $\mathbb{P}_{1}-\mathbb{P}_{0}$ element on $\Omega_{L}$.}
\end{figure}
\begin{figure}[H]
\begin{minipage}[t]{0.5\textwidth}
  \centering
  \includegraphics[width=7.3cm]{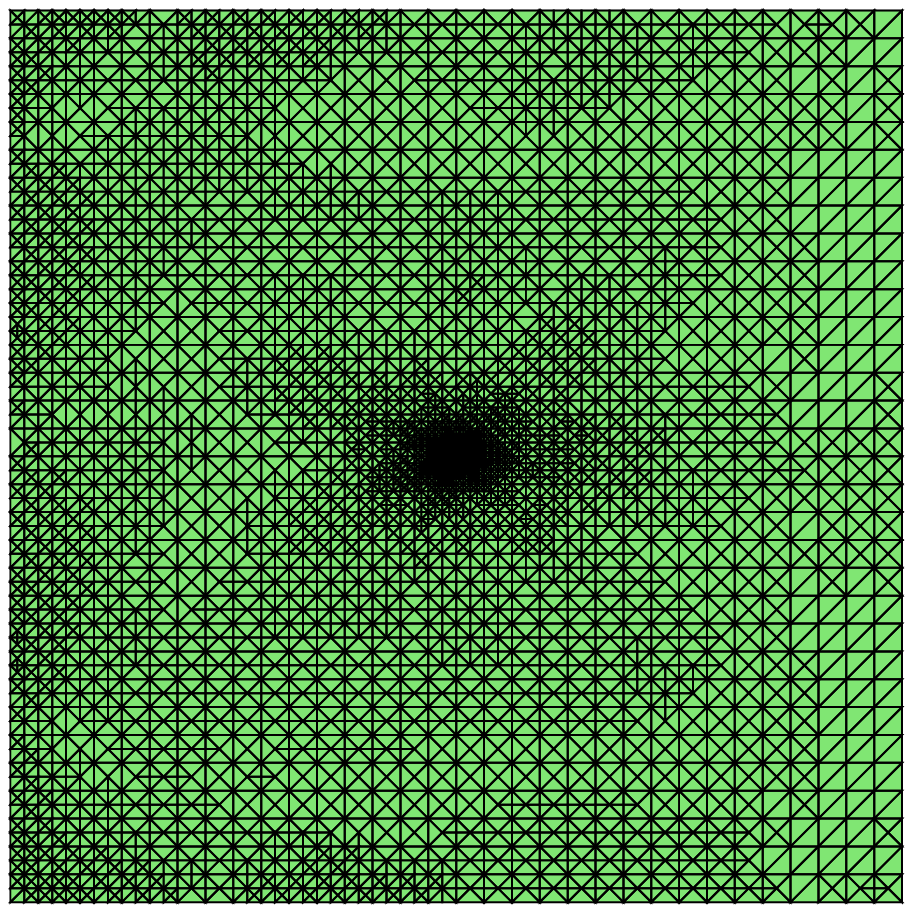}
\end{minipage}
\begin{minipage}[t]{0.5\textwidth}
  \centering
 \includegraphics[width=7.3cm]{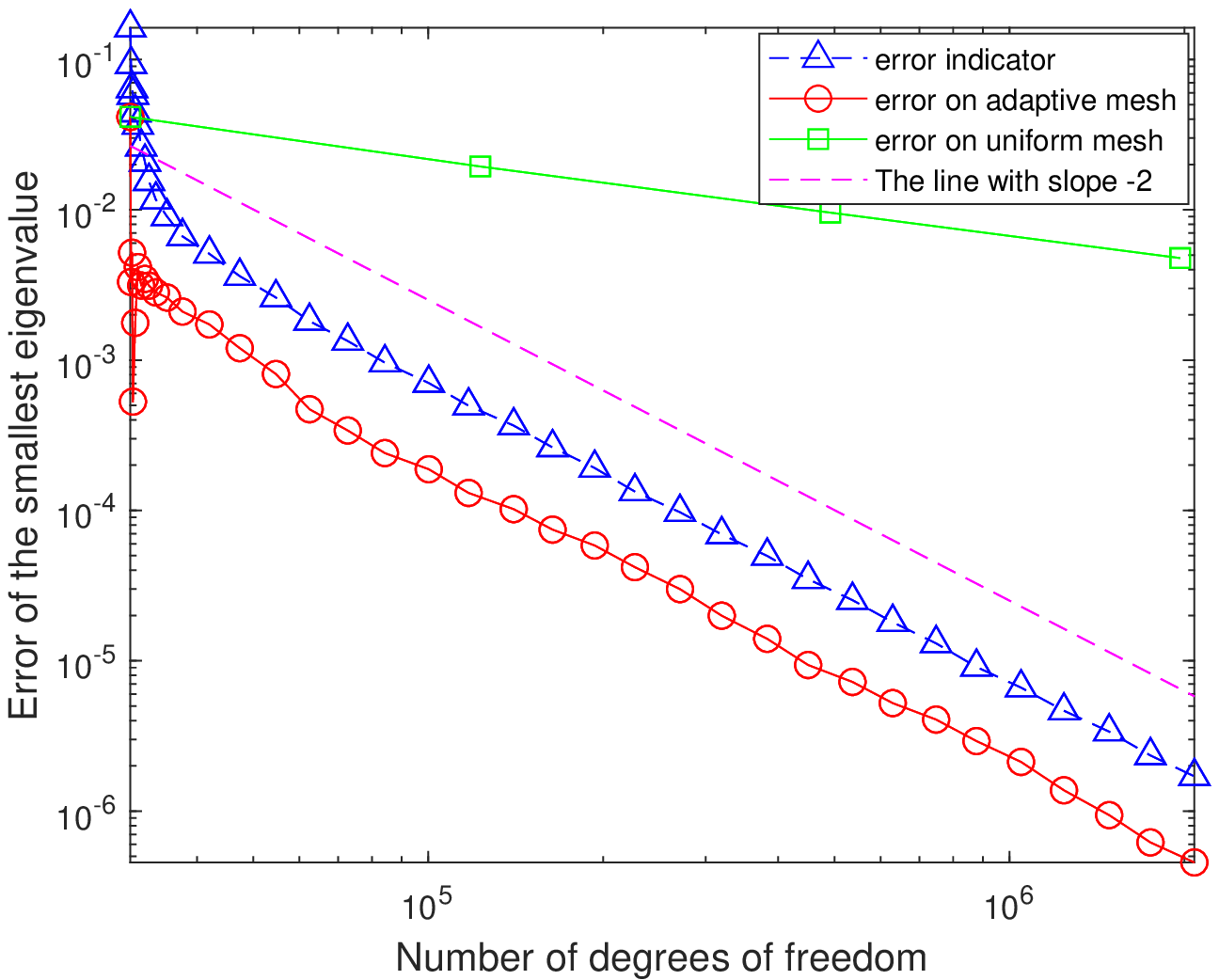}
\end{minipage}
\caption{Adaptive mesh after $l$=15 refinement times (left) and the error curves (right) of the smallest eigenvalue
by DGFEM using $\mathbb{P}_{2}-\mathbb{P}_{1}$ element on $\Omega_{C}$.}
\end{figure}
\begin{figure}[H]
\begin{minipage}[t]{0.5\textwidth}
  \centering
  \includegraphics[width=7.3cm]{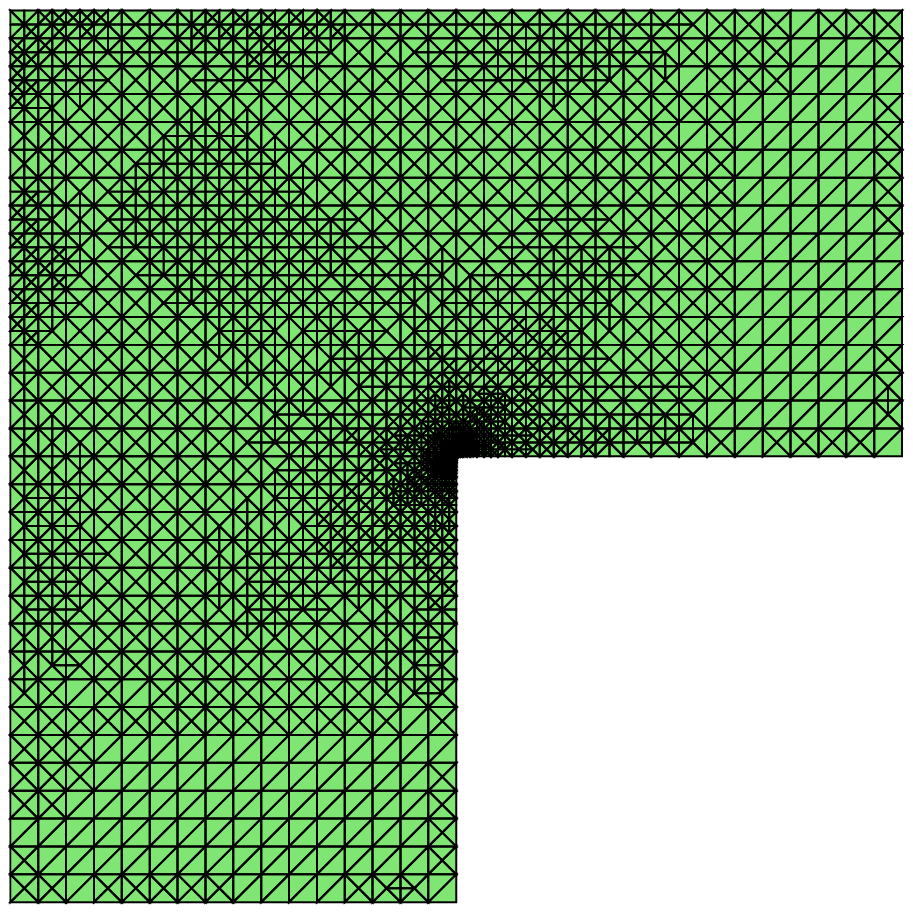}
\end{minipage}
\begin{minipage}[t]{0.5\textwidth}
  \centering
  \includegraphics[width=7.3cm]{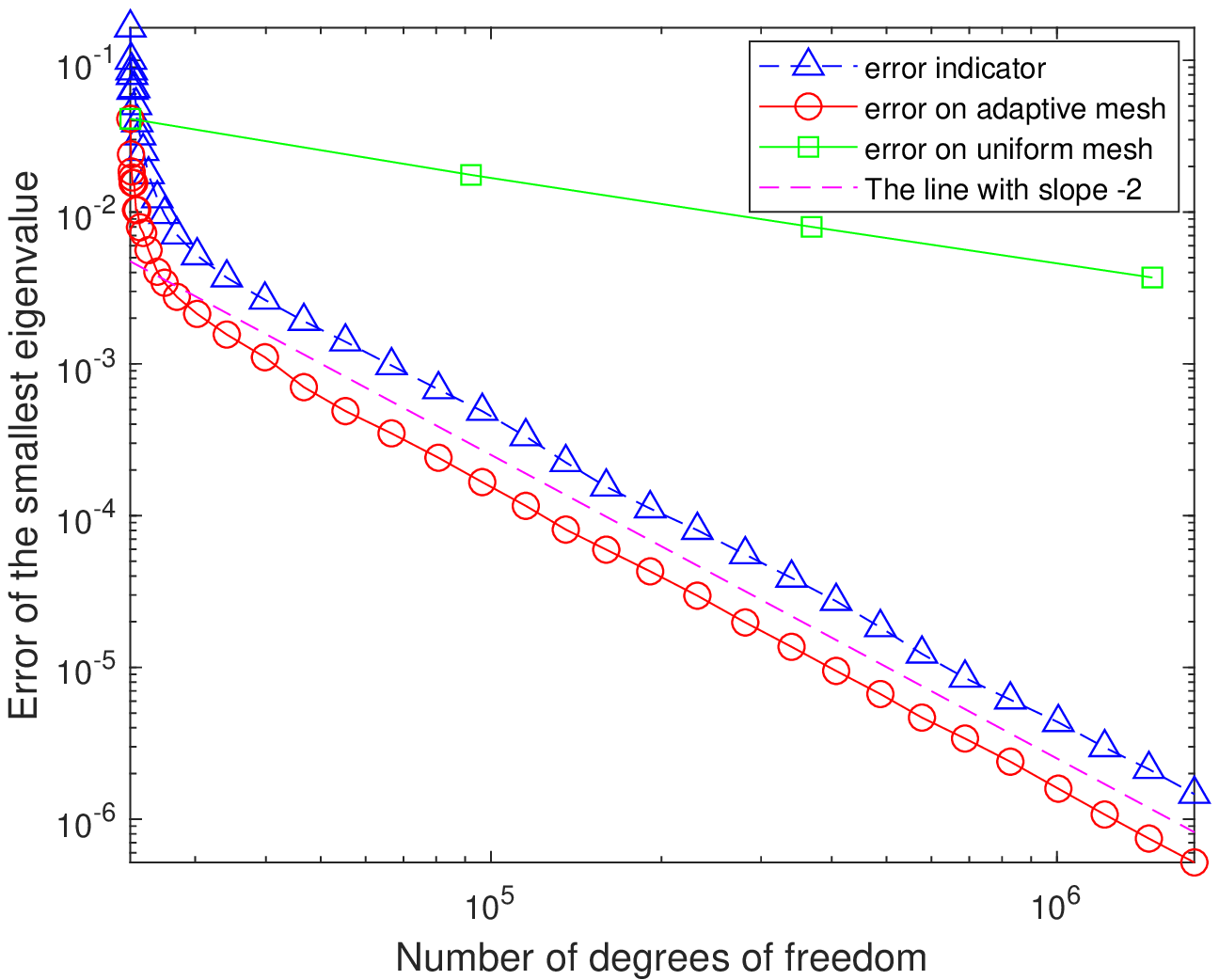}
\end{minipage}
\caption{Adaptive mesh after $l$=15 refinement times (left) and error curves (right) of the smallest eigenvalue
by DGFEM using $\mathbb{P}_{2}-\mathbb{P}_{1}$ element on $\Omega_{L}$.}
\end{figure}
\begin{figure}[H]
\begin{minipage}[t]{0.5\textwidth}
  \centering
  \includegraphics[width=7.3cm]{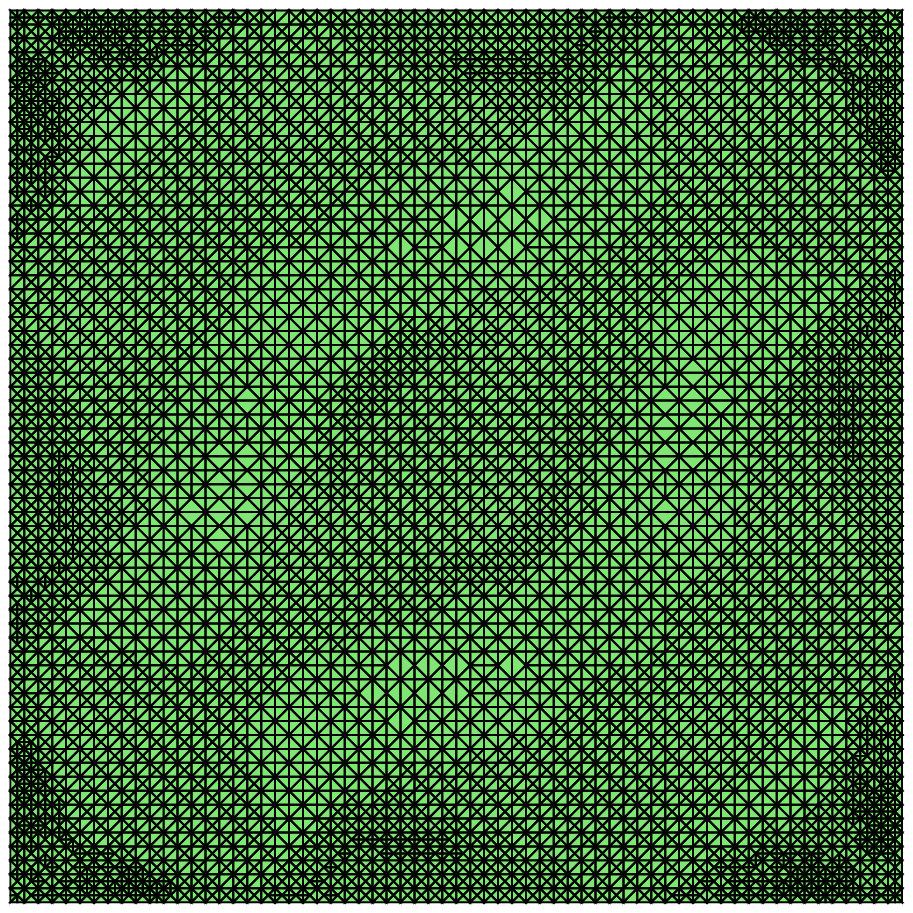}
\end{minipage}
\begin{minipage}[t]{0.5\textwidth}
  \centering
  \includegraphics[width=7.3cm]{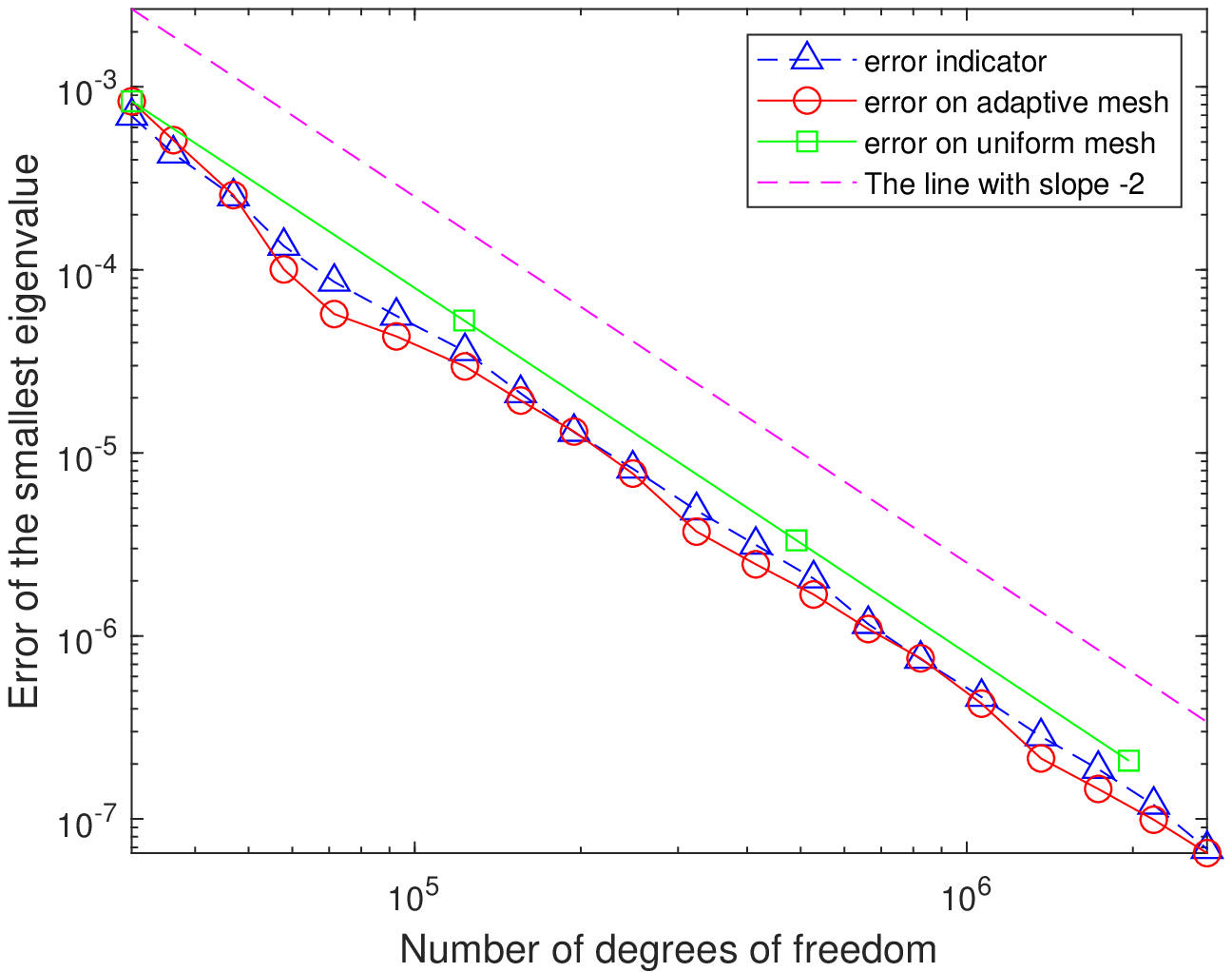}
\end{minipage}
\caption{Adaptive mesh after $l$=7 refinement times (left) and the error curves (right) of the smallest eigenvalue
by DGFEM using $\mathbb{P}_{2}-\mathbb{P}_{1}$ element on $\Omega_{S}$.}
\end{figure}
\begin{figure}[H]
\begin{minipage}[t]{0.5\textwidth}
  \centering
  \includegraphics[width=7.3cm]{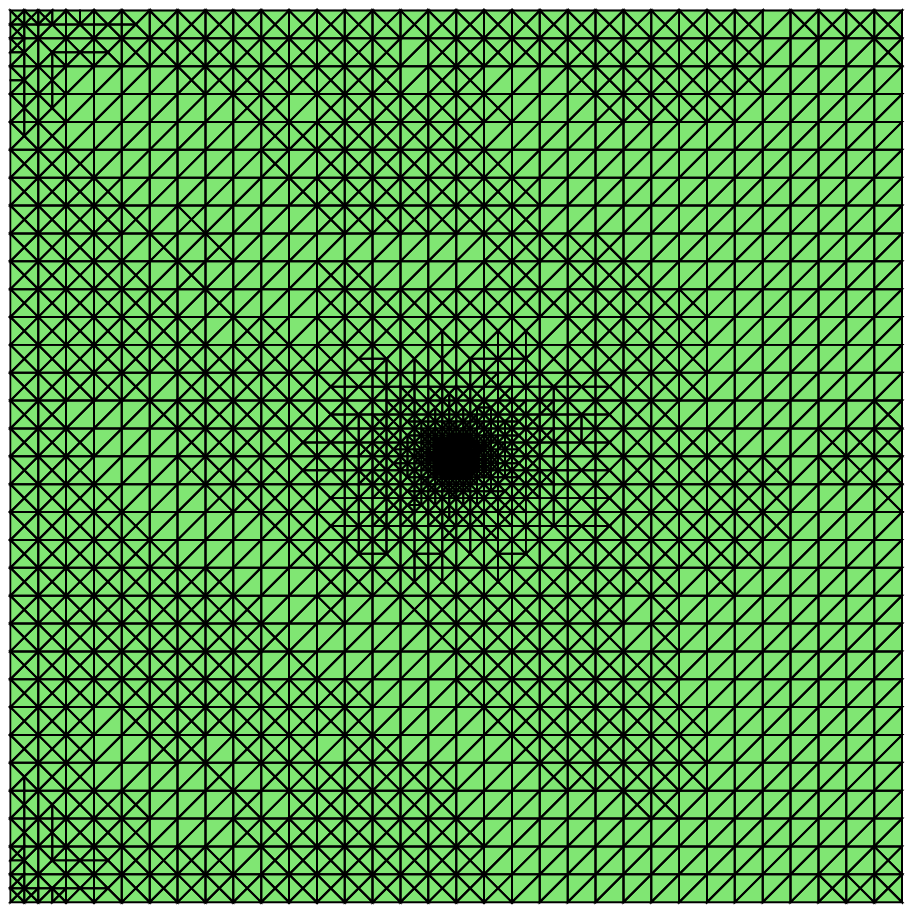}
\end{minipage}
\begin{minipage}[t]{0.5\textwidth}
  \centering
  \includegraphics[width=7.3cm]{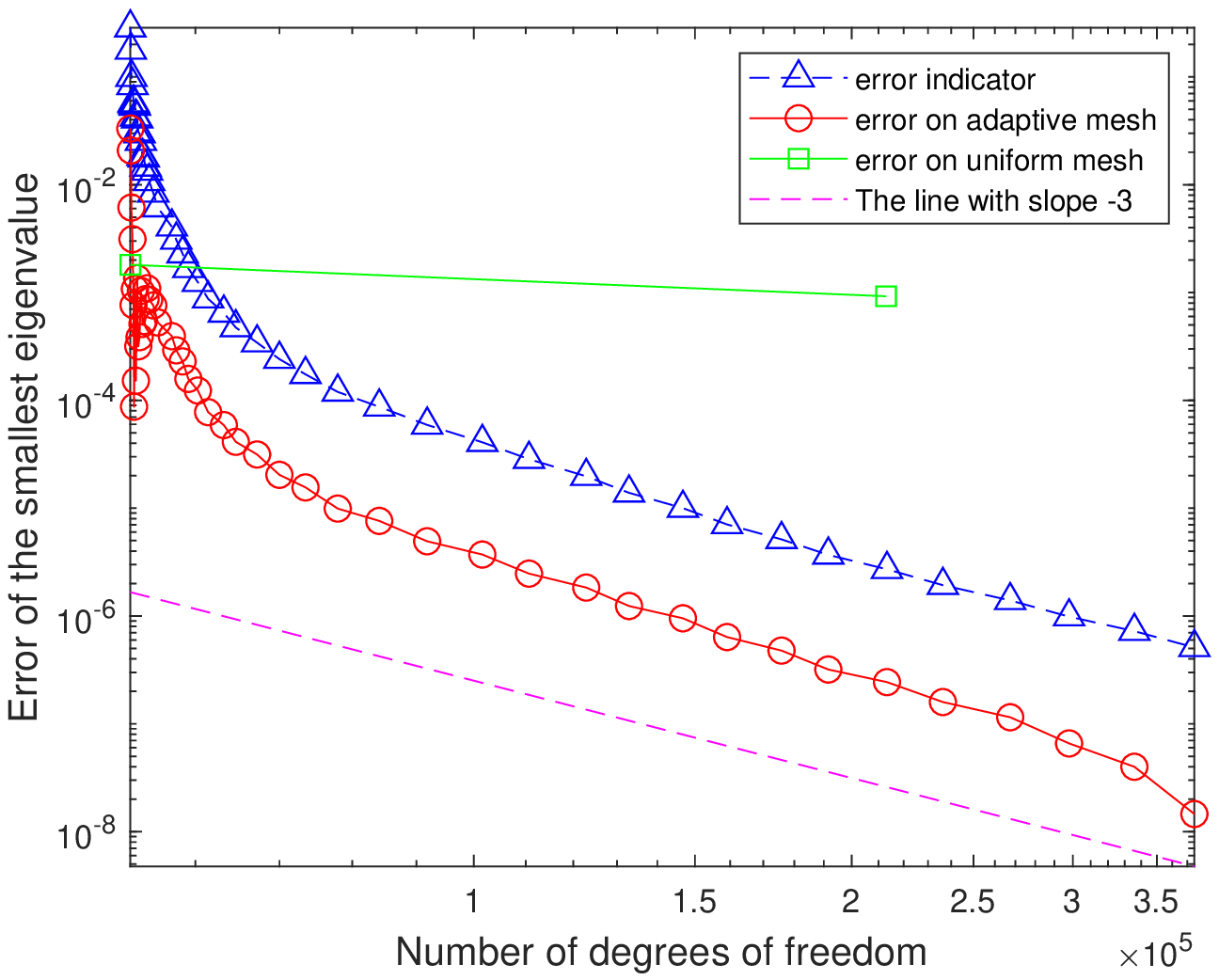}
\end{minipage}
\caption{Adaptive mesh after $l$=25 refinement times (left) and the error curves (right) of the smallest eigenvalue
by DGFEM using $\mathbb{P}_{3}-\mathbb{P}_{2}$ element on $\Omega_{C}$ }
\end{figure}
\begin{figure}[H]
\begin{minipage}[t]{0.5\textwidth}
  \centering
  \includegraphics[width=7.3cm]{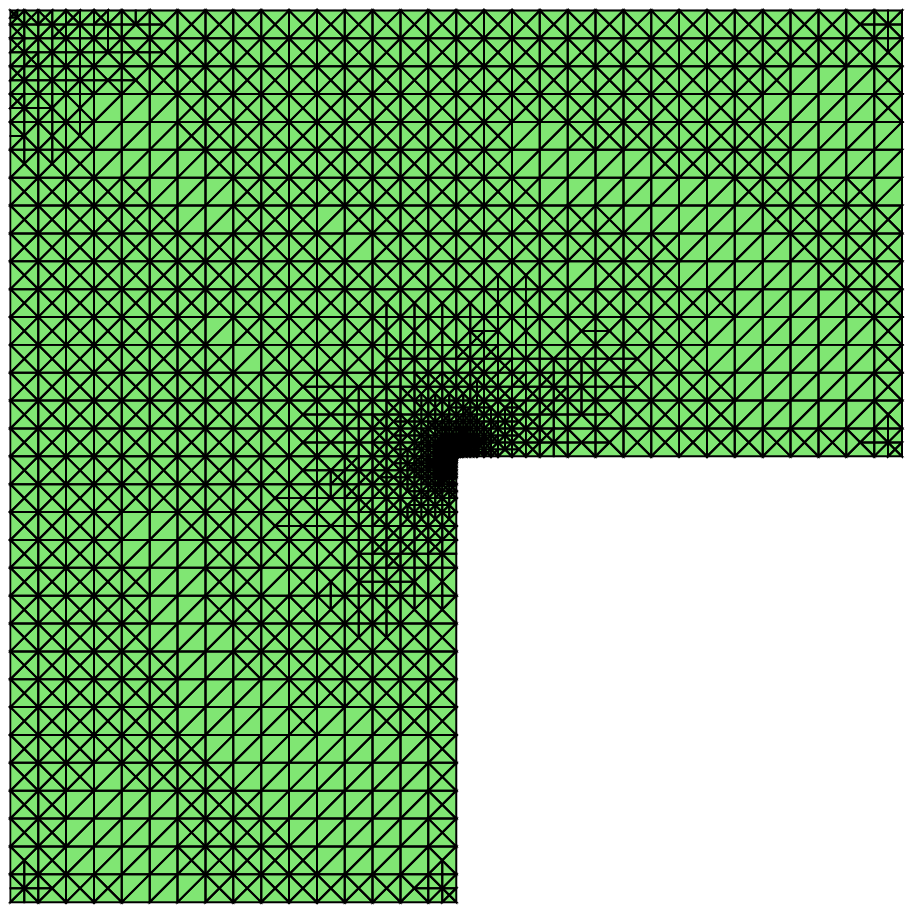}
\end{minipage}
\begin{minipage}[t]{0.5\textwidth}
  \centering
  \includegraphics[width=7.3cm]{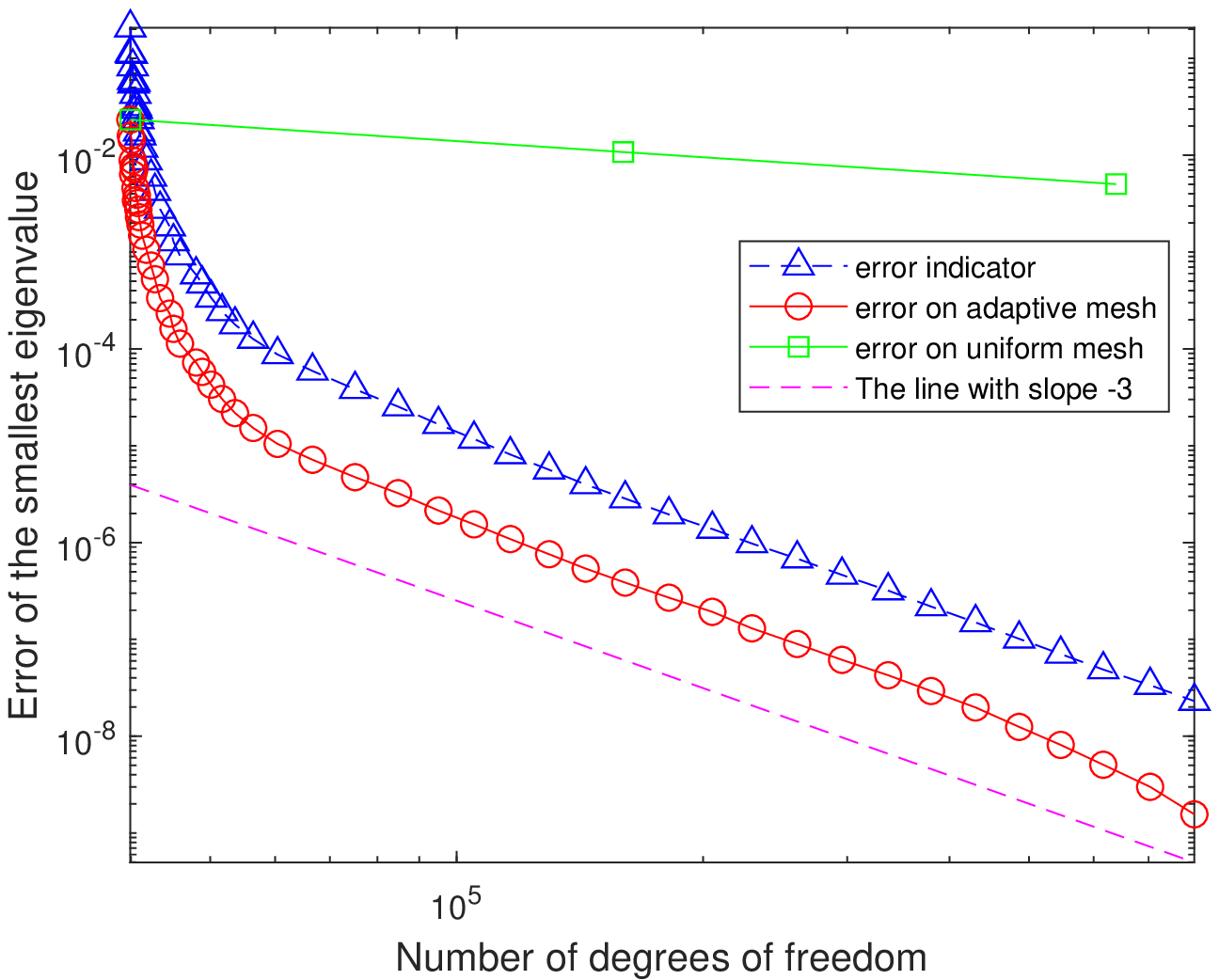}
\end{minipage}
\caption{Adaptive mesh after $l$=25 refinement times (left) and the error curves (right) of the smallest eigenvalue
by DGFEM using $\mathbb{P}_{3}-\mathbb{P}_{2}$ element on $\Omega_{L}$ .}
\end{figure}
\begin{figure}[H]
\begin{minipage}[t]{0.5\textwidth}
  \centering
  \includegraphics[width=7.3cm]{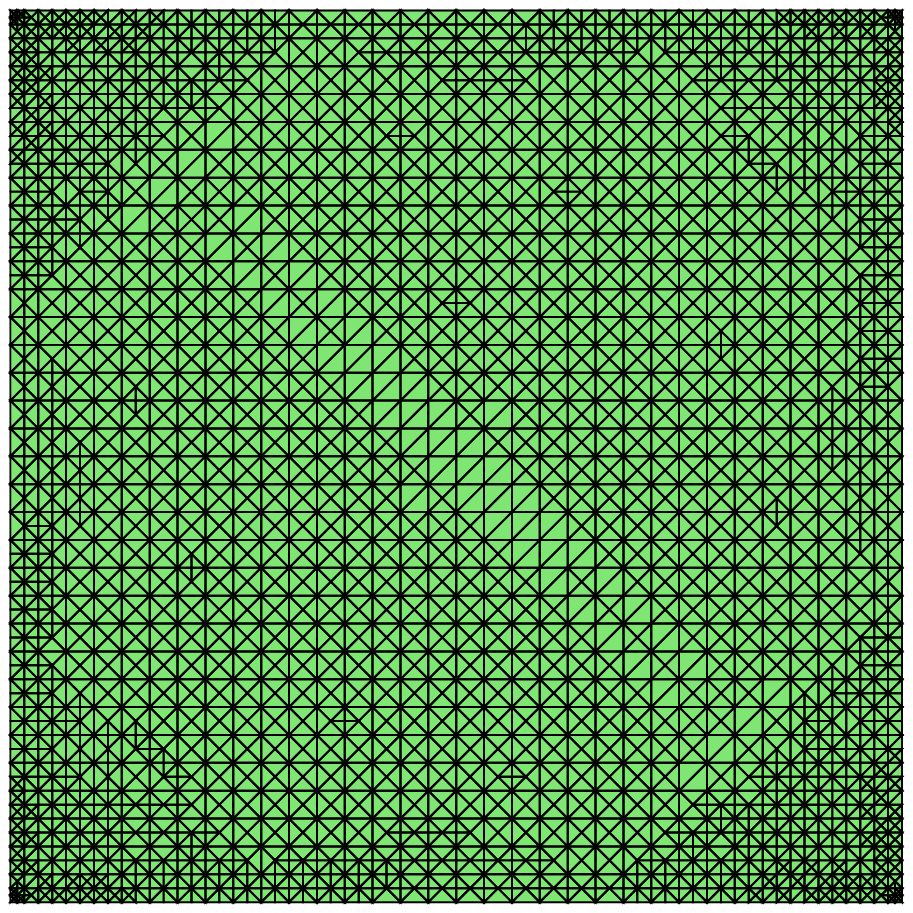}
\end{minipage}
\begin{minipage}[t]{0.5\textwidth}
  \centering
 \includegraphics[width=7.3cm]{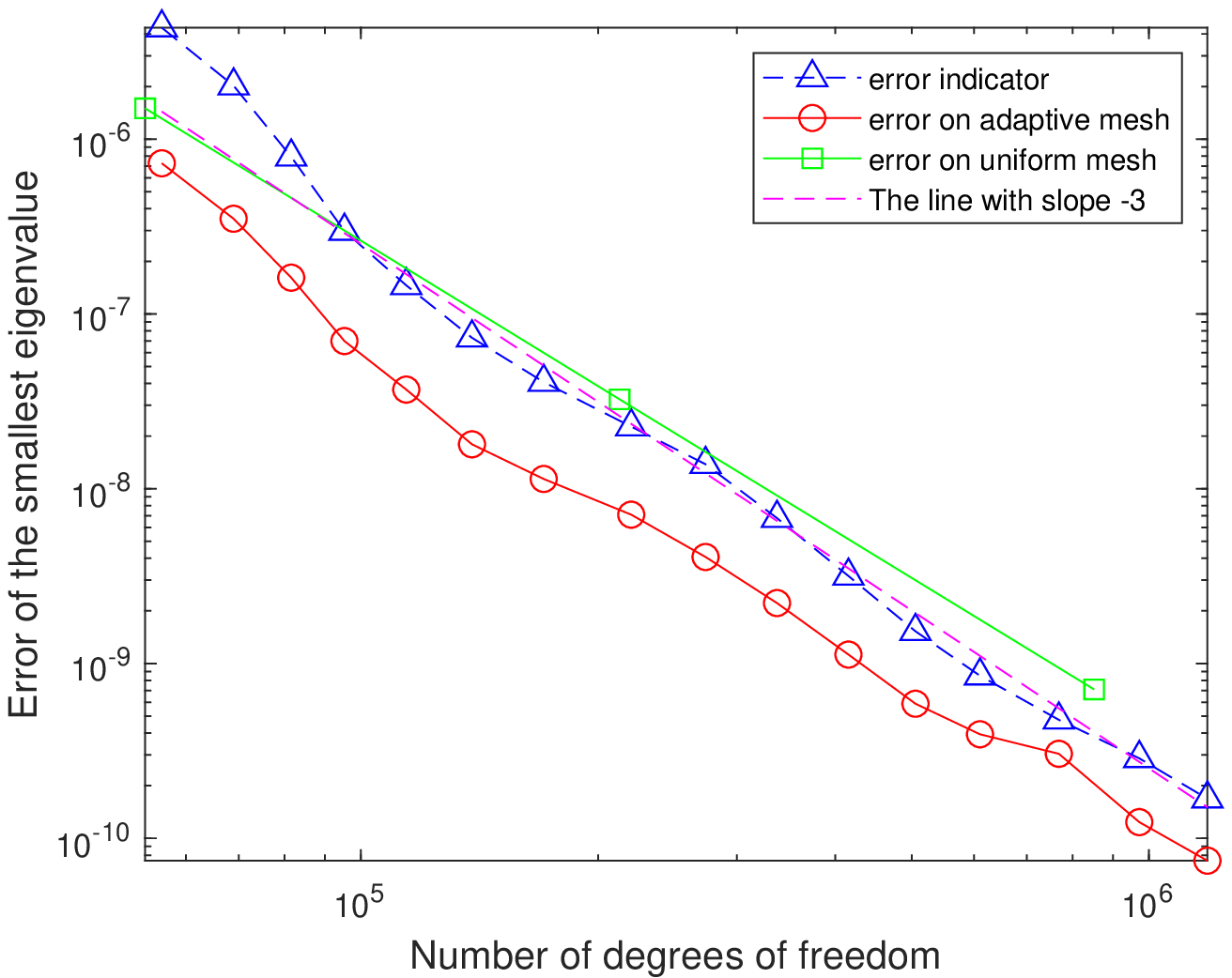}
\end{minipage}
\caption{Adaptive mesh after $l$=5 refinement times (left) and the error curves (right) of the smallest eigenvalue
by DGFEM using $\mathbb{P}_{3}-\mathbb{P}_{2}$ element on $\Omega_{S}$. }
\end{figure}

\begin{table}[H]
\caption{ \ The smallest eigenvalue using $\mathbb{P}_{3}-\mathbb{P}_{2}$ element on adaptive mesh on $\Omega_{C}$. }
\label{tab:1}
\begin{center}
\begin{tabular}{{cccccccccccccc}}
  \hline
        $l$ & $dof$ & $\lambda_{1,h_{l}}$  & $l$ & $dof$ & $\lambda_{1,h_{l}}$   \\ \hline
        1 & 53248 & 29.950023991  & 26 & 63206 & 29.916921865   \\
        2 & 53300 & 29.937640600  & 27 & 64610 & 29.916904165   \\
        5 & 53560 & 29.917626784  & 30 & 73424 & 29.916878484   \\
        10 & 54028 & 29.917180037  & 35 & 110630 & 29.916865373   \\
        11 & 54158 & 29.917248497  & 36 & 122876 & 29.916864735   \\
        14 & 54574 & 29.917406356  & 39 & 159094 & 29.916863538   \\
        15 & 54756 & 29.917731006  & 40 & 175812 & 29.916863378   \\
        21 & 57928 & 29.917153931  & 46 & 335842 & 29.916862940   \\
        22 & 58604 & 29.917092477  & 47 & 374920 & 29.916862915   \\
        23 & 59228 & 29.917021023  & 48 & 424814 & 29.916862902   \\
        24 & 60268 & 29.916986202  & 49 & 475852 & 29.916862889   \\
        25 & 61412 & 29.916940636  & 50 & 537862 & 29.916862882   \\
\hline
\end{tabular}\end{center}
\end{table}

\begin{table}[H]
\caption{ \ The smallest eigenvalue using $\mathbb{P}_{3}-\mathbb{P}_{2}$ element on adaptive mesh on $\Omega_{L}$. }
\label{tab:1}
\begin{center}
\begin{tabular}{{cccccccccccccc}}
  \hline
        $l$ & $dof$ & $\lambda_{1,h_{l}}$  & $l$ & $dof$ & $\lambda_{1,h_{l}}$   \\ \hline
        1 & 39936 & 32.155997914  & 27 & 53612 & 32.132716405   \\
        2 & 39988 & 32.148565928  & 28 & 56420 & 32.132709908   \\
        3 & 40092 & 32.147074075  & 29 & 60424 & 32.132705093   \\
        4 & 40196 & 32.141536961  & 30 & 66664 & 32.132701814   \\
        5 & 40248 & 32.139031080  & 31 & 75140 & 32.132699385   \\
        13 & 40924 & 32.134988686  & 39 & 181662 & 32.132694920   \\
        14 & 41080 & 32.134620945  & 40 & 205244 & 32.132694843   \\
        15 & 41288 & 32.134171324  & 41 & 229424 & 32.132694780   \\
        23 & 48048 & 32.132766985  & 49 & 616304 & 32.132694655   \\
        24 & 48880 & 32.132752576  & 50 & 703092 & 32.132694653   \\
        25 & 50128 & 32.132737367  & 51 & 796276 & 32.132694652   \\
        26 & 51688 & 32.132725042  & ~ & ~ &   \\
\hline
\end{tabular}\end{center}
\end{table}

\begin{table}[H]
\caption{ \ The smallest eigenvalue using $\mathbb{P}_{3}-\mathbb{P}_{2}$ element on adaptive mesh on $\Omega_{S}$. }
\label{tab:1}
\begin{center}
\begin{tabular}{{cccccccccccccc}}
  \hline
        $l$ & $dof$ & $\lambda_{1,h_{l}}$  & $l$ & $dof$ & $\lambda_{1,h_{l}}$   \\ \hline
        1 & 53248 & 52.3446926681  & 10 & 273780 & 52.3446911721   \\
        2 & 55900 & 52.3446918954  & 11 & 337324 & 52.3446911702   \\
        3 & 68952 & 52.3446915184  & 12 & 415636 & 52.3446911691   \\
        4 & 81588 & 52.3446913292  & 13 & 505544 & 52.3446911686   \\
        5 & 95316 & 52.3446912380  & 14 & 610376 & 52.3446911684   \\
        6 & 114192 & 52.3446912049  & 15 & 768560 & 52.3446911683   \\
        7 & 138372 & 52.3446911859  & 16 & 972192 & 52.3446911681   \\
        8 & 170612 & 52.3446911794  & 17 & 1186328 & 52.3446911679   \\
        9 & 220324 & 52.3446911751  & ~ & ~ &   \\
\hline
\end{tabular}\end{center}
\end{table}
\subsection{The results in three-dimensional domains}
~~~~~~~~~~~~~~~~~~~\\
\indent We also carry out numerical experiments on
 two three-dimensional domains:
$\Omega_{1}=(0,1)^{3}\setminus \{0\leq x\leq 0.5, 0\leq y\leq 0.5, 0.5\leq z\leq 1\}$ and $\Omega_{2}=(0,1)^{3}$.
In computation we take $\theta=0.25$ and initial mesh $\pi_{h_{0}}$ ($h_{0}=\frac{\sqrt{3}}{8}$ ).
To compute the error of the first eigenvalue for the Stokes eigenvalue problem,
we choose the values $\lambda_{1}=70.98560$  and $62.17341$ which are obtained by adaptive procedure with as much degrees of freedom as possible as the reference values for the domains $\Omega_{1}$ and $\Omega_{2}$, respectively.
%
~~~~~~~~~~~~~~~~~~~~~~~~~~~~~~~~\\
\indent The initial meshes are shown in Figures 9-11, and the adaptive refined meshes and the error curves are shown in Figures 12-17. The numerical results on adaptive mesh are listed in Tables 1-3.\\
\indent From Figures 12-17 we can see that the error curves and error indicators curves for DG methods using $\mathbb{P}_{k}-\mathbb{P}_{k-1}$ $(k=1,2,3)$ element are both approximately parallel to the line with slope $-\frac{2k}{3}$ , which indicates the error indicators
are reliable and efficient and the adaptive algorithm can reach the optimal convergence order. It coincides with our theoretical analysis.
It also can be seen from the error curves that under the same $dof$, the approximations obtained by high order elements are more accurate than those computed by low order elements on both uniform meshes and adaptive meshes.\\
\indent We can also see from Table 2 that it provides an upper bound for the exact eigenvalue by $\mathbb{P}_{k}-\mathbb{P}_{k-1}$ element. Note that the numerical results in Table 5 in \cite{Sun2022}  provide a lower bound on the cubic domain by the CR element. Thus we get a range for the exact eigenvalue of Stokes eigenvalue problem on the cubic.

\begin{figure}[H]
\begin{minipage}[t]{0.5\textwidth}
  \centering
  \includegraphics[width=7.3cm]{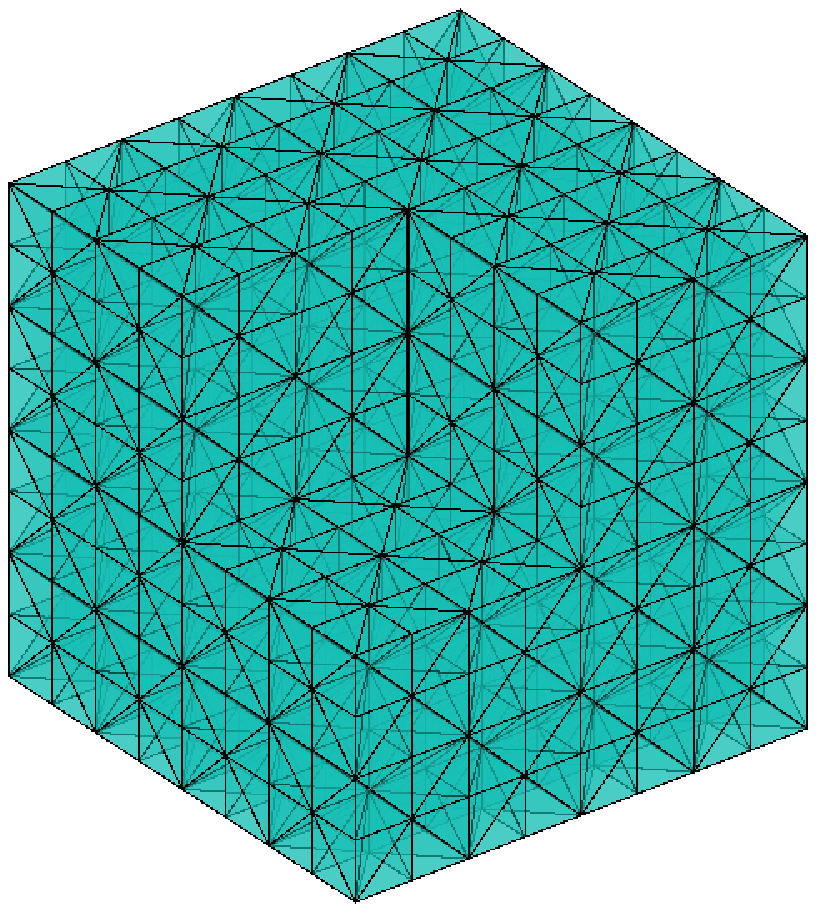}
\end{minipage}
\begin{minipage}[t]{0.5\textwidth}
  \centering
  \includegraphics[width=7.3cm]{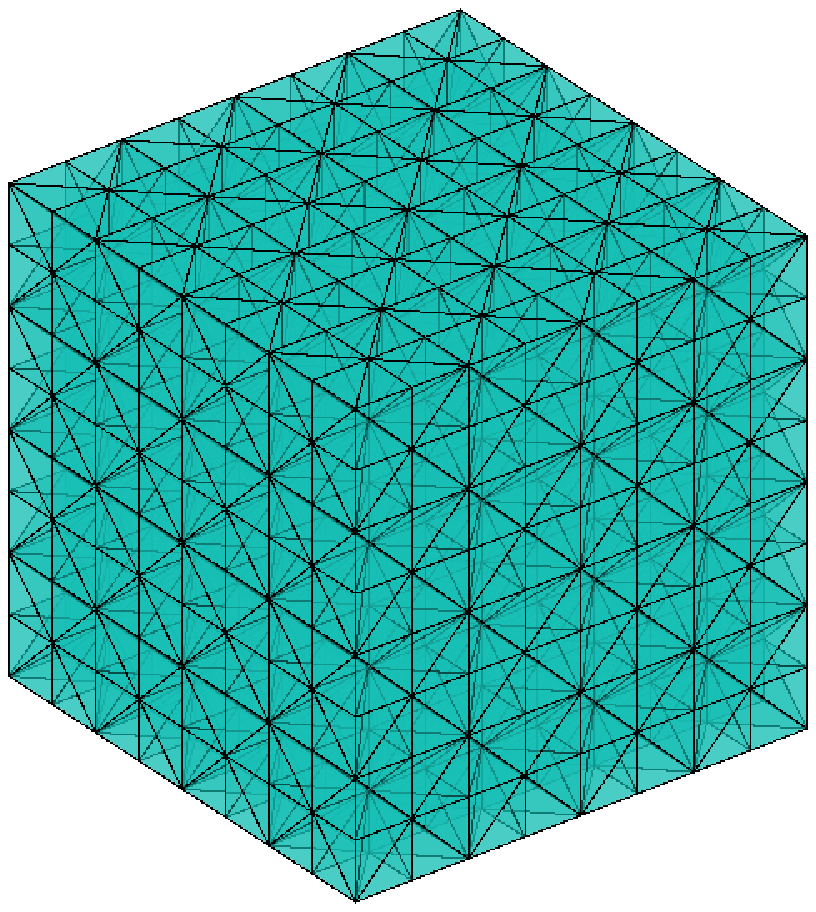}
\end{minipage}
\caption{The initial mesh on $\Omega_{1}$(left) and $~\Omega_{2}$(right). }
\end{figure}
\begin{figure}[H]
\begin{minipage}[t]{0.5\textwidth}
  \centering
    \includegraphics[width=7.3cm]{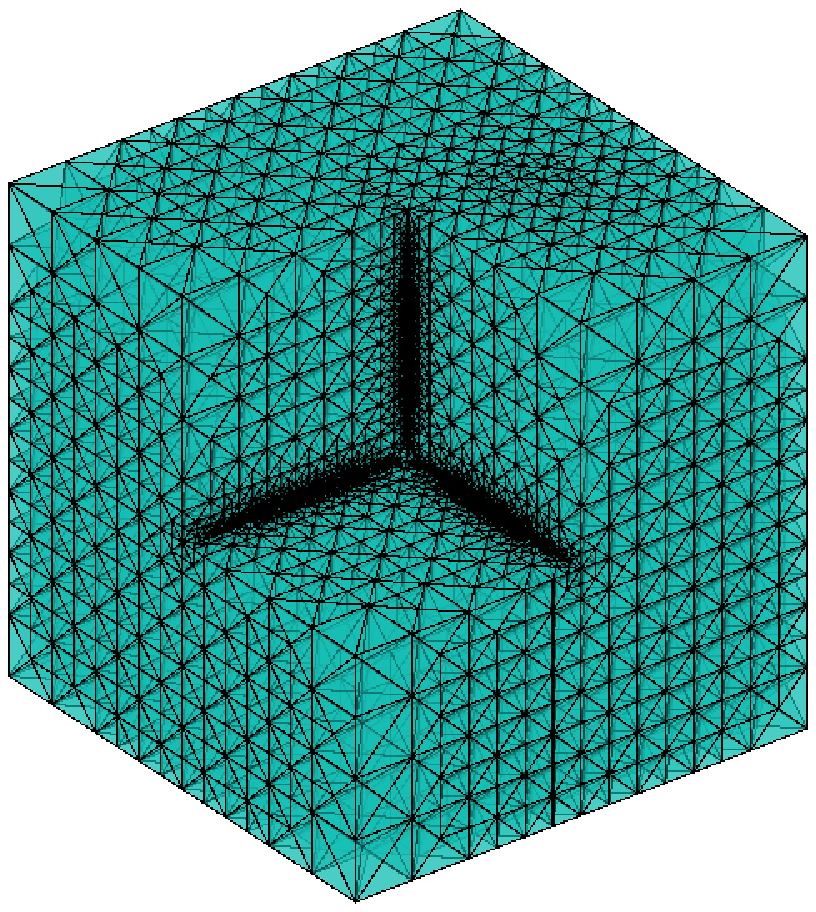}
\end{minipage}
\begin{minipage}[t]{0.5\textwidth}
  \centering
  \includegraphics[width=7.3cm]{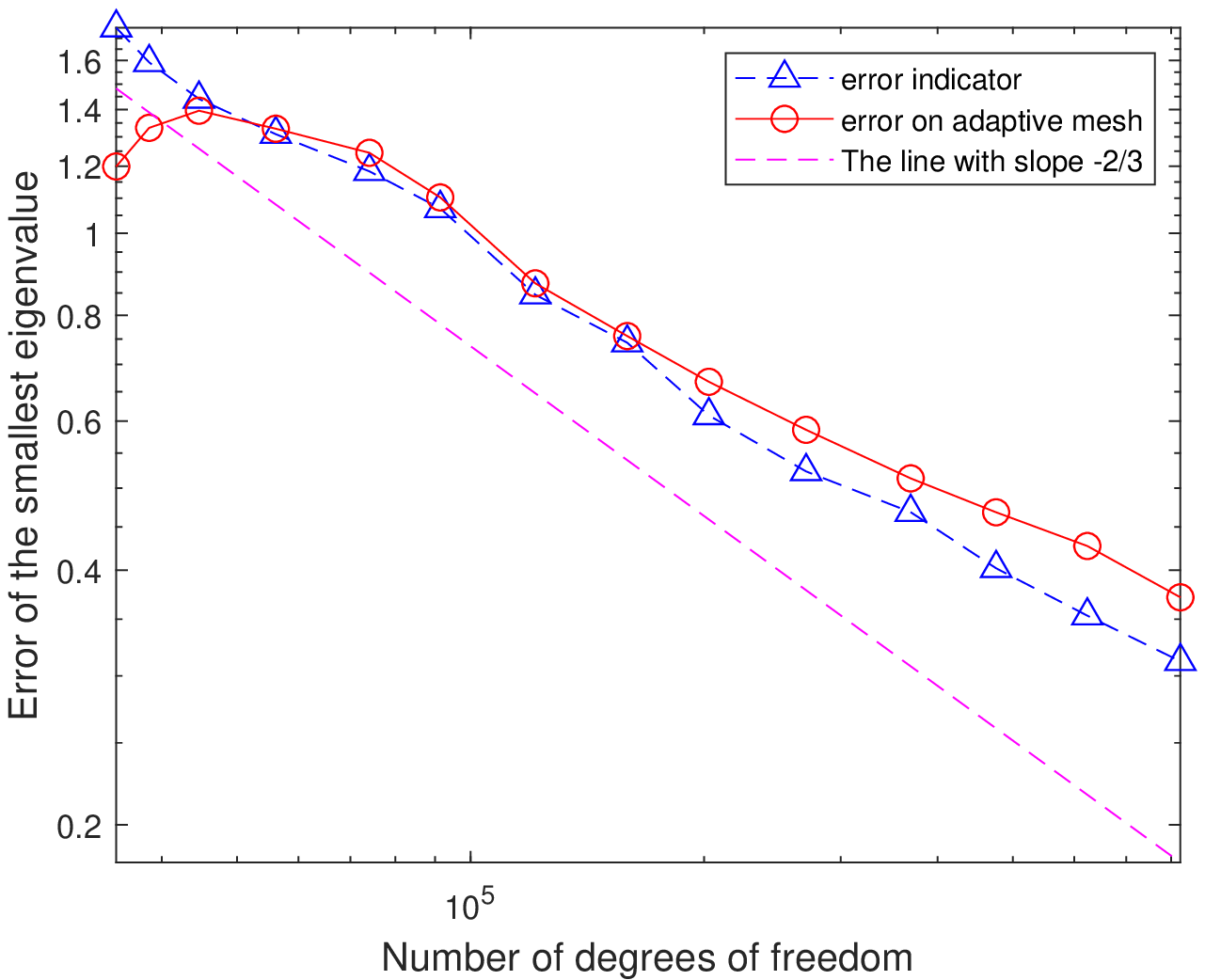}
\end{minipage}
\caption{Adaptive mesh after $l$=10 refinement times (left) and the error curves (right) of the smallest eigenvalue
by DGFEM using $\mathbb{P}_{1}-\mathbb{P}_{0}$ element on $\Omega_{1}$ }
\end{figure}

\begin{figure}[H]
\begin{minipage}[t]{0.5\textwidth}
  \centering
    \includegraphics[width=7.3cm]{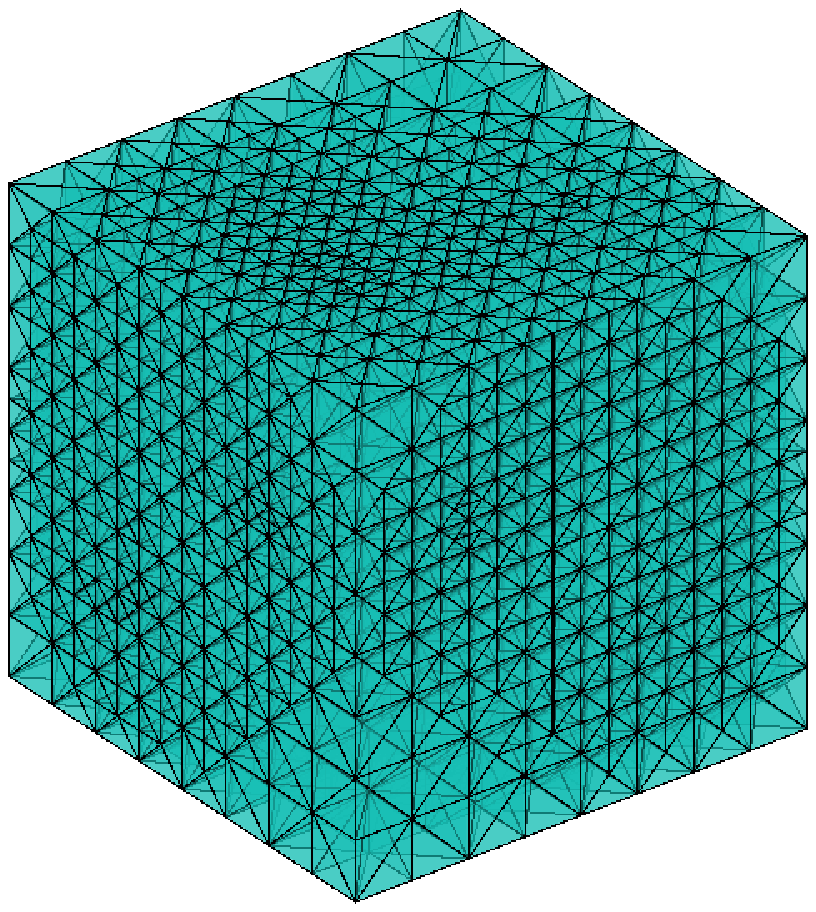}
\end{minipage}
\begin{minipage}[t]{0.5\textwidth}
  \centering
  \includegraphics[width=7.3cm]{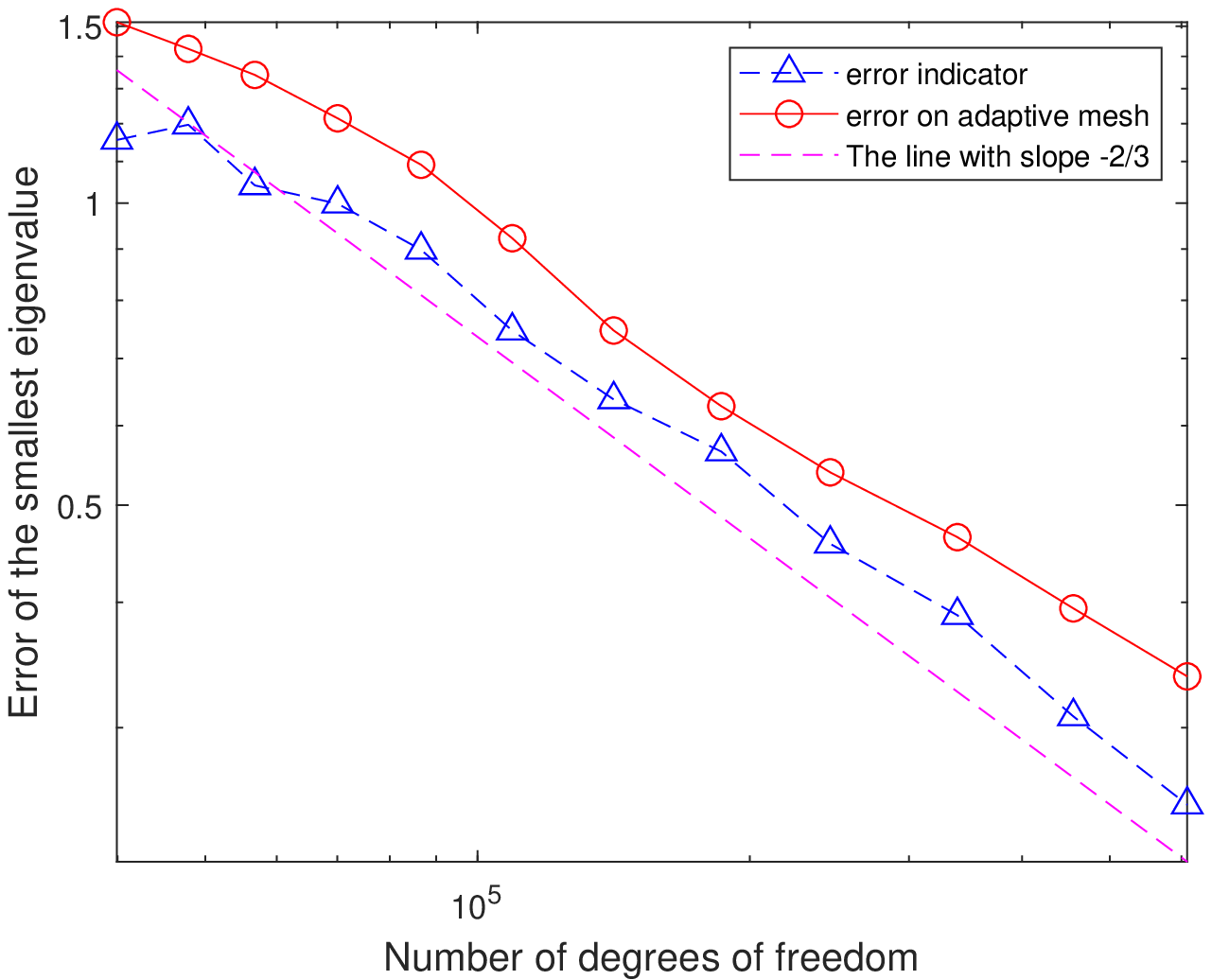}
\end{minipage}
\caption{Adaptive mesh after $l$=10 refinement times (left) and the error curves (right) of the smallest eigenvalue
by DGFEM using $\mathbb{P}_{1}-\mathbb{P}_{0}$ element on $\Omega_{2}$ }
\end{figure}
\begin{figure}[H]
\begin{minipage}[t]{0.5\textwidth}
  \centering
    \includegraphics[width=7.3cm]{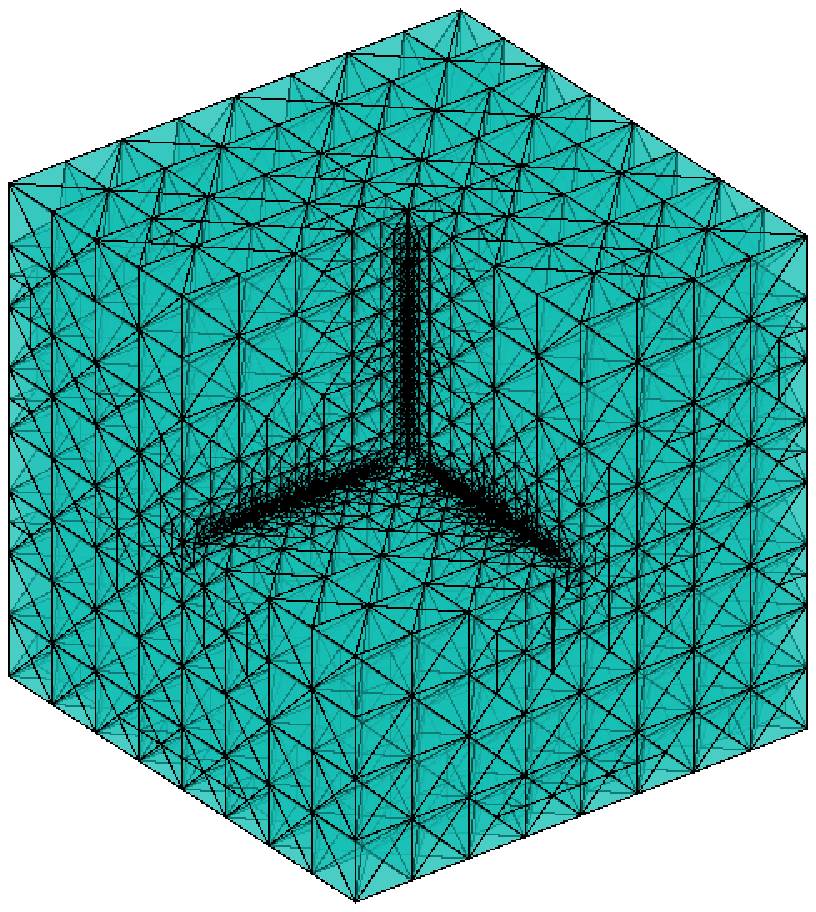}
\end{minipage}
\begin{minipage}[t]{0.5\textwidth}
  \centering
  \includegraphics[width=7.3cm]{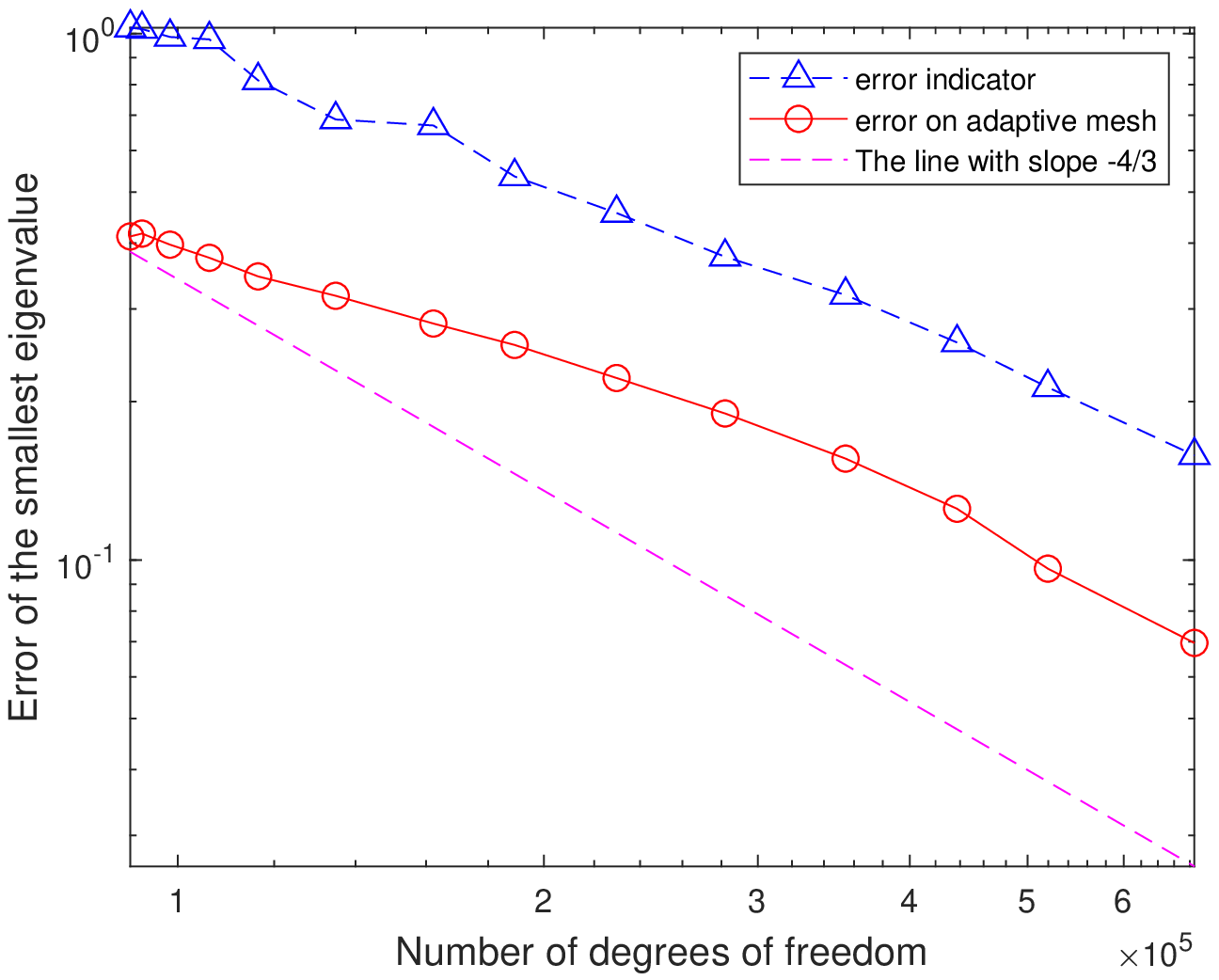}
\end{minipage}
\caption{Adaptive mesh after $l$=14 refinement times (left) and the error curves (right) of the smallest eigenvalue
by DGFEM using $\mathbb{P}_{2}-\mathbb{P}_{1}$ element on $\Omega_{1}$ }
\end{figure}

\begin{figure}[H]
\begin{minipage}[t]{0.5\textwidth}
  \centering
    \includegraphics[width=7.3cm]{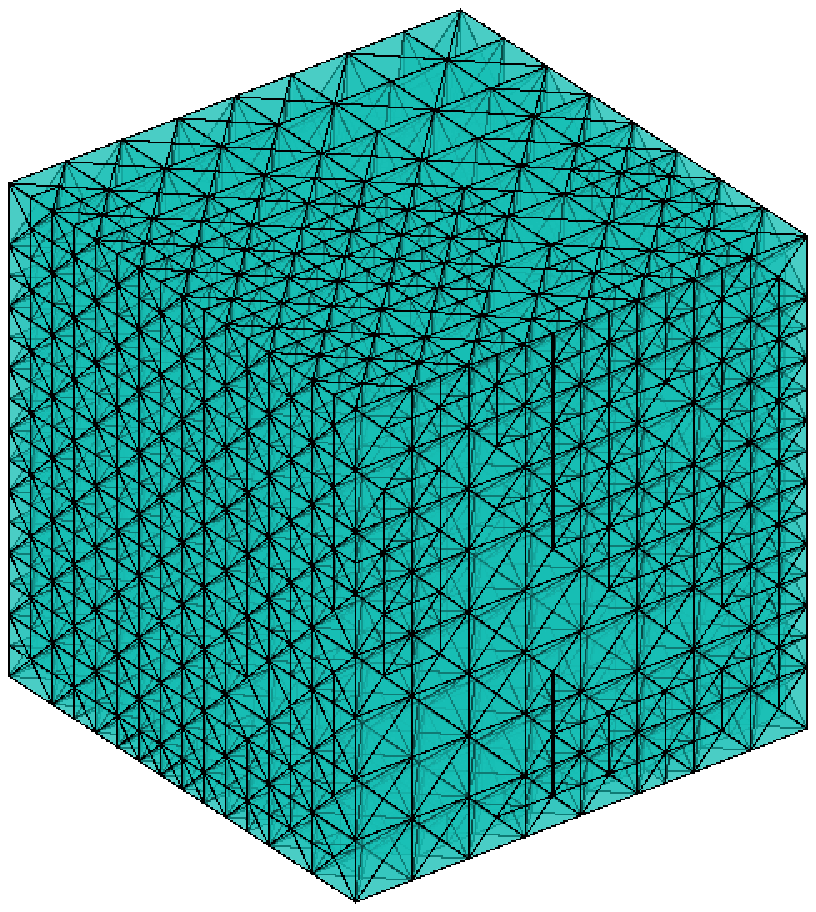}
\end{minipage}
\begin{minipage}[t]{0.5\textwidth}
  \centering
 \includegraphics[width=7.3cm]{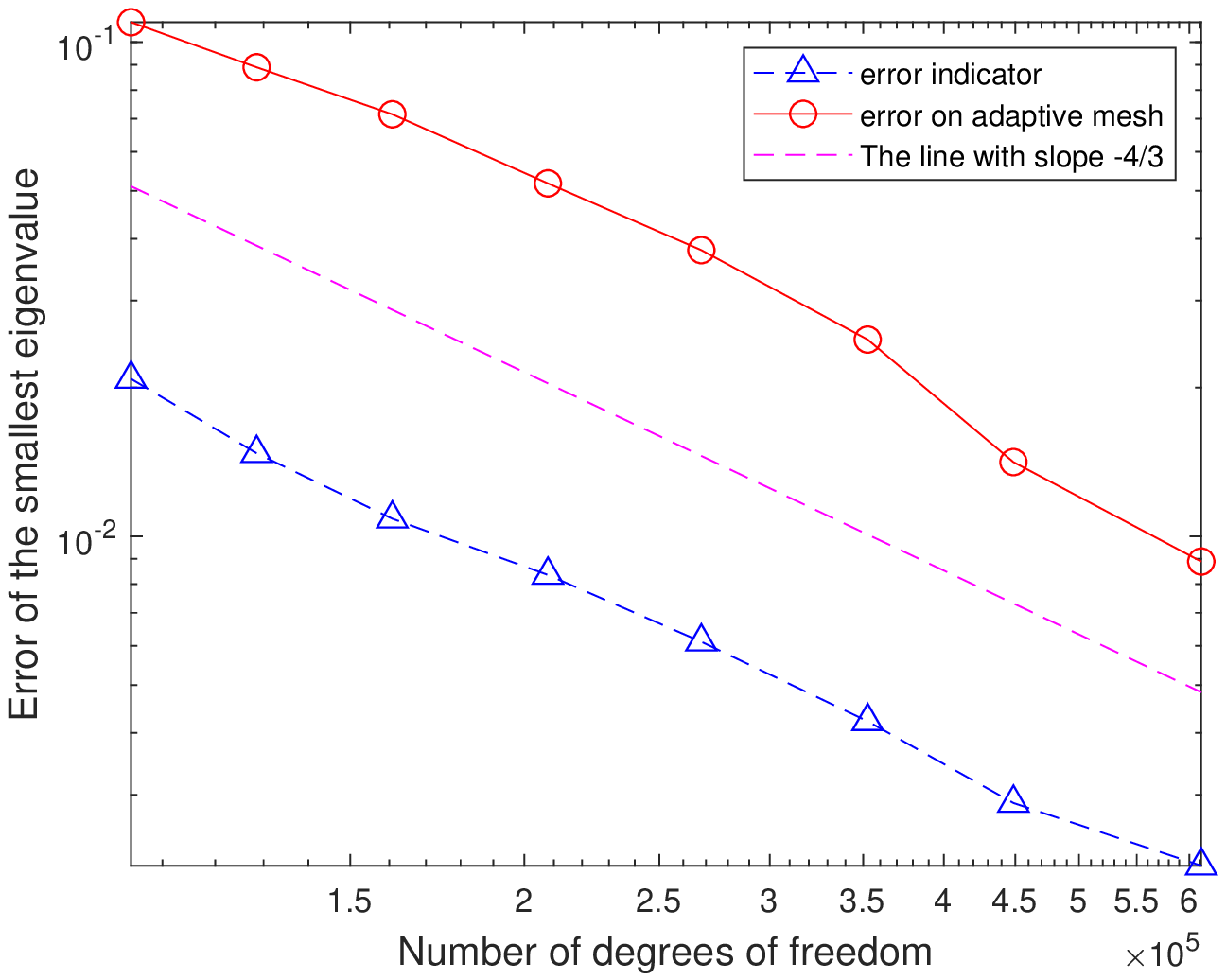}
\end{minipage}
\caption{Adaptive mesh after $l$=8 refinement times (left) and the error curves (right) of the smallest eigenvalue
by DGFEM using $\mathbb{P}_{2}-\mathbb{P}_{1}$ element on $\Omega_{2}$ }
\end{figure}
\begin{figure}[H]
\begin{minipage}[t]{0.5\textwidth}
  \centering
    \includegraphics[width=7.3cm]{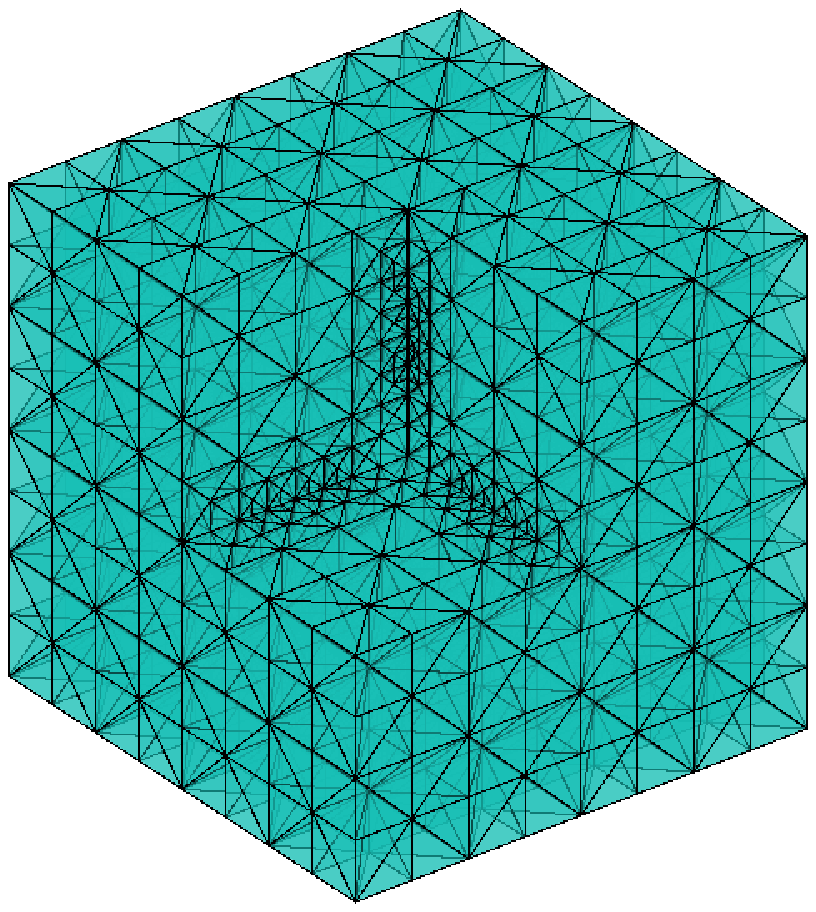}
\end{minipage}
\begin{minipage}[t]{0.5\textwidth}
  \centering
 \includegraphics[width=7.3cm]{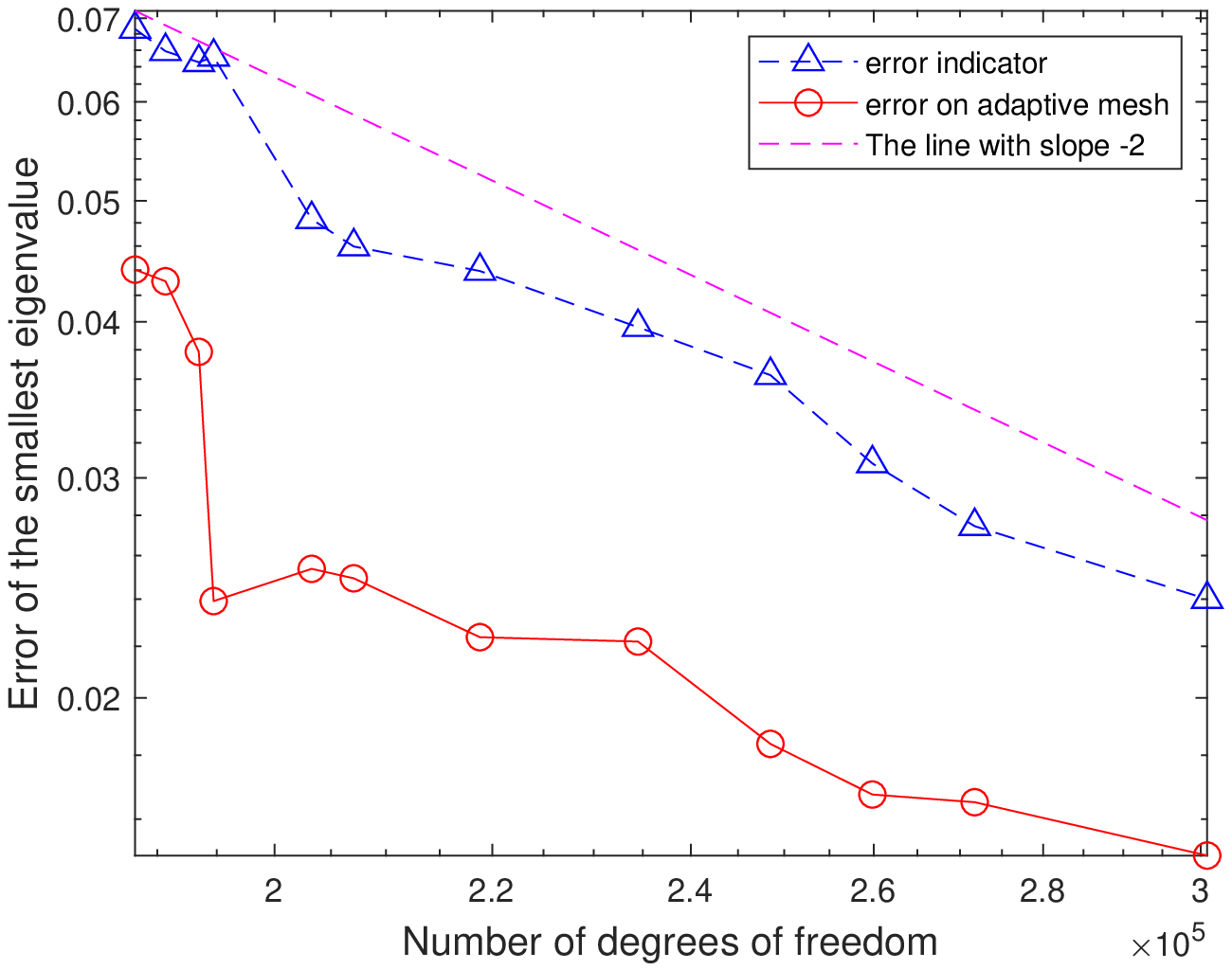}
\end{minipage}
\caption{Adaptive mesh after $l$=12 refinement times (left) and the error curves (right) of the smallest eigenvalue
by DGFEM using $\mathbb{P}_{3}-\mathbb{P}_{2}$ element on $\Omega_{1}$ }
\end{figure}

\begin{figure}[H]
\begin{minipage}[t]{0.5\textwidth}
  \centering
  \includegraphics[width=7.3cm]{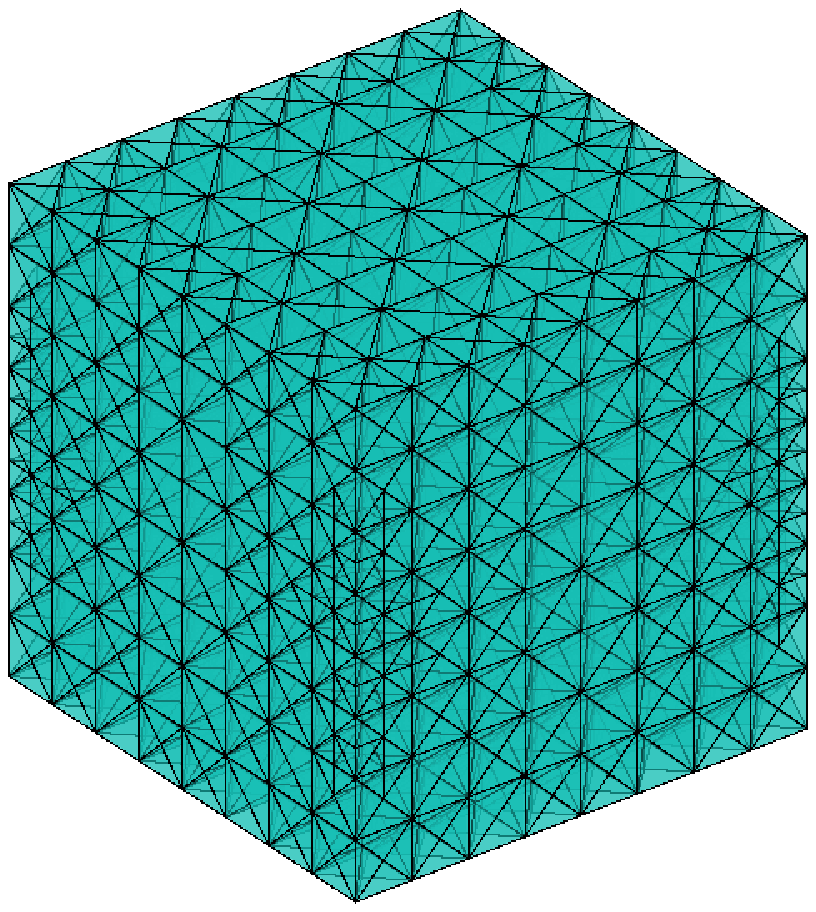}
\end{minipage}
\begin{minipage}[t]{0.5\textwidth}
  \centering
 \includegraphics[width=7.3cm]{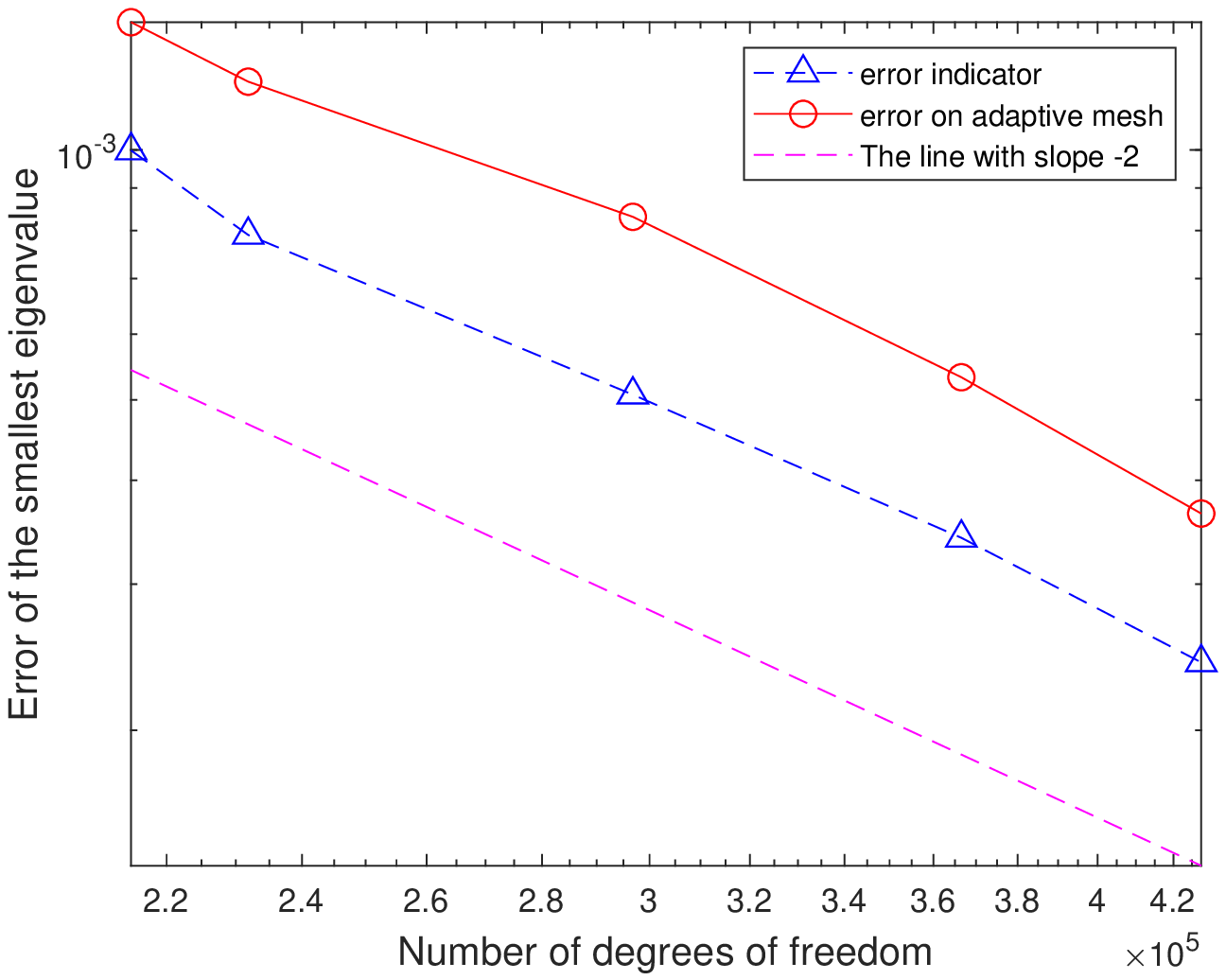}
\end{minipage}
\caption{Adaptive mesh after $l$=5 refinement times (left) and the error curves (right) of the smallest eigenvalue
by DGFEM using $\mathbb{P}_{3}-\mathbb{P}_{2}$ element on $\Omega_{2}$ }
\end{figure}
\begin{table}[H]
\caption{ \ The smallest eigenvalue on adaptive mesh on $\Omega_{1}$. }
\label{tab:1}
\begin{center}
\begin{tabular}{{cccccccccccccc}}
  \hline
        ~ & ~ & $\mathbb{P}_{1}-\mathbb{P}_{0}$ & ~ & $\mathbb{P}_{2}-\mathbb{P}_{1}$ & ~ & $\mathbb{P}_{3}-\mathbb{P}_{2}$ &   \\ \hline
        $l$ & $dof$ & $\lambda_{1,h_{l}}$  & $dof$ & $\lambda_{1,h_{l}}$  & $dof$ & $\lambda_{1,h_{l}}$   \\ \hline
        1 & 34944 & 72.17893    & 91392 & 71.39751    & 188160 & 70.94156     \\
        2 & 38532 & 72.31180    & 93432 & 71.40261    & 190680 & 70.94249     \\
        3 & 44668 & 72.37533   & 98532 & 71.38286   & 193480 & 70.94775     \\
        4 & 56108 & 72.30857    & 106148 & 71.36073    & 194740 & 70.96169     \\
        5 & 74074 & 72.22436    & 116416 & 71.33137    & 203280 & 70.96023     \\
        6 & 91416 & 72.08199    & 134844 & 71.30330    & 207060 & 70.96067     \\
        7 & 121212 & 71.85217    & 162248 & 71.26700    & 218820 & 70.96323     \\
        8 & 159224 & 71.73561   & 189244 & 71.24186    & 234500 & 70.96341     \\
        9 & 202878 & 71.64715    & 229568 & 71.20735    & 248500 & 70.96722     \\
        10 & 270738 & 71.56537    & 281928 & 71.17551    & 259840 & 70.96886     \\
        11 & 369278 & 71.49289    & 354212 & 71.14148    & 271740 & 70.96910     \\
        12 & 475852 & 71.44762    & 437512 & 71.11079    & 300860 & 70.97064    \\
        13 & 623792 & 71.40660    & 519588 & 71.08194    & 412653 & 70.98560 \\
        14 & 822172 & 71.35097   & 685780 & 71.05522   & - & -  \\
\hline
\end{tabular}\end{center}
\end{table}
\begin{table}[H]
\caption{ \ The smallest eigenvalue on adaptive mesh on $\Omega_{2}$. }
\label{tab:1}
\begin{center}
\begin{tabular}{{cccccccccccccc}}
  \hline
        ~ & ~ & $\mathbb{P}_{1}-\mathbb{P}_{0}$ & ~ & $\mathbb{P}_{2}-\mathbb{P}_{1}$  & ~ & $\mathbb{P}_{3}-\mathbb{P}_{2}$ &   \\ \hline
        l & $dof$ & $\lambda_{1,h_{l}}$  & $dof$ & $\lambda_{1,h_{l}}$  & $dof$ & $\lambda_{1,h_{l}}$   \\ \hline
        1 & 39936 & 63.68761    & 104448 & 62.28426    & 215040 & 62.17483     \\
        2 & 47892 & 63.59785    & 128520 & 62.26349    & 231840 & 62.17461    \\
        3 & 56732 & 63.51520    & 160820 & 62.24599    & 296800 & 62.17423     \\
        4 & 70044 & 63.38799    & 207944 & 62.22633    & 366520 & 62.17393     \\
        5 & 86632 & 63.26514    & 267920 & 62.21249    & 427560 & 62.17376     \\
        6 & 109252 & 63.09579    & 352648 & 62.19957    & 491060 & 62.17341  \\
        7 & 141440 & 62.91980    & 448664 & 62.18869    & - & -   \\
        8 & 186004 & 62.80083    & 611864 & 62.18345    & - & -  \\
        9 & 245440 & 62.71246    & - & - & - & -  \\
        10 & 339248 & 62.63796    & - & - & - & -   \\
        11 & 455754 & 62.56792    & - & - & - & -   \\
        12 & 608998 & 62.51093    & - & - & - & -   \\
\hline
\end{tabular}\end{center}
\end{table}

\section*{Acknowledgements}

\indent This work was supported by the National Natural Science Foundation of China (Grant Nos. 11561014,  11761022)
and the Science and Technology Planning Project of Guizhou Province (Guizhou Kehe Talent Platform [2017] No.5726).\\
\indent The authors sincerely thank Professor Jiayu Han of Guizhou Normal University for guiding the numerical experiments. \\


\begin{thebibliography}{99}
\bibitem{babuska1978}I. Babuska, W. C. Rheinboldt, Error estimates for adaptive finite element computations, SIAM .J. Numer. Anal. vol. 15 (1978)pp. 736-754.
\bibitem{Verfurth1996}R. Verf$\ddot{u}$rth, A Posteriori Error Estimation Techniques, Oxford University Press, New York, 1996.
\bibitem{dorfler1996}W. Dorfler, A convergent adaptive algorithm for Poisson's equation, SIAM J. Numer. Anal. vol. 33 (1996) pp. 1106-1124.
\bibitem{morin2002}P. Morin, R. H. Nochetto, K. Siebert, Convergence of adaptive finite element methods, SIAM Rev. vol. 44 (2002)pp. 631-658.
\bibitem{oden2011}M. Ainsworth, J. T. Oden, A Posteriori Error Estimation in Finite Element Analysis, Wiley-Interscience, New York, 2011.
\bibitem{brenner2012}S.C. Brenner, $C^{0}$ interior penalty methods, In Frontiers in Numerical Analysis-Durham 2010,
Lecture Notes in Computational Science and Engineering, Springer-Verlag, vol. 85 (2012)pp. 79-147.
\bibitem{Shi2013}Z. Shi, M. Wang, Finite Element Methods, Scientific Publishers, Beijing, 2013.
\bibitem{sunzhou2016}J. Sun, A. Zhou, Finite Element Methods for Eigenvalue Problems, CRC Press, Taylor Francis
Group, Boca Raton, London, New York, 2016.
\bibitem{Lovadina2009}C. Lovadina, M. Lyly, Stenberg, R.: A posteriori estimates for the Stokes eigenvalue problem. Numer.
Methods Partial Differ. Equ. vol. 25 (1) (2009)pp. 244-257.
\bibitem{Han2015}J. Han, Z. Zhang, Y. Yang, A new adaptive mixed finite element method based on residual type a
posterior error estimates for the Stokes eigenvalue problem. Numer. Methods Partial Differ. Equ. vol. 31 (1) (2015)pp. 31-53.
\bibitem{Armentano2014}M. G. Armentano, V. Moreno, A posteriori error estimates of stabilized low-order mixed
finite elements for the Stokes eigenvalue problem, J. Comput.  Appl. Math. vol. 269 (2014)pp. 132-149.
\bibitem{Liu2013} H. Liu, W. Gong, S. Wang, N. Yan, Superconvergence and a posteriori error estimates for the Stokes
eigenvalue problems. BIT. vol. 53 (3) (2013)pp. 665-687.
\bibitem{Xie2014}S. Jia, F. Lu, H. Xie, A Posterior error analysis for the nonconforming discretization of Stokes eigenvalue problem. Acta
Mathematica Sinica (English Series), 2014.
\bibitem{Gedicke2018}J. Gedicke, A. Khan, Arnold-Winther mixed finite elements for Stokes eigenvalue problems. SIAM
J. Sci. Comput. vol. 40 (5) (2018)pp. A3449-A3469.

\bibitem{onder2016}$\ddot{O}$nder, T$\ddot{u}$rk, Daniele, Ramon Codina. A stabilized finite element method for the two-field and three-field Stokes eigenvalue problems. Comput. Methods Appl. Mech. Engrg. 310 (2016) 886-905
\bibitem{reed1973}W. H. Reed and T. R. Hill, Triangular mesh methods for the neutron transport equation,
Technical Report LA-UR-73-479, Los Alamos Scientifik Laboratory, 1973.

\bibitem{Cockburn1999}B.Cockburn, G.E.Karniadakis, C. Shu, Discontinuous Galerkin Methods,Thoery, Computation and Applications, Springer-Verlag, 1999.

\bibitem{Hesthaven2008}Jan S. Hesthaven, Tim Warburton, Nodal Discontinuous Galerkin Methods, Algorithms, Analysis, and Applications. Springer-Verlag, New York, 2008.
\bibitem{Riviere2008}B. Rivi$\grave{e}$re, Discontinuous Galerkin Methods for Solving Elliptic and Parabolic Equations. Theory and Implementation, SIAM. 2008.
\bibitem{Pietro2009}D.A. Di Pietro and A. Ern, Discrete functional analysis tools for discontinuous Galerkin methods with application to the incompressible Navier-Stokes equations, Math. Camp. vol 79(271) (2010) pp.1303-1330.
\bibitem{Cangianl2010}A.Cangianl, Z.Dong, E.H.Georgoulis, P.Houston, hp-Version Discontinuous Galerkin Method on Polygonal and Polyhedral Meshes, Springer, 2010.
\bibitem{Antonio2012}D. Antonio, E. Alexandre, Mathematical Aspects of  Discntinuous Galerkin Methods, Springer-Verlag, 2012.

\bibitem{Ern2021} A. Ern, L J Ond Rm Gue, Finite Elements IIGalerkin approximation, elliptic and mixed PDEs.  2021.

\bibitem{Brezzi2000}F. Brezzi, G. Manzini, D. Marini, P. Pietra, A. Russo, Discontinuous Galerkin approximations for
elliptic problems. Numer. Methods Partial Differ. Equ. 16 (4) (2000)pp. 365-378.




\bibitem{Antonietti2006}P. F. Antonietti, A. Buffa, I. Perugia, Discontinuous Galerkin approximation of the Laplace eigenproblem. Comput. Meth. Appl. Mech.  Engrg. vol. 195 (25/28) (2006)pp. 3483-3503.
\bibitem{Zeng2017}Y. Zeng, F. Wang, A posteriori error estimates for a discontinuous Galerkin approximation of Steklov eigenvalue problems. Appl. Math. vol. 62 (3) (2017)pp. 243-267.
\bibitem{Brenner2015}S. C. Brenner, P. Monk, J. Sun, $C^{0}$ interior penalty Galerkin method for biharmonic eigenvalue problems.
 Spectral and High Order Methods for Partial Differential Equations, Lect. Notes Comput. Sci. Eng. 106 (2015) pp.3-15.
\bibitem{Wang2019}L. Wang, C. Xiong, H. Wu, F. Luo, A priori and a posteriori analysis for discontinuous Galerkin finite element approximations of biharmonic eigenvalue
problems, Adv. Comput. Math. 45(5-6) (2019) pp. 2623-2646.
\bibitem{Geng2016}H. Geng, X. Ji, J. Sun, L. Xu, $C^{0}$IP methods for the transmission eigenvalue problem, J. Sci, Comput. 68(2016) pp.326-338.
\bibitem{Yang2017}Y. Yang, H. Bi, H. Li, J. Han, A $C^{0}$IPG method and its error estimates for the Helmholtz transmission eigenvalue problem, J. Comput. Appl. Math. vol. 326
(2017)pp. 71-86.
\bibitem{Buffa2006}A. Buffa, I. Perugia, Discontinuous Galerkin approximation of the Maxwell eigenproblem, SIAM J. Numer. Anal. vol. 44(5) (2006)pp. 2198-2226.
\bibitem{Buffa2007}A. Buffa, P. Houston, I. Perugia, Discontinuous Galerkin computation of the Maxwell eigenvalues on simplicial meshes, J. Comput. Appl. Math.
vol. 204 (2007)pp. 317-333.
\bibitem{Gedicke2019}  J. Gedicke , A. Khan , Divergence-conforming discontinuous Galerkin finite elements for Stokes eigenvalue problems, Numer. Math. 144 (3) (2019)pp.
585-611.
\bibitem{LepeF2020}F. Lepea,  Mora D . Symmetric and Nonsymmetric Discontinuous Galerkin Methods for a Pseudostress Formulation of the Stokes Spectral Problem. SIAM Journal on Scientific Computing, 42 (2) (2020)pp. A698-A722.

\bibitem{Badia2014}S. Badia, et al., Error analysis of discontinuous Galerkin methods for the Stokes problem under minimal regularity. IMA. J. Numer. Anal. vol. 34(2014)pp. 800-819.

\bibitem{Hansbo2002}P. Hansbo and M. G. Larson, Discontinuous Galerkin methods for incompressible and
nearly incompressible elasticity by Nitsche's method, Comput. Methods Appl. Mech.
Engrg. vol. 191 (2002)pp. 1895-1908.

\bibitem{Girault2005}V. Girault, B. Rivie`re, M.F. Wheeler, A discontinuous Galerkin method with non-overlapping domain decomposition for the
Stokes and Navier-Stokes problems, Math. Comput. vol. 74 (2005)pp. 53-84.

\bibitem{Riviere2006}B. Rivi$\grave{e}$re and V. Girault, Discontinuous finite element methods for incompressible
flows on subdomains with non-matching interfaces, Computer Methods in Applied Mechanics and Engineering, 195 (25/28) (2006) 3274-3292.



\bibitem {toselli2002}A. Toselli, hp discontinuous Galerkin approximations for the Stokes problem, Math. Models Methods
Appl. Sci. vol. 12 (11) (2002)pp. 1565-1597



\bibitem{Houston2005}P. Houston, D. Sch$\ddot{o}$tzau, T.P. Wihler, Energy norm a posteriori error estimation for mixed discontinuous
Galerkin approximations of the Stokes problem. J. Sci. Comput. vol. 22 (23) (2005)pp. 347-370.

\bibitem{Ohannes2003}O. A. Karakashian, F. Pascal, A posteriori error estimates for a discontinuous Galerkin approximation of second-order elliptic problems. SIAM J. Numer. Anal. vol. 41(6)(2004) pp. 2374-2399.

\bibitem{Brenner2003}S.C. Brenner, Poincar\'e-Friedrichs inequalities for piecewise $H^{1}$ functions. SIAM J. Numer. Anal. vol. 41 (2003)pp. 306-324.

 \bibitem{Schotzau2002}D. Sch$\ddot{o}$tzau, C. Schwab, A. Toselli, Mixed hp-DGFEM for incompressible flows. SIAM J. Numer.
Anal. vol. 40 (6) (2002)pp. 2171-2194.


\bibitem{Perugia2002}I. Perugia, D. Sch$\ddot{o}$tzau, The $hp$-local discontinuous Galerkin method for low-frequency time-harmonic Maxwell equations. Math. Comp. vol. 72 (2002)pp. 1179-1214.

\bibitem{Houston2007}P. Houston, I. Perugia, D. Sch$\ddot{o}$tzau. An a posteriori error indicator for discontinuous Galerkin discretizations of H(curl)-elliptic partial differential equations. IMA J. Numer. Anal. vol. 27 (2007)pp. 122-150.

\bibitem{Babuska1991book}I. Babuska, J.E. Osborn, Eigenvalue problems, in: P.G. Ciarlet, J.L. Lions (Eds.), Finite Element Methods (Part I), in: Handbook of Numerical Analysis, vol.2, North-Holland: Elsevier Science Publishers, (1991)pp. 641-787.

\bibitem{D.Boffi2010}D. Boffi, Finite element approximation of eigenvalue problems, Acta Numer. vol. 19 (2010)pp. 1-120.

\bibitem{yzl2010}Y.Yang, Z. Zhang, F. Lin: Eigenvalue approximation from below using non-conforming
finite elements. Sci. China, Math. vol. 53 (2010)pp. 137-150.

\bibitem{ScottZhang1990} L. R. Scott, S. Zhang, Finite element interpolation of non-smooth functions satisfying boundary conditions, Math. Comp. vol. 54 (1990)pp. 483-493.



\bibitem{Kanschat2008}G. Kanschat, D. Sch$\ddot{o}$tzau. Energy norm a posteriori error estimation for divergence-free discontinuous Galerkin approximations of the Navier-Stokes equations. Int. J. Numer. Meth. in Fluids, 2008.

\bibitem{Wu2007}H. Wu, Z. Zhang, Can we have superconvergent gradient recovery under adaptive meshes, SIAM J. Numer. Anal. vol. 45 (2007)pp. 1701-1722.

\bibitem{Yang2020}Y. Yang, Y. Zhang, H. Bi, A type of adaptive $C^{0}$ non-conforming finite element method
for the Helmholtz transmission eigenvalue problem, Comput. Methods Appl. Mech. Engrg,
360 (2020), Doi: 10.1016/j.cma.2019.112697.

\bibitem{Dai2008}X. Dai, J. Xu, A. Zhou, Convergence and optimal complexity of adaptive finite element eigenvalue computations, Numer. Math. vol. 110 (2008)pp. 313-355.

\bibitem{Chen2009}L. Chen, iFEM: An integrated finite element method package in matlab, Technical Report, University of California at Irvine, 2009.

\bibitem{Sun2022}L. Sun, Y. Yang, The a posteriori error estimates and adaptive computation of nonconforming mixed finite elements for the Stokes eigenvalue problem, Applied Mathematics and Computation, 421 (2022) 126951.

\bibitem{LepeF2021}F. Lepe, G. Rivera, A virtual element approximation for the pseudostress formulation of
the Stokes eigenvalue problem, Comput. Methods Appl. Mech. Engrg. 379 (2021) 113753.



\bibitem{girault1979}V. Girault, P.A. Raviart, Finite Element Approximation of the Navier-Stockes Equations. Springer-Verlag, 1979.

\bibitem{Riviere2006+}B. Rivi$\grave{e}$re and V. Girault, Discontinuous finite element methods for incompressible
flows on subdomains with non-matching interfaces, Comput. Meth. Appl.
Mech. Engrg. vol. 195(2006)pp. 3274-3292.








\end{thebibliography}

\end{document}